\documentclass[12pt,a4paper]{book}
\usepackage{amssymb,amsfonts,amsmath,amscd}
\usepackage[all]{xy}

\usepackage{cite}
\usepackage[latin1]{inputenc}
\usepackage{amsmath, amssymb, amsthm}
\usepackage{graphics, booktabs, ctable, array}
\usepackage{colortbl, xcolor}
\usepackage{color}
\usepackage{booktabs,ctable,array}
\usepackage{fancyhdr}
\usepackage{atbeginend}
\usepackage{algorithmic}

\definecolor{DarkBlue}{rgb}{0.1,0.05,0.7}
\definecolor{CaptionColor}{rgb}{0.6,0.25,0.1}
\usepackage[bookmarks=true,pdfborder={0 0 0},colorlinks=true,urlcolor=DarkBlue,linkcolor=DarkBlue,citecolor=DarkBlue]{hyperref}

\def\sqr#1#2{{\vcenter{\hrule height.#2pt \hbox{\vrule width.#2pt
height#1pt \kern#1pt \vrule width.#2pt} \hrule height.#2pt}}}
\def\square{\mathchoice\sqr64\sqr64\sqr{4}3\sqr{3}3}
\def\qed{\hfill$\square$}
\def\demo{\noindent{\bf Proof: }}

\AfterBegin{itemize}{\setlength{\itemsep}{-1.5mm}}
\AfterBegin{enumerate}{\setlength{\itemsep}{-1.5mm}}

\newtheorem{Theorem}{\bf Theorem}[section]
\newtheorem{Proposition}[Theorem]{\bf Proposition}
\newtheorem{Lemma}[Theorem]{\bf Lemma}
\newtheorem{Corollary}[Theorem]{\bf Corollary}

\newtheorem{Remark}[Theorem]{\bf Remark}
\newtheorem{Example}[Theorem]{\bf Example}
\newtheorem{Definition}[Theorem]{\bf Definition}

\newtheorem{Discussion}[Theorem]{\bf Discussion}
\newtheorem{Algorithm}[Theorem]{\bf Algorithm}
\newtheorem{Numbered}[Theorem]{}
\newtheorem*{Theorem A}{\bf Theorem A}
\newtheorem*{Theorem A*}{\bf Theorem A*}
\newtheorem*{summary:exemple-classic}{\bf Example \ref{exemple-classic}}
\newtheorem*{Theorem B}{\bf Theorem B}
\newtheorem*{summary:exemple-article1}{\bf Example \ref{exemple}}

\newcommand{\fm}{\mbox{$\mathfrak{m}$}}

\newcommand{\fp}{\mbox{$\mathfrak{p}$}}
\newcommand{\fq}{\mbox{$\mathfrak{q}$}}

\newcommand{\cm}{\mbox{$\mathcal{M}$}}

\newcommand{\lsigma}{\mbox{$\lambda_\sigma$}}
\newcommand{\xsigma}{\mbox{${x}^\sigma$}}
\newcommand{\Xsigma}{\mbox{${X}^\sigma$}}

\newcommand{\kdim}[1]{\mbox{${\rm
dim}\hspace{0.02cm}_k\hspace{0.05cm}{(#1)}$}}

\newcommand{\rees}{\mbox{$\mathbf{R}$}}
\newcommand{\graded}{\mbox{$\mathbf{G}$}}
\newcommand{\fiber}{\mbox{$\mathbf{F}$}}
\newcommand{\symmetric}{\mbox{$\mathbf{S}$}}

\newcommand{\reltype}{\mbox{$\rm{rt}$}}
\newcommand{\rednumber}{\mbox{$\rm{r}$}}
\newcommand{\spread}{\mbox{$\ell$}}
\newcommand{\height}{\mbox{$\rm{ht}$}}
\newcommand{\grade}{\mbox{$\rm{grade}$}}

\newcommand{\lcm}{\mbox{$\rm{lcm}$}}
\newcommand{\adeviation}{\mbox{$\rm{ad}$}}
\newcommand{\sadeviation}{\mbox{$\rm{sd}$}}

\newcommand{\koszul}[3]{\mbox{$H_{#1}(#2\,;#3)$}}

\newcommand{\cycle}[3]{\mbox{$Z_{#1}(#2\,;#3)$}}
\newcommand{\boundary}[3]{\mbox{$B_{#1}(#2\,;#3)$}}

\begin{document}

\begin{titlepage}
\begin{center}

{\large UNIVERSITAT POLITÈCNICA DE CATALUNYA} \\[3.5cm]
 
{\bf \large THE EQUATIONS OF REES ALGEBRAS OF IDEALS OF ALMOST-LINEAR TYPE} \\[1cm]
by \\[1cm]

{\Large Ferran Muiños} \\[11.0cm]

\end{center}

\normalsize

\begin{center}
Programa de Doctorat de Matemàtica Aplicada\\
Universitat Politècnica de Catalunya (UPC)
\end{center}
\end{titlepage}

\newpage
\thispagestyle{empty}
Advisor: Francesc Planas-Vilanova;

Document version: 19 November 2011. Corrections on the PhD Thesis with same title, presented and approved on 3 October 2011, Facultat de Matemàtiques i Estadística (FME), Universitat Politècnica de Catalunya (UPC). 

\newpage

\setcounter{page}{1}
\vspace*{12cm}

\begin{quote}
Certainly the best times were when I was alone with the mathematics, free of ambition and of pretence, and indifferent to the world.
\end{quote}
\begin{flushright}
Robert P. Langlands\\
\emph{Mathematical retrospections}
\end{flushright}
\vspace*{1cm}
\begin{quote}
Algebra is the offer made by the devil to the mathematician. The devil says: I will give you this powerful machine, it will answer any question you like. All you need to do is give me your soul: give up geometry and you will have this marvelous machine. 
\end{quote}
\begin{flushright}
Michael F. Atiyah \\ 
From Barry Mazur's foreword to the 2005 reprint of\\ 
Tobias Dantzig's \emph{Number: The Language of Science}.
\end{flushright}

\chapter*{Acknowledgements}

Many people have supported this work in many different ways, making this journey easier, deeper and worthier. 

\medskip

\noindent
A few words cannot adequately express the thanks I owe to Francesc Planas for his constant support, valuable insights, critical comments, guidance and affection.

\medskip

\noindent
As regards the present work, I would also like to acknowledge the following contributors:

\medskip
José M. Giral, for carefully reading early versions of most of which is included in this work and also for his well-written and useful lecture notes; 

\medskip
Liam O'Carroll and Santiago Zarzuela, for carefully reading this work in its final form; 

\medskip
Miguel Ángel Barja, for hosting me as a PhD student from the very beginning; 

\medskip
All the people who contributed to the courses and seminars I attended as a PhD student, about a wealth of exciting topics; special thanks to Jaume Amorós, Miguel Ángel Barja, Xavier Cabré, Carles Casacuberta, Pere Pascual and Francesc Planas; 

\medskip
Colleagues and employers at Centre de Recerca i Innovació de Catalunya (CRIC) and Theia Innovation Consulting (THEIA), since they all set the conditions to make a professional path out of the academic bounds compatible with my mathematical interests; in particular, thanks to Eva Romagosa, for the mentoring and her share in promoting a healthy working environment; thanks also to Irene Larroy, for the first opportunity and long-standing friendship;

\medskip
Andreu Làzaro, for his friendship and the interesting discussions we held while digesting the pearls that we stumbled upon around every corner of Atiyah-MacDonald's book;

\medskip
Guillem Godoy, for his friendship and valuable personal contribution to my problem solving simulator, shaped through a wealth of energyzing conversations and sportive adventures of all sorts.

\bigskip

\bigskip

\noindent
From a broader perspective, personal thanks to my loved ones, family, friends, ancestors and people I admire. After all, they hold the largest share as sources of energy to accomplish the most demanding personal projects.

\tableofcontents
\chapter*{Abstract}
\addcontentsline{toc}{chapter}{Abstract}

The Rees algebra $\rees(I)=R[It]\subset R[t]$ of an ideal $I\subset R$ of a Noetherian local ring $R$ plays a major role in commutative algebra and in algebraic geometry, since ${\rm Proj}(\rees(I))$ is the blowup of the affine scheme ${\rm Spec}(R)$ along the closed subscheme ${\rm Spec}(R/I)$.

So far, the problem of describing the equations of Rees algebras of ideals, as well as other related algebras, has shown to be relevant in order to further understand these major algebraic objects. The equations of $\rees(I)$ arise as the elements in the kernel of a symmetric presentation $\varphi:V\to\rees(I)$. While this kernel may differ from one presentation to another, the degrees of a minimal generating set of homogeneous elements are known to be independent of $\varphi$. The top degree among such generating sets, known as the \emph{relation type} and denoted by $\reltype(I)$, is a coarse measurement of the complexity of the underlying Rees algebra which is nonetheless a useful numerical invariant. The ideals $I$ such that $\reltype(I)=1$, known as ideals of \emph{linear type}, have been intensely studied so far.

In this dissertation, we tackle the problem of describing the equations of $\rees(I)$ for $I=(J,y)$, with $y\notin J$ and $J$ being of linear type, i.e., for ideals of linear type up to one minimal generator. Throughout, such ideals will be referred to as ideals of \emph{almost-linear type}.

The main results of this work stem from two different approaches towards the problem. 

In Theorem~A, we give a full description of the equations of Rees algebras of ideals of the form $I=(J,y)$, with $J=(x_1,\ldots,x_s)$ satisfying an homological vanishing condition. Theorem~A permits us to recover and extend well-known results about families of ideals fulfilling the almost-linear type condition due to Vasconelos, Huckaba, Trung, Heinzer and Kim, among others. 

Let $\alpha:\symmetric(I)\to\rees(I)$ be the canonical morphism from the symmetric algebra of $I$ to the Rees algebra of $I$. In Theorem~B, we prove that the injectivity of a single component of $\alpha:\symmetric(I)\to\rees(I)$ propagates downwards, provided $I$ is of almost-linear type. In particular, this result gives a partial answer to a question posed by Tchernev.

Finally, packs of examples are introduced, which illustrate the scope and applications of each of the results presented. The author also gives a collection of computations and examples which motivate ongoing and future research.
\chapter*{Notations}
\addcontentsline{toc}{chapter}{Notations}

\noindent
$R$: commutative ring with unit $1\neq 0$, most often Noetherian and local.

\noindent
$k$: field.

\noindent
$I$, $J$: finitely generated ideals of $R$, with $J\subset I$.

\noindent
$J_i=(x_1,\ldots,x_i)$.

\noindent
$M$: $R$-module.

\noindent
$(M_1:_R M_2)=\{a\in R \;|\; aM_2\subseteq M_1\}$.

\noindent
$(J:_R y^\infty)=\cup_{n\geq 1} (J:_R y^n)$: saturation.

\noindent
$M_1\otimes_R M_2$: tensor product of $R$-modules.

\noindent
$\fp$, $\fq$: prime ideals.

\noindent
$\fm$: maximal (or maximal homogeneous) ideal.

\noindent
$VV_J(I)_n=(J\cap I^n)/JI^{n-1}$: $n$-th Valabrega-Valla module of $I$ with respect to $J$.


\bigskip

\noindent
${\rm Ass}(R/I)$: associated prime ideals of $R/I$.

\noindent
${\rm Min}(R/I)$: minimal prime ideals of $R/I$.

\noindent
$\mathfrak{N}(R)$: nilradical of $R$.

\bigskip

\noindent
$U=\oplus_{n\geq 0} U_n$: standard $R$-algebra.

\noindent
$U_+=\oplus_{n>0} U_n$: irrelevant ideal of $U$.

\noindent
$V,W$: polynomial rings over $R$.

\noindent
$\varphi$, $\psi$: polynomial presentations.

\noindent
$Q=\oplus_{n\geq 0} Q_n$: homogeneous ideal.

\noindent
$Q\langle r\rangle$: subideal of $Q$ generated by its homogeneous 
elements of degree $\leq r$.

\bigskip

\noindent
$\partial$: differential operator of a complex of $R$-modules.

\noindent
$\koszul{i}{z}{R}$: $i$-th Koszul homology.

\noindent
$\cycle{i}{z}{R}$: $i$-th Koszul cycles.

\noindent
$\boundary{i}{z}{R}$: $i$-th Koszul boundaries.

\bigskip

\noindent
$\symmetric(M)=\oplus_{n\geq 0}\symmetric_n(M)$: symmetric algebra of $M$.

\noindent
$\rees(I)=\oplus_{n\geq 0} I^n$: Rees algebra of $I$.

\noindent
$\graded(I)=\oplus_{n\geq 0} I^n/I^{n+1}$: associated graded ring of $I$.

\noindent
$\fiber(I)=\oplus_{n\geq 0} I^n/\fm I^n$: fiber cone of $I$ with respect to $\fm$.

\bigskip

\noindent
$\alpha_I$, $\beta_I$, $\gamma_I$ (or simply denoted $\alpha$, $\beta$ and $\gamma$) are the 
canonical blowing-up morphisms:

\noindent 
$\alpha_I:\symmetric(I)\to\rees(I)$; 

\noindent
$\beta_I:\symmetric(I/I^2)\to\graded(I)$;

\noindent
$\gamma_I:\symmetric(I/\fm I)\to\fiber(I)$.

\noindent
$E(I)_n=\ker\alpha_n/\symmetric_1(I)\cdot\ker\alpha_{n-1}$: $n$-th module of effective relations.

\bigskip

\noindent
$X^\sigma, x^\sigma$: monomials in multi-index notation.

\noindent
${\rm supp}(\sigma)$: support of an integer vector.

\bigskip

\noindent
$\dim(R)$: Krull dimension.

\noindent
$\lambda(M)$: length of an $R$-module $M$.

\noindent
$\mu(I)$: minimal number of generators of $I$.

\noindent
$\height(I)$: height of $I$.

\noindent
${\rm pd}_R(M)$: projective dimension of an $R$-module $M$.

\noindent
$\reltype(I)$: relation type of $I$.

\noindent
$\rednumber_J(I)$: reduction number with respect to $J$.

\noindent
$\spread(I)$: analytic spread of $I$. 

\noindent
$\mu(I)-\height(I)$: deviation of $I$.

\noindent
${\rm ad}(I)=\spread(I)-\height(I)$: analytic deviation of $I$.

\noindent
${\rm sd}(I)=\mu(I)-\spread(I)$: second analytic deviation of $I$.

\onecolumn
\chapter*{Summary of results}\label{summary}
\addcontentsline{toc}{chapter}{Summary of results}

For the sake of clarity, let us first introduce the basic assumptions 
and state the main results presented along this work. 

\subsection*{The equations of Rees algebras of ideals of almost-linear type}

Recall that the reduction number of $I$ with respect to $J$, $J\subseteq I$, 
is the least integer $r\geq 0$ such that $I^{r+1}=JI^r$, denoted by 
$\rednumber_J(I)$. Set $V=R[X_1,\ldots,X_s,Y]$ and let $\varphi:V\to\rees(I)$ 
be the polynomial presentation of $\rees(I)$, the Rees algebra of $I$, 
sending $X_i$ to $x_it$ and $Y$ to $yt$. Let $Q=\oplus_{n\geq 1}Q_n$ 
be the kernel of $\varphi$. Given an integer $m\geq 0$, set 
$Q\langle m\rangle\subset Q$ the ideal generated by the homogeneous 
elements of $Q$ of degree at most $m$ in $X_1,\ldots,X_s,Y$. The 
relation type of $I$, denoted by $\reltype(I)$, is the least integer 
$N\geq 1$ such that $Q=Q\langle N\rangle$. Let $\graded(I)=\oplus_{n\geq 0} I^n/I^{n+1}$ 
be the associated graded ring of $I$. If $z\in R\,\backslash \cap_{n\geq 0} I^n$, we will denote by $z^*$
the initial form of $z$ in $\graded(I)$, i.e., $z^*=z+I^{n+1}$, 
where $z\in I^n\backslash I^{n+1}$.

\begin{Theorem A}
Let $(R,\fm)$ be a Noetherian local ring and let $I$ be an ideal of
$R$. Let $x_1,\ldots,x_s,y$ be a minimal generating set of $I$, where
$J=(x_1,\ldots,x_s)$ is a reduction of $I$ with reduction number
$r=\rednumber_J(I)$. Assume that $x_1,\ldots,x_s$ verify 
the following condition for all $n\geq 2$:
\[ 
((x_1,\ldots,x_{i-1})I^{n-1}: x_i)\cap I^{n-1}=(x_1,\ldots,x_{i-1})I^{n-2}, 
\textrm{ for all } i=1,\ldots,s. \label{eq:trung-condition} \tag{$\mathcal{T}_n$}
\]
Then, for each $n\geq 2$, the map sending $F\in Q_{n}$ to
$F(0,\ldots,0,1)\in (JI^{n-1}: y^{n})$ induces an isomorphism of
$R$-modules
\begin{eqnarray*}
\left[ \frac{Q}{Q\langle n-1\rangle}\right]_{n}\cong
\frac{(JI^{n-1}: y^{n})}{(JI^{n-2}: y^{n-1})}.
\end{eqnarray*}
In particular, $\reltype(I)=\rednumber_J(I)+1$ and there is a form
$Y^{r+1}-\sum X_iF_i\in Q_{r+1}$, with $F_i\in V_r$, such that
$Q=(Y^{r+1}-\sum X_iF_i)+Q\langle r\rangle$. Moreover, if 
$x_1,\ldots,x_s$ is an $R$-sequence and $x_1^*,\ldots,x_{s-1}^*$ is
a $\graded(I)$-sequence, then $x_1,\ldots,x_s$ verify condition $(\mathcal{T}_n)$
for all $n\geq 2$.
\end{Theorem A}

Theorem~A permits us to recover, even extend, a set of well known results 
regarding classes of almost-linear type ideals due to 
Vasconelos \cite[Theorem~2.3.3]{vasconcelos1}, 
Trung \cite[Theorem~1.2, Proposition~5.1]{trung1}, 
Huckaba \cite[Theorems~1.4, 1.5]{huckaba2} 
and Heinzer-Kim \cite[Theorem~5.6]{hk}.

\subsection*{Examples}

Let $(R,\fm)$ be a Noetherian local ring. Let $a,b\in R$ be an $R$-sequence. 
As an application of Theorem~A, we give a full description of the 
equations of $\rees(I)$, where $I$ is any of the following ideals:
\begin{itemize}
\item[$(i)$] $I=(a^p,b^p,ab^{p-1})$ with $p\geq 2$ an integer (see Example \ref{exemple-classic});
\item[$(ii)$] $I=(a^p,b^p,a^2b^{p-2})$ with $p\geq 2$ an odd integer (see Example \ref{exemple-pseudo-classic}).
\end{itemize}
These examples may be considered folklore and they are easily computable, 
although this computation may be expensive, for a given $p\geq 2$. However, 
our approach is new. 

\subsection*{The injectivity of the canonical blowing-up morphism}

Let $\alpha_{I}:\symmetric(I)\rightarrow \rees(I)$ be the canonical
morphism from the symmetric algebra of $I$ to the Rees algebra of
$I$. We write $\alpha_{I,p}$ the $p$-th graded component of $\alpha_I$.
Fixing $p\geq 2$, Tchernev asked in \cite{tchernev} whether 
$\alpha_{I,p}$ being an isomorphism implies that $\alpha_{I,n}$ is an 
isomorphism for each $2\leq n\leq p$. Observe that, 
as regards the following result, we do not need the 
Noetherian hypothesis.

\begin{Theorem B}
Let $R$ be a commutative ring. Let $I=(x_1,\ldots ,x_s,y)$ be an ideal
of $R$ and let $p\geq 2$ be an integer. Suppose that the ideal 
$J=(x_1,\ldots ,x_s)$ verifies that $\alpha_{J,n} :\symmetric_n(J)\rightarrow J^n$ 
is an isomorphism for all $2\leq n\leq p$. Then the following conditions are 
equivalent:
\begin{enumerate}
\item[$(i)$] $\alpha_{I,p}:\symmetric_p(I)\to I^p$ is an isomorphism;
\item[$(ii)$] $\alpha_{I,n} :\symmetric_n(I)\rightarrow I^n$ is an
isomorphism for each $2\leq n\leq p$.
\end{enumerate}
\end{Theorem B}

Remark that the ideals of almost-linear type fulfil the hypotheses of Theorem~B, 
hence the equivalence of $(i)$ and $(ii)$ holds true for such ideals.

\subsection*{Example}

K\"uhl \cite[Example~1.4]{kuhl} gave an example 
of a finitely generated ideal $I$ with $\alpha_{I,n}$ being an isomorphism 
for $n\geq 3$, while $\ker\alpha_2\neq 0$. In particular, 
$\widetilde{\alpha}_{I}:{\rm Proj}(\rees(I))\to {\rm Proj}(\symmetric(I))$ 
is an isomorphism of schemes, but $\alpha_{I}$ is not an isomorphism 
of $R$-algebras. The following example shows, for a given integer $p\geq 2$, 
an ideal $I$ such that $\alpha_{I,n}$ is an 
isomorphism for $n\geq p+1$, whereas $\alpha_{I,p}$ is not.

\begin{summary:exemple-article1}{\rm
Let $k$ be a field and let $p\geq
2$. Let $S=k[U_0,\ldots,U_p,X,Y]$ be a polynomial ring and let $Q$ be
the ideal of $S$ defined as $Q=Q_1 + (U_0X^p)$, where
\begin{eqnarray*}
Q_1=(U_0Y,U_0X-U_1Y,U_1X-U_2Y,\ldots, U_{p-1}X-U_pY, U_pX).
\end{eqnarray*}
Let $R$ be the factor ring $S/Q=k[u_0,\ldots,u_p,x,y]$ and consider
the ideal $I=(x,y)\subset R$. Then, $\alpha_{I,n}$ is an isomorphism
for all $n\geq p+1$, whereas $\alpha_{I,p}$ is not.
}\end{summary:exemple-article1}
\chapter{Introduction}\label{intro}

\section{Rees algebras and blowing-up}

So far, several graded rings associated with a (Noetherian local) 
ring $R$ and an ideal $I$ have been given much attention within
the bounds of commutative algebra and algebraic geometry: 
the Rees algebra $\rees(I)=R[It]=\oplus_{n\geq 0} I^{n}$, 
the associated graded ring $\graded(I)=\oplus_{n\geq 0} I^n/I^{n+1}$, 
the fiber cone $\fiber(I)=\oplus_{n\geq 0}I^n/\fm I^n$, 
the extended Rees algebra $R[It,t^{-1}]$, the symbolic Rees 
algebra $\rees_s(I)=\oplus_{n\geq 0} I^{(n)}$, as well as
approximations to these rings such as the symmetric algebras 
$\symmetric_R(I)$ and $\symmetric_{R/I}(I/I^2)$. All these 
graded rings are best known as blowing-up algebras. 
Today, the Rees algebra plays a central role among them, 
although originally more attention was focused on the associated 
graded ring and the extended Rees algebra. 

\medskip

Rees algebras were originally introduced by David Rees to study several 
problems in commutative algebra, thus yielding, according to \cite{hh2}, 
to three most remarkable results that have become classical over time: 
the construction of a counterexample to Zariski's generalisation of 
Hilbert's $14$-th problem \cite{rees2}, the Artin-Rees lemma \cite{rees1} 
and a celebrated criterion for a local ring to be analytically unramified 
\cite{rees3}. Rees algebras of ideals in Noetherian rings capture in 
another Noetherian ring a great deal of information about how the ideal 
$I$ sits in the ring $R$ and how its powers change. The study of Rees algebras 
is important to several fields of study and open problems in algebra as well 
as in geometry.

As it is claimed in \cite[Introduction]{hh2} ever since the introduction 
of the Hilbert function and the discovery of its key numerical invariants, 
the importance of understanding how the ideals behave as we raise them to 
powers is widely acknowledged. The study of Cohen-Macaulayness, normality 
or torsion-freeness of blowing-up algebras has played a major role so far, 
and such topics have been often studied in terms of numerical invariants.

\medskip

From an algebraic geometry perspective, given the affine scheme $X={\rm Spec}(R)$ 
and a closed subscheme $Z={\rm Spec}(R/I)\subset X$, the morphism of schemes 
$Y={\rm Proj}(\rees(I))\to {\rm Spec}(R)=X$ is the construction known as
\emph{blowing-up of} $X$ \emph{along} $Z$. The 
procedure of blowing-up is a standard, explicit way of modifying a variety 
(see \cite[I.4, II.7]{hartshorne1}). In the review \cite{abramovich} the 
blowing-up is presented as a surgery operation in which a 
subvariety $Z\subset X$ is removed and replaced by the sets of normal 
directions of $Z$ in $X$. When $X$ and $Y$ are smooth, $Z$ is replaced by 
$\mathbb{P}(N_{Z/X})$ the projectivization of the normal bundle. In 
particular, if $Z$ is a point, $Z$ is replaced by $\mathbb{P}^{\dim X - 1}$.

Satisfactory answers to open questions related to the birrational 
geometry of varieties have been given using techniques relying on 
blowing-up. As a first example, recall that any birational map 
between algebraic surfaces factors as a sequence of blowing-ups 
modulo isomorphism \cite[Chapter II]{beauville} (see also 
\cite[Remark V.5.6.1]{hartshorne1} for generalisations in higher 
dimensions). Blowing-up plays also a central role in the problem of 
\emph{resolution of singularities}, which can be stated as follows: 
given a variety $X$, find a smooth variety $Y$ birational to $X$ 
(the birational map $\pi:Y\to X$ being known as a resolution of singularities). 
While analysing nonsingular varieties in characteristic zero one can use, 
under mild conditions, the tools and techniques of complex manifold theory. 
In characteristic zero a resolution of singularities can be fully achieved 
as a finite composition of blowing-ups. While the existence of a resolution of 
singularities is not known in arbitrary characteristic, it has already 
been shown for curves, surfaces and $3$-folds (see, for instance, 
\cite{hauser} and \cite{kollar}).

A subject recently introduced in the context of computer geometric 
modelling are the methods of implicitization by moving curves and surfaces 
(see e.g. \cite{csch}). Such methods essentially provide with a description 
of the image of a rational map between projective spaces in 
terms of implicit equations. Such equations can be 
obtained by means of the computation of an elimination ideal, which turns out 
to be an ideal of equations of a certain Rees algebra. For a taste of 
this approach see, for instance, \cite{cox}, \cite{chw} and \cite{dac}, 
as well as the references therein.

\medskip

Among the many sources of the work presented in this dissertation, we want to 
highlight the classical work of Herzog, Simis and Vasconcelos in \cite{hsv}. 
This paramount piece of work, as well as the many references therein, keep being 
a fundamental source of valuable insights for the algebraic study of blowing-up.

\section{Scope of the present work}

We can split the scope of our work into a list of main achievements: 
\begin{enumerate}

\item[$(i)$] Establish conditions leading to procedures 
to describe explicitely the equations of blowing-up algebras 
of ideals of almost-linear type. The criteria 
displayed extend a set of well-known results about the 
equations of blowing-up algebras due to Vasconcelos, Trung, 
Huckaba, Heinzer and Kim.

\item[$(ii)$] Establish conditions under which the injectivity 
of a single component of the canonical morphism 
$\alpha:\symmetric(I)\to\rees(I)$ propagates downwards in the 
graded structure. These criteria lead to sufficient 
conditions for an ideal to be of linear type and provide a 
partial positive answer to a question posed by Tchernev.

\item[$(iii)$] Give examples and applications of the results
presented that may help to better understand the extent of the 
hypotheses considered.

\item[$(iv)$] Present additional work in progress, including
some preliminary results and ongoing experimentation, suggesting
future research pathways.

\end{enumerate}

We will make use of elementary linear and commutative algebra 
as well as standard results on (graded) Koszul homology. Among the 
specific tools used, the \emph{module of effective relations} of a 
standard algebra will play a major role. Although not so widespread 
across the literature, the module of effective relations is a very 
natural invariant to be considered in this context: it encodes useful 
information about the fresh generators of the equations of standard 
algebras and it can be defined canonically. This module can be seen 
as the first graded Koszul homology of a sequence of elements of a 
standard $R$-algebra. Playing with homological properties, 
we are able to establish the obstructions for such homology to 
vanish in terms of colon ideals.

Special attention will be paid to the degrees of the equations 
of blowing-up algebras: these are coarse numerical invariants. 
Low cost methods to know in advance what such degrees are for a 
given ideal, would be of great help for designing more efficient 
algorithms for computing these equations.

The conditions to be considered in the present work will involve 
most often hypotheses on the generating sets of ideals. We tackle the 
problem of describing the equations of $\rees(I)$ for $I=(J,y)$, 
with $y\notin J$ and $J$ being of linear type, i.e., for ideals of 
linear type up to one minimal generator. Throughout, such ideals will be referred 
to as ideals of \emph{almost-linear type}. Remark that taking an ideal $I$ of the 
form $(J,y)$ and assuming additional hypotheses on the ideal $J$ and the relation 
between $J$ and the element $y$, is a natural and common approach in the context 
(see, for instance, \cite{costa2}, \cite[Theorem~4.7]{hmv}, 
\cite[Proposition~3.9]{hsv}, \cite[Theorem~2.3 and Proposition~2.5]{valla}).

The numerical invariants of ideals which we will 
make use of explicitely are the deviation $\mu(I)-\height(I)$, 
the analytic spread $\spread(I)$, the analytic 
deviation $\adeviation(I)=\spread(I)-\height(I)$ and the second analytic deviation 
$\sadeviation(I)=\mu(I)-\spread(I)$. The techniques we will rely upon and the subsequent 
conditions on the generators of $I=(J,y)$ and $J$, as well as their 
relationship, will often imply small upper bounds on such invariants. 
Equimultiple ideals, i.e., ideals with analytic deviation $=0$, are among the 
classes of ideals we will come across. Such ideals include, 
for instance, $\fm$-primary ideals in a Noetherian 
local ring $(R,\fm)$.

\section{Preliminary definitions and results}\label{preliminary-results}

For the sake of self containment we devote a 
condensed section to summarise a set of well-known 
definitions and results that will support the proofs 
and discussions along the core exposition in Chapters
\ref{article2}, \ref{article1} and \ref{oddsandends}. 
Throughout, $R$ will denote a commutative
ring with unit $1\neq 0$, referred to as a commutative 
ring, for short.

\subsection{Blowing-up algebras}\label{sec:blowing-up-algebras}

By a standard $R$-algebra we mean a $\mathbb{N}$-graded commutative 
algebra over $R$ of the form $U=\oplus_{n\geq 0}U_n$, with $U_0=R$ and such that 
$U=R[U_1]$, with $U_{1}$ finitely generated as an $R$-module. We will 
denote by $z_{1},\ldots ,z_{s}$ a minimal generating set of $U_1$ as an
$R$-module. The irrelevant ideal of $U$, denoted by $U_+$, is the graded 
ideal generated by the elements of positive degree: $U_+=(U_1)=\oplus_{n>0}U_n$.

Let $I=(x_1,\ldots,x_s)$ be an ideal of $R$. Among the several standard algebras
associated to $I$, those known as blowing-up algebras will play a paramount role: 
the Rees algebra of $I$ is the $R$-algebra $\rees(I)=\oplus_{n\geq 0}I^{n}=R[It]\subset R[t]$; 
the associated graded ring of $I$ is the $R/I$-algebra 
$\graded(I)=\oplus_{n\geq 0}I^{n}/I^{n+1}=\rees(I)/I\rees(I)$; 
the fiber cone of $I$ with respect to a maximal ideal $\fm$ is the $R/\fm$-algebra 
$\fiber(I)=\oplus_{n\geq 0}I^n/\fm I^n=\rees(I)/\fm\rees(I)=\graded(I)/\fm\graded(I)$.

The symmetric algebra of an $R$-module $M$ is the standard $R$-algebra 
defined as $\symmetric_R(M)=\oplus_{n\geq 0}\symmetric_n(M)=T(M)/\mathfrak{S}$, 
where $T(M)=\oplus_{n\geq 0}M^{\otimes n}$ stands for the tensor 
algebra of $M$ and $\mathfrak{S}$ stands for the bilateral ideal of $T(M)$ 
generated by all the elements of the form $x\otimes y-y\otimes x$ for $x,y\in M$. 
Since $\mathfrak{S}$ is generated by homogeneous elements of degree $2$, 
$\symmetric_0(M)=R$ and $\symmetric_1(M)=M$ and there is a canonical 
injection $\iota:M\to\symmetric(M)$. For an account of basic properties of
symmetric algebras, see \cite[Chapter III, Section 6]{bourbaki}.

If $U$ is an $R$-standard algebra, we can canonically define a surjective 
morphism of standard $R$-algebras $\alpha:\symmetric(U_1)\to U$ induced by 
the inclusion $U_1\hookrightarrow U$ in degree one. When $U=\rees(I)$ we will say
that $\alpha:\symmetric(I)\to\rees(I)$ is the canonical blowing-up morphism. We 
define $\beta=\alpha\otimes \textbf{1}_{R/I}:\symmetric(I/I^2)\to \graded(I)$ 
and $\gamma=\alpha\otimes \textbf{1}_{R/\mathfrak{m}}:\symmetric(I/\mathfrak{m}I)\to\fiber(I)$.

\subsection{Sequential conditions}

In the sequel we will consider different assumptions on sequences of elements
in $R$. For a given sequence of elements $x_1,\ldots,x_s\in R$, we will 
consistently denote $J_i=(x_1,\ldots,x_i)$ and set $J_0=(0)$. 
Let us recall some well-known definitions and results that will be useful 
for our aims.

\begin{Definition}
A sequence of elements $x_1,\ldots,x_s\in R$ is said to be 
an $R$-sequence (or regular sequence in $R$) if the two 
following conditions hold:
\begin{enumerate}
\item[$(i)$] $((x_1,\ldots,x_i): x_{i+1})=(x_1,\ldots,x_i)$, for all 
$0\leq i\leq s-1$;
\item[$(ii)$] $(x_1,\ldots,x_s)\neq R$.
\end{enumerate}
\end{Definition}

\begin{Numbered}\label{regular-sequence-koszul}{\rm
If $I$ is an ideal of $R$ generated by an $R$-sequence
$x_1,\ldots,x_s$, then $\alpha:\symmetric(I)\to\rees(I)$ 
is an isomorphism. Moreover, $\symmetric(I)\cong \rees(I) \cong 
R[X_1,\ldots,X_s]/(x_iX_j-x_jX_i\;|\; 1\leq i<j\leq s)$. 
See \cite[Chapitre~Premier, Théorème~1 and Lemme~2 ]{micali}.
}\end{Numbered}

For many purposes, the elements of an $R$-sequence behave just like the 
variables in a polynomial ring over a field. This is exactly the case for 
a Noetherian local ring containing a field (see \cite{hartshorne2}). 

Recall that in a Noetherian local ring, $R$-sequences can be detected by 
means of the Koszul homology.

\begin{Numbered}{\rm
Let $(R,\fm)$ be a Noetherian local ring and let $x_1,\ldots,x_s$ 
be a sequence of elements in $R$. Then $x_1,\ldots,x_s$ is an 
$R$-sequence if and only if $\koszul{i}{x_1,\ldots,x_s}{R}=0$ 
for all $i\geq 1$. In particular, the Koszul complex gives rise 
to a free resolution of $R/(x_1,\ldots,x_s)$ of length $s$.
}\end{Numbered}

We will say that a sequence $x_1,\ldots,x_s\in R$ is a permutable 
$R$-sequence if $x_{\sigma(1)},\ldots,x_{\sigma(s)}$ is an $R$-sequence,
for every permutation $\sigma\in \mathfrak{S}_s$.
If an $R$-sequence is contained in the Jacobson radical then
it is permutable (in particular, in a local ring, every 
$R$-sequence is permutable). Monomials in permutable $R$-sequences 
behave just like monomials in indeterminates of a polynomial ring, 
in the sense that the basic arithmetic rules for monomial ideals 
hold true in a way that is made precise in the following statements
(see \cite{ks}).

\begin{Numbered}{\rm
Let $\underline{x}=x_1,\ldots,x_s$ be a permutable $R$-sequence. 
Let $x^a=x_1^{a_1}\cdots x_s^{a_s}$ and $x^b=x_1^{b_1}\cdots x_s^{b_s}$
be monomials in $\underline{x}$. Then we can define the greatest 
common divisor (gcd) and the least common multiple (lcm) of $x^a$ and $x^b$ as follows:
\begin{enumerate}
\item[] ${\rm gcd}(x^a,x^b)=x^c$, where $c_i={\rm min}\{a_i,b_i\}$ for all $1\leq i\leq s$;
\item[] ${\rm lcm}(x^a,x^b)=x^d$, where $d_i={\rm max}\{a_i,b_i\}$ for all $1\leq i\leq s$.
\end{enumerate}
}\end{Numbered}

\begin{Numbered}\label{monomial-ideal-operations}{\rm
Let $\underline{x}$ be a permutable $R$-sequence. Let $\mathfrak{a}=(m_1,\ldots,m_r)$,
$\mathfrak{b}=(n_1,\ldots,n_t)$ be ideals generated by monomials in $\underline{x}$ 
(referred to as monomial ideals in $\underline{x}$). Let 
$m$ be a monomial in $\underline{x}$. 
Then the following basic computational properties follow:
\begin{enumerate}
\item[$a)$] $\mathfrak{a}:m=\sum_{j=1}^r({\rm lcm}(m_j,m)/m)R$;
\item[$b)$] $\mathfrak{a}\cap mR=\sum_{j=1}^r{\rm lcm}(m_j,m)R$;
\item[$c)$] $\mathfrak{a}\cap\mathfrak{b}=\sum_{j=1}^r\sum_{i=1}^t {\rm lcm}(m_j,n_i)R$;
\item[$d)$] $\mathfrak{a}:\mathfrak{b}=\cap_{i=1}^t\sum_{j=1}^r ({\rm lcm}(m_j,n_i)/n_i)R$
\end{enumerate}
Moreover, $\mathfrak{a}$ is minimally generated by a unique generating set 
formed by monomials in $\underline{x}$.
}\end{Numbered}

If $x\in R\backslash \{0\}$, the leading form of $x$ in $\graded(I)=\oplus_{n\geq 0} I^n/I^{n+1}$, 
denoted by $x^*$, is the image of $x$ in $I^n/I^{n+1}$, where $x\in I^n \backslash I^{n+1}$. 
The following result was proved in \cite[Theorem 2.3]{vv} and will be used throughout:

\begin{Numbered}\label{valabrega-valla}{\rm
Let $R$ be a Noetherian ring, $I$ an ideal of $R$ and let 
$x_1,\ldots,x_s$ be an $R$-sequence, with $(x_1,\ldots,x_s)$ and $I$ being not comaximal. 
Then the leading forms $x_1^*,\ldots,x_s^*$ form a $\graded(I)$-sequence if and 
only if, for all $1\leq i\leq s$ and all $m\geq 1$, it is 
verified that 
\[ 
(x_1,\ldots,x_i)\cap I^m = \sum_{j=1}^i{x_jI^{m-\nu(x_j)}},
\label{valabrega-valla1} \tag{$VV$}
\]
where $\nu(x_j)$ stands for the degree of the initial form $x_j^*$ in $\graded(I)$ (see \cite[Section 2]{vv}).
Further, if $R$ is Noetherian local and $x_1,\ldots,x_s$ are part of a minimal
generating set of $I$, then $\nu(x_i)=1$ for all $1\leq i\leq s$. Thus,
setting $J_i=(x_1,\ldots,x_i)$, the condition (\ref{valabrega-valla1})
has the form 
\begin{eqnarray*}\label{valabrega-valla2}
J_i\cap I^m = J_iI^{m-1}.
\end{eqnarray*}
}\end{Numbered}

\begin{Definition}
If $J\subset I$ are ideals of $R$, we will define the $n$-th 
Valabrega-Valla module of $I$ with respect to $J$, 
denoted by $VV_J(I)_n$, as the quotient $J\cap I^n/JI^{n-1}$.
\end{Definition}

We introduce a notion that generalises regular sequences, the notion of $d$-sequence, 
which was introduced by Huneke and has shown to be useful in a wealth of scenarios.

\begin{Definition}\label{def:d-sequence}
A sequence $x_1,\ldots,x_s\in R$ is said to be a 
$d$-sequence if the following two conditions are fulfilled:
\begin{enumerate}
\item[$(i)$] $x_1,\ldots,x_s$ minimally generate the ideal $(x_1,\ldots,x_s)$ 
\item[$(ii)$] $((x_1,\ldots,x_i):x_{i+1}x_k)=((x_1,\ldots,x_i):x_k)$ for all
$0\leq i\leq s-1$ and $k\geq i+1$.
\end{enumerate}
\end{Definition}

It is readily seen that condition $(ii)$ in Definition \ref{def:d-sequence}
is equivalent to
\begin{enumerate}
\item[$(ii^*)$] $((x_1,\ldots,x_i):x_{i+1})\cap I=(x_1,\ldots,x_i)$ for $0\leq i\leq s-1$, 
where $I=(x_1,\ldots,x_s)$.
\end{enumerate}

\begin{Numbered}\label{d-sequences-are-linear-type}{\rm
If $I$ is an ideal of $R$ generated by a $d$-sequence, then
$\alpha:\symmetric(I)\to\rees(I)$ is an isomorphism (see the 
works of \cite[Theorem~3.1]{huneke2} and \cite[Theorem~3.15]{valla}).
}\end{Numbered}

\begin{Numbered}\label{length-of-d-sequences}{\rm
If $\height(I)=n$, then $I$ contains a $d$-sequence of both length 
and height $=n$ (see \cite[Proposition 6.2]{hsv}). But an ideal $I$ 
may also contain $d$-sequences of length $> \height(I)$: for instance, 
if $\fp$ is a prime ideal of $R$ with $\height(\fp)=n$ such that 
$\grade(\fp,R)=n$, $R_\mathfrak{p}$ is a regular local ring and 
$\fp$ has $n+1$ generators, then $\fp$ is generated by a $d$-sequence 
(see \cite[Exercise 5.24]{sh}). In any case, if $\dim R=r$ then every 
$d$-sequence has at most $r+1$ elements (see \cite[6.1.3]{hsv}). 
}\end{Numbered}

\medskip

Within a polynomial ring over a field, there is a manageable criterion to 
decide whether a given sequence consisting of monomials is a $d$-sequence 
(see \cite{tang}).

\begin{Numbered}{\rm
Let $R=k[X_1,\ldots,X_s]$ be a polynomial ring over a field $k$ and 
let $f_1,\ldots,f_s\in R$ be a sequence consisting of monomials. Then
$f_1,\ldots,f_s$ is a $d$-sequence if and only the following three
conditions hold:
\begin{enumerate}
\item $f_1,\ldots,f_s$ is a minimal generating set of $I=(f_1,\ldots,f_s)$;
\item ${\rm gcd}(f_i,f_j)\; |\; f_k$, for all $1\leq i<j<k\leq s$;
\item ${\rm gcd}(f_i,f_j)={\rm gcd}(f_i,f_j^2)$, for all $1\leq i<j\leq s$.
\end{enumerate}
}\end{Numbered}

\subsection{Reduction of ideals}

\begin{Definition}
The ideal $J$, $J\subset I$, is said to be a reduction of $I$ if 
there is an integer $r\geq 0$ such that $I^{r+1}=JI^r$. 
The reduction number of $I$ with respect to $J$ is the
least such integer $r$, denoted by $\rednumber_J(I)$. It is said 
that $J$ is a minimal reduction of $I$ if it is a 
reduction and it is minimal among all reductions of $I$ with respect 
to inclusion, i.e., no ideal strictly contained in $J$ is a 
reduction of $I$.
\end{Definition}

\begin{Numbered}{\rm
Remark that if $J$ is a reduction of $I$ then 
${\rm Min}(R/I)={\rm Min}(R/J)$, thus 
$\sqrt{I}=\sqrt{J}$ and $\height(I)=\height(J)$. 
}\end{Numbered}

Reductions of ideals can be detected by means of extension properties 
of the corresponding blowing-up algebras.

\begin{Numbered}\label{reductions-blowing-up}{\rm
Let $(R,\fm)$ be a Noetherian local ring with residue field $k=R/\fm$. 
Let $J\subset I$ be ideals in $R$. Then $J$ is a
reduction of $I$ if and only if any of the following 
equivalent conditions hold:
\begin{enumerate}
\item[$(i)$] $\rees(I)$ is a finitely generated $\rees(J)$-module, i.e., 
the extension $\rees(J)\subset \rees(I)$ is module-finite;
\item[$(ii)$] If $k[J+\fm I/\fm I]$ stands for the subalgebra of $\fiber(I)$ generated 
over $k$ by $J+\fm I/\fm I$, the extension $k[J+\fm I/\fm I]\subset \fiber(I)$ is module-finite.
\end{enumerate}
See \cite[Theorem~8.2.1 and Proposition~8.2.4]{sh}.
}\end{Numbered}

\begin{Definition}
The analytic spread of $I$, denoted by $\spread(I)$, is defined
as the Krull dimension of $\fiber(I)$.
\end{Definition}

It immediately follows that the degree of the Hilbert polynomial 
of $\fiber(I)$ is $\spread(I)-1$, thus $\spread(I)$ can be 
interpreted as the rate of growth of 
$\lambda(I^n/\mbox{\fm}I^n)=\mu(I^n)$ as $n$ increases. 

\begin{Numbered}{\rm
Let $(R,\fm)$ be a Noetherian local ring with infinite residue field $k$.
Let $J$ be a minimal reduction of $I$. Then,
\begin{enumerate}
\item[$(i)$] $\mu(J)=\spread(I)$;
\item[$(ii)$] any minimal generating set of $J$ extends to a minimal 
generating set of $I$.
\item[$(iii)$] $\fiber(J) \cong \symmetric(J/\fm J)$.
\end{enumerate}
}\end{Numbered}

Let $(R,\fm)$ be a Noetherian local ring. Throughout, 
we will raise almost-linear type assumptions, 
where $I=(J,y)$, with $y\notin J$ and $J\subset I$ will be 
a reduction generated either by an $R$-sequence or by a 
$d$-sequence. Notice that if $J$ is generated by an $R$-sequence of length $s$, 
then $s\leq \grade(I)\leq \height(I)\leq \spread(I)\leq \mu(J)=s$, 
whence $\adeviation(I)=\spread(I)-\height(I)=0$ and 
$\sadeviation(I)=\mu(I)-\spread(I)=1$. On the other hand, 
if $I$ is an $\fm$-primary ideal with $J$ minimally generated 
by a $d$-sequence $x_1,\ldots,x_s$ of length $s$, 
we have $s\leq \dim R+1$ by \ref{length-of-d-sequences}.
Then, since $I$ is $\mathfrak{m}$-primary, 
we have $s\leq \dim R+1=\height(I)+1=\spread(I)+1 \leq \mu(J)+1=s+1$.
Consequently, $s-1\leq\height(I)=\spread(I)\leq s$.

\subsection{The equations of blowing-up algebras}

\subsubsection{Relation type}
Let $R$ be a Noetherian ring and let $U=R[z_1,\ldots,z_s]$ be a 
standard $R$-algebra. Let $V=R[T_1,\ldots ,T_s]$ be a 
polynomial ring with variables $T_{1},\ldots ,T_s$ and let 
$\varphi:V\rightarrow U$ be the polynomial presentation of $U$ 
sending $T_{i}$ to $z_{i}$. Let $Q=\oplus_{n\geq 1}Q_{n}$ be 
the kernel of $\varphi$, whose elements will be referred to as 
the equations of $U$. Let $Q\langle n\rangle$ be the ideal generated 
by the homogeneous equations of $U$ of degree at most $n$. 
Notice that there is an ascending chain 
$Q\langle 1 \rangle \subset Q\langle 2 \rangle \subset \ldots$
that stabilises at $Q\langle n\rangle$, for some integer $n\geq 1$.

\begin{Numbered}{\rm
Remark that the following isomorphisms of $R$-algebras hold: 
$V/Q\langle 1 \rangle \cong \symmetric_R(U_1)$ and $V/Q \cong U$,
where $\symmetric_R(U_1)$ is the symmetric algebra of $U_1$.
}\end{Numbered}

\begin{Definition}
The relation type of $U$, denoted by $\reltype(U)$, is 
the least integer $N\geq 1$, such that $Q\langle N\rangle=Q$, 
i.e., the maximum degree in a minimal generating set of 
homogeneous elements of $Q$.
\end{Definition}

The case where $U=\rees(I)$ will be of 
special importance for our aims. We will denote 
$\reltype(\rees(I))$ more simply as $\reltype(I)$, 
the relation type of $I$.

\begin{Numbered}{\rm
While $Q$ clearly depends on the polynomial presentation $\varphi:V\to U$, 
the relation type does not depend on the presentation of $U$. Moreover, even
the degrees of a minimal generating set of homogeneous elements of $Q$ remain
invariant with respect to the presentation 
(see Remarks \ref{explicit-presentation} and \ref{min-gen-set-of-equations}).
}\end{Numbered}

\medskip

In what follows, let $U=\rees(I)$, $\varphi:V\rightarrow \rees(I)$ 
be the polynomial presentation of $\rees(I)$ sending $T_{i}$ to $x_{i}t$ 
and $Q=\ker\varphi$. 
\begin{Numbered}\label{rees-elimination}{\rm
The ideal $Q$ can be computed by means of elimination 
theory. Let $I=(x_1,\ldots,x_s)$, $V=R[T_1,\ldots,T_s]$ 
and let $T$ be an indeterminate over $V$. Then,
\begin{enumerate}
\item[$(i)$] $Q=(T_1-x_1T,\ldots,T_s-x_sT)V[T]\cap V$;
\item[$(ii)$] if $x_1$ is a non-zero-divisor of $R$, 
\begin{align*}
Q & = ((x_2T_1-x_1T_2,\ldots,x_sT_1-x_1T_s):_V x_1^\infty) \\
  & = (x_2T_1-x_1T_2,\ldots,x_sT_1-x_1T_s,x_1T-1)V[T]\cap V.
\end{align*}
\end{enumerate}
}\end{Numbered}

In the context of polynomial rings these computations can be done via 
Gröbner basis techniques (see e.g. \cite{froberg}, \cite{clo}). 
Recall that if $L\subset k[x_1,\ldots,x_s,y_1,\ldots,y_r]$ is an ideal and 
$\leq_e$ is a monomial ordering such that ${\bf x}^{\sigma}\leq_e y_i$ for 
any $\sigma\in\mathbb{N}^s$ and $1\leq i\leq r$, and 
$G=\{f_1,\ldots,f_k\}$ is a Gröbner basis of $L$ with respect to 
$\leq_e$, then $G\cap k[x_1,\ldots,x_s]$ is a Gröbner basis of 
$L\cap k[x_1,\ldots,x_s]$ (see \cite[Proposition~8, p. 95]{froberg}). 
A monomial ordering $\leq_e$ fulfilling the conditions above is said to be 
an \emph{elimination ordering}.

The expected complexity of computing $Q$ is high. The underlying 
rationale is that the theoretical complexity of computing Gröbner bases 
is generically exponential and doubly-exponential 
in the worst case (see \cite{lazard} and \cite{mm}), 
as a function of the number of variables (see also \cite[Appendix]{bw} 
for a further account on Gröbner bases and complexity issues). 
On the other hand, the cost of the computations in \ref{rees-elimination} 
heavily relies on the number of generators of $I$ 
(see \cite{bsv}). 

\medskip

The following definition canonically encloses the condition that 
$Q=Q\langle 1 \rangle$, i.e., $\reltype(I)=1$ (see comments 
in \ref{effective-and-q}).

\begin{Definition}\label{definition-of-linear-type}
An ideal $I$ of $R$ is said to be of linear type if 
$\alpha:\symmetric(I)\to\rees(I)$ is an isomorphism.
\end{Definition}

\begin{Numbered}\label{arithmetic-linear-type}{\rm
Let $I$ be an ideal of linear type and let $\fp\supset I$ a prime 
ideal containing $I$. The isomorphism $\alpha:\symmetric(I)\to\rees(I)$ 
induces an isomorphism $\beta_\mathfrak{p}:\symmetric(I_\mathfrak{p}/I_\mathfrak{p}^2)\to\graded(I_\mathfrak{p})$, 
hence the Krull dimensions of $\symmetric(I_\mathfrak{p}/I_\mathfrak{p}^2)$ 
and $\graded(I_\mathfrak{p})$ coincide and are equal to $\dim R_\mathfrak{p}$
(see \cite[Proposition~5.1.6]{sh}). By the surjectivity of
$\symmetric(I_\mathfrak{p}/I_\mathfrak{p}^2)\to \symmetric(I_\mathfrak{p}/\fp R_\mathfrak{p}I_\mathfrak{p})$ 
it follows that $\mu(I_\mathfrak{p})\leq \height(\fp)$ (see also \cite[Proposition 2.4]{hsv}).
}\end{Numbered}

\begin{Numbered}{\rm
Let $I$ be an ideal of linear type. Then, for each prime ideal 
$\fp\supset I$, $\spread(I_\mathfrak{p})=\mu(I_\mathfrak{p})$. 
The converse, however, is not true in general.
Consider $R=k[[X,Y]]/(X^2Y, XY^2)=k[[x,y]]$ and let $\fp=(y)\subset \mathfrak{m}=(x,y)\subset R$. 
Then $\fp$ is prime with $\spread(\fp R_\mathfrak{p})=\mu(\fp R_\mathfrak{p})=0$ and
$\spread(\fp R_\mathfrak{m})=\mu(\fp R_\mathfrak{m})=1$, but $\ker\alpha_n\neq 0$ 
for all $n\geq 2$ (see also Example~\ref{principal-g-linear-but-not-p-linear}).
}\end{Numbered}

\begin{Definition}
An ideal $I$ is said to be syzygetic if $Q\langle 1 \rangle = Q\langle 2 \rangle$, 
i.e., $\alpha_2:\symmetric_2(I)\to I^2$ is an isomorphism. 
\end{Definition}

\begin{Numbered}{\rm
A syzygetic ideal need not be an ideal of linear type. For instance,
any prime ideal of height $2$ in a $3$-dimensional regular local ring is syzygetic 
(see \cite[Proposition~2.7 and subsequent Remark]{hsv}). In the formal power 
series ring $k[[X,Y,Z]]$, there exist prime ideals $\fp$ of height $2$ with arbitrarily 
large minimal number of generators (primes of Moh, see \cite{moh}). 
Consequently, such prime ideals are syzygetic but they do not fulfil the inequality 
$\mu(\fp)=\mu(\fp_\mathfrak{m})\leq \height(\fm)=3$ (see claim \ref{arithmetic-linear-type} 
above), thus they are not of linear type.
}\end{Numbered}

\subsubsection{Effective relations}\label{sec:effective-relations}

Let $R$ be a Noetherian ring, let $U$ be a standard $R$-algebra and let 
$\alpha:\symmetric(U_1)\to U$ be the canonical morphism. 

\begin{Definition}\label{definicio-relacions-efectives}
If $n\geq 2$ is an integer, define the $n$-th \emph{module of effective relations} 
of $U$ as $E(U)_n=\ker\alpha_n/\symmetric_1(U_1)\cdot\ker\alpha_{n-1}$. If $I$ is an
ideal of $R$ and $U=\rees(I)$, we set $E(I)_n=E(\rees(I))_n$.
\end{Definition}

\begin{Numbered}\label{effective-and-q}{\rm
The module $E(U)_n$ can also be described by means of a polynomial presentation 
$\varphi:V\to U$. If $Q=\ker\varphi$, then $E(U)_n\cong(Q/Q\langle n-1\rangle)_{n}=Q_{n}/V_{1}Q_{n-1}$, 
for $n\geq 2$ (see \cite{planas}). In particular, taking $U=\rees(I)$, one can see that 
$I$ is of linear type if and only if $E(U)_n=0$ for all $n\geq 2$. Therefore, $I$ is of linear type
if and only if $Q=Q\langle 1\rangle$, i.e., $\reltype(I)=1$ (see 
Definition~\ref{definition-of-linear-type}).
}\end{Numbered}

The information encoded by the module of effective relations accounts for 
the fresh generators of $\ker\alpha$ or $Q$. The non-zero elements of 
$E(U)=\oplus_{n\geq 2}E(U)_n$ correspond to the equations of $U$ 
of degree $\geq 2$ which are not combination of the equations in lower degree. 
Consequently, if $R$ is Noetherian, $E(U)_n=0$ for $n$ large enough.

\begin{Numbered}\label{effective-relations-homological}{\rm
Let $R$ be a Noetherian ring, let $U=R[z_1,\ldots,z_s]$ be a standard $R$-algebra and let
$\underline{z}=z_1,\ldots,z_s$. Then: 
\begin{enumerate}
\item[$(i)$] The module of effective relations can be expressed in terms of the graded 
Koszul homology: $E(U)_n\cong\koszul{1}{\underline{z}}{U}_n$, for $n\geq 2$.
\item[$(ii)$] The relation type $\reltype(U)$ can be given in terms of the module of 
effective relations: $\reltype(U)=\min\{s\geq 1\; |\; E(U)_k=0,\;\; \textrm{for each}\;\; k\geq s+1\}$.
\end{enumerate}
}\end{Numbered}

\begin{Numbered}\label{defining-equations-associated-graded}{\rm
Let $R$ be a Noetherian ring, $I$ an ideal of $R$. The following properties hold:
\begin{enumerate}
\item[$(i)$] There is an exact sequence $E(I)_{n+1}\to E(I)_n\to E(\graded(I))_n\to 0$, 
for all $n\geq 2$;
\item[$(ii)$] As a consequence of $(i)$, $\reltype(I)=\reltype(\graded(I))=r$ 
and $E(I)_r\cong E(\graded(I))_r$.
\end{enumerate}

See \cite[Theorem~1.3]{valla} for the case of $\reltype(I)=1$; see 
\cite[Section~3]{planas} for the general case using effective relations; 
and see \cite[Discussion~2.2]{hku} for the same result but using 
extended Rees and symmetric algebras. One of the advantages of using 
effective relations is that one is able to deduce at once that the top 
degree equations of $\rees(I)$ are in correspondence with the top degree 
equations of $\graded(I)$.
}\end{Numbered}

\subsubsection{Geometric conditions close to linear type}\label{geometric-conditions-linear-type}

We shall introduce two notions which are close to the linear type condition 
that will play a role in our study: geometric linear type (or $g$-linear type) 
and projectively of linear type (or $p$-linear type). These notions come from
considering geometric properties linked to the canonical epimorphism 
$\alpha:\symmetric(I)\to \rees(I)$.

In Hermann-Moonen-Villamayor \cite{hmv}, the authors consider 
the condition of $g$-linear type, which is weaker than that of linear type 
in the following sense: while linear type imply that the normal cone 
and the normal bundle of a closed subscheme of a given scheme coincide, 
the condition $g$-linear type refers only to the reduced structures of the 
normal cone and normal bundle. 

\begin{Definition}
The ideal $I$ is said to be of geometric linear type ($g$-\emph{linear type}, 
for short) if the morphism induced by $\alpha$, 
$\alpha^*:{\rm Spec}(\rees(I))\to{\rm Spec}(\symmetric(I))$, is a
homeomorphism of topological spaces.
\end{Definition}

It is readily seen that since $\rees(I)\cong \symmetric(I)/\ker\alpha$, if $I$ is of
$g$-linear type, then $\alpha^*$ induces a homeomorphism 
$V(\ker\alpha)\to {\rm Spec}(\symmetric(I))$. Consequently, the minimal
primes in $V(\ker\alpha)$ correspond to the minimal primes in $\symmetric(I)$, 
thus $\ker\alpha\subset \bigcap_{\mathfrak{p}\in {\rm Min}(S(I))} \fp$, 
i.e., $\ker\alpha$ is nilpotent.

\begin{Numbered}{\rm
The ideal $I$ is of $g$-linear type if and only if the following equivalent 
conditions hold (see \cite[Proposition~1.1]{hmv}):
\begin{itemize}
\item[$(i)$] $\ker\alpha$ is nilpotent;
\item[$(ii)$] $\ker\beta$ is nilpotent (see \cite[Theorem~3.2]{hsv});
\item[$(iii)$] $\mu(I_{\mathfrak{p}})=\spread(I_{\mathfrak{p}})$, for all $\fp\in{\rm Spec}(R)$;
\end{itemize}
}\end{Numbered}
Clearly, if $I$ is of linear type, then it is also of $g$-linear type, 
whereas the converse is not true in general. One of the advantages of 
introducing the notion of $g$-linear type is that this property can be 
checked by numerical conditions, which provides a general strategy of 
proving linear type by first proving $g$-linear type and then looking for 
additional conditions which assure that linear type and $g$-linear 
type coincide \cite{hmv}.

\begin{Definition}
An ideal $I$ is said to be projectively of linear type (of $p$-linear type, for short) 
if the morphism $\widetilde{\alpha}:{\rm Proj}(\rees(I))\to{\rm Proj}(\symmetric(I))$ 
is an isomorphism of schemes, where $\widetilde{\alpha}$ is the morphism induced 
by $\alpha$.
\end{Definition}

Recall that for a standard $R$-algebra $U$,
${\rm Proj}(U)=\{ \fq\in {\rm Spec}(U)\; |\; \fq\;{\rm homogeneous},\, U_+\not\subset \fq \}$ 
and, if $\alpha$ is surjective, $\widetilde{\alpha}: {\rm Proj}(V)\to{\rm Proj}(U)$ 
is a morphism of schemes given locally by the morphisms of affine schemes 
$\alpha^*_{u,0}:{\rm Spec}((V_{\alpha(u)})_0)\to{\rm Spec}((U_u)_0)$, 
with $u\in U_+$. It is clear that an ideal of linear type is of $p$-linear type. 

\begin{Numbered}\label{p-linear-type-giral}{\rm
Let $\alpha:U\to V$ be a surjective morphism of standard $R$-algebras. Then, 
one can show that $\widetilde{\alpha}:{\rm Proj}(V) \to {\rm Proj}(U)$ is an 
isomorphism of schemes if and only if any of the following equivalent statements 
hold:
\begin{enumerate}
\item[$(i)$] there is an integer $r\geq 0$ such that 
$(U_+)^r\ker\alpha=0$; 
\item[$(ii)$] $\ker\alpha_n=0$ for $n$ sufficiently large. 
\end{enumerate}
Consequently, if $I$ is of $p$-linear type, it is of $g$-linear type.
}\end{Numbered}

\medskip

It is known that linear type, $p$-linear type and $g$-linear type
do not coincide in general. Kühl in \cite{kuhl} was the first to highlight
an example of an ideal of $p$-linear type, but not of linear type. 

\begin{Example}\label{exemple-kuhl}{\rm (see \cite[Example 1.4]{kuhl}). 
Let $A=k[U_0,U_1,U_2]$ and put 
$$R=A[X,Y]/Q=A[x,y],\,{\rm where}\; Q=(U_0X,U_1Y,U_0Y^2,U_1X^2,U_1X-U_2Y,U_0Y+U_2X).$$ 
Then, the ideal $I=(x,y)$ of $R$ is not of linear type, i.e., 
$\alpha:\symmetric(I)\to\rees(I)$ is not an isomorphism; however, 
$I$ is of $p$-linear type, i.e., 
${\rm Proj}(\alpha):{\rm Proj}(\rees(I))\to {\rm Proj}(\symmetric(I))$ 
is an isomorphism of schemes.
}\end{Example}

One can see that if $I$ is principal, then $I$ is of $p$-linear type if and 
only if $I$ is of linear type. However, a similar claim is not true when 
comparing $g$-linear type and $p$-linear type. Later, in 
Example~\ref{principal-g-linear-but-not-p-linear}, we will give a principal 
ideal which is of $g$-linear type, but not of $p$-linear type.
\chapter{The equations of Rees algebras of ideals of almost-linear type}\label{article2}

\section{Introduction}\label{article2:introduction}

The aim of this chapter is to explicitly describe the equations of 
Rees algebras of classes of ideals of almost-linear type. Let us first 
recall the basic definitions in order to properly introduce the main 
result of this chapter (see also Section~\ref{preliminary-results}).

Recall that the reduction number of $I$ with respect to $J$, $J\subset I$,
denoted by $\rednumber_J(I)$, is the least integer $r\geq 0$ such that 
$I^{r+1}=JI^r$. Set $V=R[X_1,\ldots,X_s,Y]$ and let $\varphi:V\to\rees(I)$ 
be the polynomial presentation of $\rees(I)$, the Rees algebra of $I$, 
sending $X_i$ to $x_it$ and $Y$ to $yt$. Let $Q=\oplus_{n\geq 1}Q_n$ 
be the kernel of $\varphi$. Given an integer $m\geq 0$, set 
$Q\langle m\rangle\subset Q$ the ideal generated by the homogeneous 
elements of $Q$ of degree at most $m$ in $X_1,\ldots,X_s,Y$. The 
relation type of $I$, denoted by $\reltype(I)$, is the least integer 
$N\geq 1$ such that $Q=Q\langle N\rangle$. Let $\graded(I)=\oplus_{n\geq 0} I^n/I^{n+1}$ 
be the associated graded ring of $I$. If $z\in R\backslash \{0\}$, we will denote by $z^*$
the initial form of $z$ in $\graded(I)$, i.e., $z^*=z+I^{n+1}$, 
where $z\in I^n\backslash I^{n+1}$.

\begin{Theorem A}
Let $(R,\fm)$ be a Noetherian local ring and let $I$ be an ideal of
$R$. Let $x_1,\ldots,x_s,y$ be a minimal generating set of $I$, where
$J=(x_1,\ldots,x_s)$ is a reduction of $I$ with reduction number
$r=\rednumber_J(I)$. Assume that $x_1,\ldots,x_s$ verify 
the following condition for all $n\geq 2$:
\[ 
((x_1,\ldots,x_{i-1})I^{n-1}: x_i)\cap I^{n-1}=(x_1,\ldots,x_{i-1})I^{n-2}, 
\textrm{ for all } i=1,\ldots,s. \label{eq:trung-condition} \tag{$\mathcal{T}_n$}
\]
Then, for each $n\geq 2$, the map sending $F\in Q_{n}$ to
$F(0,\ldots,0,1)\in (JI^{n-1}: y^{n})$ induces an isomorphism of
$R$-modules
\begin{eqnarray*}
\left[ \frac{Q}{Q\langle n-1\rangle}\right]_{n}\cong
\frac{(JI^{n-1}: y^{n})}{(JI^{n-2}: y^{n-1})}.
\end{eqnarray*}
In particular, $\reltype(I)=\rednumber_J(I)+1$ and there is a form
$Y^{r+1}-\sum X_iF_i\in Q_{r+1}$, with $F_i\in V_r$, such that
$Q=(Y^{r+1}-\sum X_iF_i)+Q\langle r\rangle$. Moreover, if 
$x_1,\ldots,x_s$ is an $R$-sequence and $x_1^*,\ldots,x_{s-1}^*$ is
a $\graded(I)$-sequence, then $x_1,\ldots,x_s$ verify the condition $(\mathcal{T}_n)$
for all $n\geq 2$.
\end{Theorem A}


Roughly speaking, the theorem says how to obtain a minimal generating 
set of the equations of $\rees(I)$. For the equations of degree 1, 
pick a minimal generating set of the first syzygies of $I$, viewed as 
elements of $Q_{1}$. For the equations of higher degree $=n$, with 
$2\leq n\leq r+1$, take representatives of the inverse images of a 
minimal generating set of $(JI^{n-1}:y^{n})/(JI^{n-2}:y^{n-1})$ 
(see Remark~\ref{min-gen-set-of-equations} and Example~\ref{exemple-classic}).

\medskip

Thus far, the study of the equations of $\rees(I)$ has produced a vast
literature. Some of this work is focused on ideals having small
deviations as well as on the interplay between the reduction number
and the relation type (see, just to mention some references,
\cite{huckaba1}, \cite{huckaba2}, \cite{hh1}, \cite{hh2},
\cite{trung1}, \cite{trung2}, \cite{vasconcelos1},
\cite{vasconcelos2}). The particular hypotheses and interests in this 
chapter owe much to the works of Cortadellas and Zarzuela in \cite{cz},
Heinzer and Kim in \cite{hk}, Huckaba in \cite{huckaba1} and
\cite{huckaba2}, Trung in \cite{trung1} and \cite{trung2}, and
Vasconcelos in \cite{vasconcelos1}. In fact, Theorem~A sprouted as an 
attempt to understand \cite[Theorem~2.3.3]{vasconcelos1}. That is the 
reason for considering ideals $I$ of the form $I=(J,y)$, with $y\notin J$ 
and $J$ being a reduction of $I$ generated by a sufficiently good sequence that 
makes $J$ an ideal of linear type. The reader may also consult 
\cite{hsv1} for a recent account on the equations of 
$\rees(I)$ with similar assumptions.

\medskip

Let $R$ be a Noetherian local ring. Even for ideals of the principal 
class, i.e., $\mu(I)-\height(I)=0$, the equations of $\rees(I)$ may 
be difficult to describe. Remarkably, any ideal $I$ of the principal 
class is generated by an $R$-sequence $x_{1},\ldots,x_{s}$, provided 
that $I$ is prime (see \cite{davis}) or $R$ is Cohen-Macaulay. In both 
cases the equations are generated by the Koszul relations 
$x_iX_j-x_jX_i$, with $1\leq i<j\leq s$ (see \ref{regular-sequence-koszul}). 
In particular, $I$ is of linear type. However, if $I$ is not prime and $R$ is not
Cohen-Macaulay, this is no longer true. For instance, consider the
ideals generated by systems of parameters. Huneke asked in
\cite{huneke} whether there is a uniform bound for the relation type
of these ideals in a complete local equidimensional Noetherian ring
$R$. The full answer to this question was given in \cite{wang},
\cite{lai} and \cite{agh}. Concretely, in \cite[Example 2.1]{agh}, it
was shown that if the non-Cohen-Macaulay locus of $R$ has dimension $\geq 2$, 
there exist families of parametric ideals of $R$ with
unbounded relation type. This gives an idea of the complexity of the
structure of the equations of $\rees(I)$.

The case $I=(J,y)$, with $y\notin J$, $J=(x_1,\ldots,x_s)$ a reduction 
of $I$, $x_1,\ldots,x_s$ is an $R$-sequence and $x_1^*,\ldots,x_{s-1}^*$ 
a $\graded(I)$-sequence, considered in Theorem~A, verifies 
$\mu(I)-\height(I)=1$. In general, if $\mu(I)-\height(I)=1$ and $R$ is Cohen-Macaulay, 
$I$ is of linear type if and only if $I$ is locally of the
principal class at all minimal primes of $I$ (see
\cite[Theorem~4.8]{hmv}). However, if $I$ is not locally of the principal
class at its minimal primes, then the relation type may be arbitrarily
large. For instance, take $R=k[[a,b]]$, the power series ring
in two variables over a field $k$ and let 
$I=(a^{p},b^{p},ab^{p-1})$, with $p\geq 2$. Then, $I$ 
is $(a,b)$-primary and verifies $\mu(I)-\height(I)=1$, 
hence it is not locally of the principal class at its minimal primes. 
One can check that $\reltype(I)\geq p$. In fact, $I$ fulfils the 
hypotheses of Theorem~A, $J=(a^{p},b^{p})$ 
being a reduction of $I$. Thus, the containment 
$(ab^{p-1})^{p}\in JI^{p-1}$, induces the only equation of 
degree $=p$ in a minimal generating set of equations of $\rees(I)$ and 
$\reltype(I)=p$ (see Example~\ref{exemple-classic}).

As highlighted before, this study has its origins in the following 
result of Vasconcelos in \cite[Theorem~2.3.3]{vasconcelos1}: let $(R,\fm)$ be a
Cohen-Macaulay local ring of dimension $d$ and let $I$ be an
$\fm$-primary ideal of $R$. Let $x_{1},\ldots, x_{d},y$ be a minimal
generating set of $I$, where $J=(x_{1},\ldots ,x_{d})$ is a reduction
of $R$ with reduction number $\rednumber_J(I)=1$. Then there is a form
$Y^{2}-\sum X_iF_i\in Q_{2}$, with $F_i\in V_1$, such that
$Q=(Y^{2}-\sum X_iF_i)+Q\langle 1\rangle$. In particular,
$\reltype(I)=2$. Remark that the hypotheses imply that 
$x_{1},\ldots,x_{d}$ is an $R$-sequence and that the initial forms
$x_1^*,\ldots,x_d^*$ are a $\graded(I)$-sequence (see the result of
Valabrega and Valla in \cite[Proposition~3.1]{vv}). By
Theorem~A, it is enough to suppose that $x_1^*,\ldots,x_{d-1}^*$ 
is a $\graded(I)$-sequence. Moreover, one can consider any 
reduction number (see Corollary~\ref{crelle-generalisation}).

As a corollary of Theorem~A, we recover a result of Heinzer and Kim 
in \cite[Theorem~5.6]{hk}, where they prove that the equations of 
the fiber cone of $I$, $\fiber(I)=\rees(I)\otimes R/\fm=\oplus_{n\geq 1}I^n/\fm I^n$, 
are generated by a unique equation of degree $\rednumber_J(I)+1$.

\section{Effective relations of a standard algebra}\label{equations-of-u}

By a standard $R$-algebra $U$ we mean a graded $R$-algebra
$U=\oplus_{n\geq 0}U_{n}$, with $U_{0}=R$, $U=R[U_{1}]$ and $U_{1}$
minimally generated by $\underline{z}=z_{1},\ldots ,z_{s}\in U_1$ as an
$R$-module. For instance, the $R$-algebra $\rees(I)=\oplus_{n\geq
0}I^{n}$, the $R/I$-algebra $\graded(I)=\oplus_{n\geq
0}I^{n}/I^{n+1}$ and the $R/\fm$-algebra $\fiber(I)=\oplus_{n\geq
0}I^n/\fm I^n$, are standard algebras.

Let $V=R[T_1,\ldots ,T_s]$ be a polynomial ring with variables
$T_{1},\ldots ,T_s$, let $\varphi:V\rightarrow U$ be the induced
presentation of $U$ sending $T_{i}$ to $z_{i}$ and let $Q=\oplus_{n\geq 1}Q_{n}$ 
be the kernel of $\varphi$, the elements of which will be referred
to as the equations of $U$. Let $\symmetric(U_{1})$ be the symmetric 
algebra of $U_{1}$ and let $\alpha:\symmetric(U_{1})\rightarrow U$ be 
the canonical morphism induced by the identity in degree one. 
Given $n\geq 2$, the module of effective $n$-relations of $U$ is defined 
to be $E(U)_{n}=\ker(\alpha_{n})/U_{1}\ker(\alpha_{n-1})$. One can prove
that, for any $n\geq 2$, $(Q/Q\langle n-1\rangle)_{n}\cong E(U)_{n}$
(see \cite[Definition~2.2]{planas}). In particular, the relation type
of $U$ is the least integer $N\geq 1$, such that
$E(U)_{n}=0$ for all $n\geq N+1$.

This description of $(Q/Q\langle n-1\rangle)_{n}$ as
$E(U)_{n}=\ker(\alpha_{n})/U_{1}\ker(\alpha_{n-1})$ has the advantage
of being canonical. On the other hand, it is the bridge to think about
$(Q/Q\langle n-1\rangle)_{n}$ as a Koszul homology module 
(see \cite[Corollary~2.7]{planas}). More precisely, for $n\geq
2$, we have the isomorphisms $(Q/Q\langle n-1\rangle)_{n}\cong E(U)_{n}\cong \koszul{1}{\underline{z}}{U}_{n}$,  
where $\koszul{1}{\underline{z}}{U}_n$ denotes the first homology module of the complex:
\begin{eqnarray*}
\ldots \to \wedge_{2}(R^{s})\otimes
U_{n-2}\stackrel{\partial_{2,n-2}}{\longrightarrow}
\wedge_{1}(R^{s})\otimes
U_{n-1}\stackrel{\partial_{1,n-1}}{\longrightarrow}U_{n}\to 0,
\end{eqnarray*}
where the Koszul differentials are defined as follows: if $e_1,\ldots
,e_s$ stands for the canonical basis of $R^s$, $u\in U_{n-2}$ and
$v\in U_{n-1}$, then
\begin{eqnarray*}
\partial_{2,n-2}((e_i\wedge e_j)\otimes u)=e_j\otimes (z_i\cdot u)-e_i\otimes (z_j\cdot u) 
\; \mbox{ and } \; \partial_{1,n-1}(e_i\otimes
v)=z_i\cdot v.
\end{eqnarray*}
As usual, $\cycle{1}{\underline{z}}{U}$ and $\boundary{1}{\underline{z}}{U}$ 
will stand for the graded modules of $1$-cycles and $1$-boundaries, respectively, of
the Koszul complex of $\underline{z}$. Observe that $\boundary{1}{\underline{z}}{U}_{1}=0$.

\begin{Remark}\label{explicit-presentation}{\rm 
For each integer $n\geq 2$, there is a well-defined isomorphism 
$$\tau_n: (Q/Q\langle n-1\rangle)_{n}\stackrel{\cong}{\to}\koszul{1}{\underline{z}}{U}_{n}$$ 
sending the class of $F\in Q_n$ modulo $Q\langle n-1 \rangle$ to the homology
class of $(F_1(z),\ldots,F_s(z))\in
\cycle{1}{\underline{z}}{U}_{n}\subset U_{n-1}\oplus \stackrel{s)}{\ldots} \oplus U_{n-1}$, 
where $F_1,\ldots,F_s$ is any set of $s$ elements in $V_{n-1}$ satisfying 
$F=T_1F_1+\ldots+T_sF_s$. 

As regards the case $n=1$, notice that $Q\langle 0\rangle =0$ and $\boundary{1}{\underline{z}}{U}_{1}=0$, 
thus a similar construction leads to an isomorphism $\tau_1:Q_1\stackrel{\cong}{\to}\cycle{1}{\underline{z}}{U}_{1}$.  
}\end{Remark}
\demo Take $f=\varphi:V\to U$ and $g={\bf 1}_{V}$ in
\cite[Theorem~2.4]{planas}.

\begin{Remark}\label{min-gen-set-of-equations}{\rm 
Let $U$ be a standard $R$-algebra and let $Q=\oplus_{n\geq 1}Q_{n}$ be
the equations of $U$.  Suppose that $Q=Q\langle N\rangle$ for some
$N\geq 1$. Set 
\begin{eqnarray*}
\mathcal{C}=\{ F_{1,1},\ldots ,F_{1,s_{1}},\ldots ,F_{N,1},\ldots
,F_{N,s_{N}}\}, 
\end{eqnarray*}
with $F_{i,j}\in Q_{i}$. Then $\mathcal{C}$ is a minimal generating
set of $Q$ if and only if, for each $n=1,\ldots ,N$, the classes of
$F_{n,1},\ldots ,F_{n,s_{n}}$ modulo $Q\langle n-1\rangle$ are a
minimal generating set of $(Q/Q\langle n-1\rangle)_{n}$.
}\end{Remark}

\section{Effective relations through Koszul homology}
\label{effective relations of U} 

Let us set a general framework for the study of the injectivity of blowing-up
morphisms by first setting some properties of the module of effective relations 
that we will later specialise and use. We will exploit the fact that
there is an isomorphism between the $n$-th module of effective relations $E(U)_n$ 
and the $n$-th graded component $\koszul{1}{\underline{z}}{U}_n$ (see Section \ref{sec:effective-relations}).


\medskip

Let us first examine a straightforward application of the well-known
long exact sequence of the (graded) Koszul homology, which will be key for 
the subsequent arguments and which we will be able to specialise at low cost.

\begin{Lemma}\label{successio_U}
Let $R$ be a Noetherian ring and let $U=\oplus_{n\geq 0} U_n$ be a 
standard $R$-algebra. Let $\underline{z}=z_1,\ldots,z_s$ and 
$\underline{z}'=\underline{z},z_{s+1}$. Let the sequence $\underline{z}'$ 
be a minimal generating set of $U_1$ as an $R$-module. Then, for each integer 
$n\geq 2$, there is a short exact sequence
\begin{eqnarray*}
0\to \frac{\koszul{1}{\underline{z}}{U}_n}{z_{s+1}\koszul{1}{\underline{z}}{U}_{n-1}} 
\longrightarrow \koszul{1}{\underline{z}'}{U}_n \stackrel{\sigma_n}{\longrightarrow} 
\frac{((\sum_{i=1}^s z_iU_{n-1}):_R z_{s+1}^n)}{((\sum_{i=1}^s z_iU_{n-2}):_R z_{s+1}^{n-1})}\to 0.
\end{eqnarray*}
Moreover, $\sigma_n$ sends the class of a cycle
$(w_1,\ldots,w_s,w_{s+1})$, $w_{i}\in U_{n-1}$,
to the class of $a\in ((\sum_{i=1}^s z_iU_{n-1}):_R z_{s+1}^n)$, 
where $w_{s+1}=az_{s+1}^{n-1}+b$ for some $b\in \sum_{i=1}^s z_iU_{n-2}$.
\end{Lemma}
\demo
Consider the graded long exact sequence of Koszul 
homology (see e.g. \cite[Corollary~1.6.13]{bh}):
\begin{eqnarray*}
 & \koszul{1}{\underline{z}}{U}_{n-1}
  \stackrel{(\rho_1)_n}{\longrightarrow} \koszul{1}{\underline{z}}{U}_n
  \longrightarrow
  \koszul{1}{\underline{z}'}{U}_n\to
  \koszul{0}{\underline{z}}{U}_{n-1}\stackrel{(\rho_0)_n}
         {\longrightarrow} & \koszul{0}{\underline{z}}{U}_n,
\end{eqnarray*}
where $(\rho_i)_n$ is just the multiplication by $\pm z_{s+1}$.
Then we get the following short exact sequence:
\begin{eqnarray*}
0 \to\textrm{coker}(\rho_1)_n \to
\koszul{1}{\underline{z}'}{U}_n
\stackrel{\widetilde{\sigma}_n}{\to} \ker(\rho_0)_n \to 0.
\end{eqnarray*}
Clearly, ${\rm coker}(\rho_1)_n=
\koszul{1}{\underline{z}}{U}_n/z_{s+1}\koszul{1}{\underline{z}}{U}_{n-1}$. 
On the other hand, 
\begin{eqnarray*}
\koszul{0}{\underline{z}}{U}_{n-1}=U_{n-1}/\sum_{i=1}^s z_iU_{n-2}.
\end{eqnarray*}
Thus $\ker(\rho_0)_n=(U_{n-1}\cap ((\sum_i^s z_iU_{n-1}):_U z_{s+1}))/\sum_i^s z_iU_{n-2}$. 
One can check that $\widetilde{\sigma}_n$ maps the homology class of a cycle
$(w_1,\ldots,w_s,w_{s+1})\in\oplus_{i=1}^{s+1}U_{n-1}$ to the class of $w_{s+1}\in
U_{n-1}\cap ((\sum_{i=1}^s z_iU_{n-1}):_U z_{s+1})$ modulo $\sum_{i=1}^s z_iU_{n-2}$. 
Finally, consider the mapping $\theta_n$,
\begin{eqnarray*}
\ker(\rho_0)_n=\frac{U_{n-1}\cap((\sum_{i=1}^s z_iU_{n-1}):_U z_{s+1})}{\sum_{i=1}^s z_iU_{n-2}}
\stackrel{\theta_n}{\longrightarrow}\frac{((\sum_{i=1}^s z_iU_{n-1}):_R z_{s+1}^n)}{((\sum_{i=1}^s z_iU_{n-2}):_R z_{s+1}^{n-1})},
\end{eqnarray*}
defined as follows: for each $w\in U_{n-1}\cap((\sum_{i=1}^s z_iU_{n-1}):_U z_{s+1})$, 
since $U_{n-1}=z_{s+1}^{n-1}R+\sum_{i=1}^s z_iU_{n-2}$, 
take $a\in R$ and $b\in \sum_{i=1}^s z_iU_{n-2}$ such
that $w=az_{s+1}^{n-1}+b$. Clearly $a\in ((\sum_{i=1}^s z_iU_{n-1}):_R z_{s+1}^{n})$. 
Let $\overline{w}$ be the class of $w$ modulo $\sum_{i=1}^s z_iU_{n-2}$ and let
$\overline{a}$ the class of $a$ modulo $((\sum_{i=1}^s z_iU_{n-2}):_R z_{s+1}^{n-1})$. 
Set $\theta_n(\overline{w})=\overline{a}$: an elementary computation shows
that $\theta_n$ is a well-defined isomorphism. Finally, set $\sigma_n=\theta_n\circ\widetilde{\sigma}_n$. 
\qed

\medskip

Lemma \ref{successio_U} can be profitably specialised for either $\rees(I)$, 
$\graded(I)$ or $\fiber(I)$, thus characterising the vanishing of their modules of
effective relations in terms of a mixed condition on colon ideals and Koszul homology. 
In the case where $\mu(I)=2$, the vanishing of the Koszul homology itself 
admits a full description in terms of colon ideals that we make explicit
in the following lemma.

\begin{Lemma}\label{successio_U-two-generated}
Let $R$ be a Noetherian ring and let $U=\oplus_{n\geq 0} U_n$ be a 
standard $R$-algebra. Suppose that $z_1,z_2$ is a minimal generating set
of $U_1$ as an $R$-module. Then, for each integer $n\geq 2$, there is a 
short exact sequence
\begin{eqnarray*}
0\to \frac{(0:_U z_1)\cap U_{n-1}}{z_2((0:_U z_1)\cap U_{n-2})} 
\longrightarrow \koszul{1}{z_1,z_2}{U}_n \stackrel{\sigma_n}{\longrightarrow} 
\frac{(z_1U_{n-1}:_R z_2^n)}{(z_1U_{n-2}:_R z_2^{n-1})}\to 0.
\end{eqnarray*}
\end{Lemma}
\demo
Apply Lemma \ref{successio_U}, observing that $\koszul{1}{z_1}{U}_n\cong (0:_U z_1)\cap U_{n-1}$.
\qed

\medskip

Given ideals $I, L\subset R$, we will denote the standard algebra 
$\rees(I)\otimes_R R/L=\oplus_{n\geq 0}I^n/LI^n$ by $\rees_L(I)$. 
Let $\underline{z}$ denote the sequence of classes modulo $LI$ 
of the elements of a sequence $\underline{x}$ in $I$. If $L=0$, 
then $\underline{z}$ will stand for the sequence $\underline{x}t\subset It$.	 

\begin{Proposition}\label{succexcurta-RL(I)}
Let $(R,\fm)$ be a Noetherian local ring. Let $I=(x_1,\ldots,x_s,y)$,
$J=(x_1,\ldots,x_s)\subset I$ and $L$ be ideals in $R$. 
Let $\underline{x}=x_1,\ldots,x_s$ and assume that 
$\underline{x},y$ minimally generate $I$. Let $\underline{z}$ and 
$\underline{z}'=\underline{z},z_{s+1}$ be the sequence of 
classes modulo $LI$ of the elements of $\underline{x}$ and $\underline{x},y$, respectively.  
Then the following short exact sequence holds:
$$0\to \frac{\koszul{1}{\underline{z}}{\rees_L(I)}_n}{z_{s+1}\koszul{1}{\underline{z}}{\rees_L(I)}_{n-1}} 
\to \koszul{1}{\underline{z}'}{\rees_L(I)}_n \stackrel{\sigma_n}{\to} \frac{(JI^{n-1}+LI^{n}:_R y^n)}{(JI^{n-2}+LI^{n-1}:_R y^{n-1})}\to 0.$$
In particular, if $I=(x,y)$, the preceding short exact sequence has the following form:
$$0\to \frac{(LI^n:_R x)\cap I^{n-1}}{y((LI^{n-1}:_R x)\cap I^{n-2})+LI^{n-1}} \to \koszul{1}{\underline{z}'}{\rees_L(I)}_n \stackrel{\sigma_n}{\to} \frac{(xI^{n-1}+LI^{n}:_R y^n)}{(xI^{n-2}+LI^{n-1}:_R y^{n-1})}\to 0.$$
\end{Proposition}
\demo
One just needs to check that the expressions in the short exact sequences are 
isomorphic to the ones that appear in Lemma \ref{successio_U} and Lemma 
\ref{successio_U-two-generated}.

In order to establish the second short exact sequence, let $z_1,z_2\in I/LI$ be the 
classes of $x,y$, respectively. Let $\epsilon_{n-1}:(LI^n:x)\cap I^{n-1}\to (0:_{R_L(I)} z_1)\cap (I^{n-1}/LI^{n-1})$ 
be the homomorphism taking classes modulo $LI^{n-1}$. Then we have the 
following commutative diagram with exact rows:
\begin{equation*}
\begin{CD}
0 @>>> LI^{n-2} @>>> (LI^{n-1}:_R x)\cap I^{n-2} @>{\epsilon_{n-2}}>> (0 :_{R_L(I)} z_1)\cap \frac{I^{n-2}}{LI^{n-2}} @>>> 0\\
&& @VV{y}V @VV{y}V @VV{z_2}V \\
0 @>>> LI^{n-1} @>>>  (LI^n:_R x)\cap I^{n-1}  @>{\epsilon_{n-1}}>> (0:_{R_L(I)} z_1)\cap \frac{I^{n-1}}{LI^{n-1}} @>>> 0\\
&& @VVV @VVV @VVV \\
&& LI^{n-1}/yLI^{n-2} @>>> \frac{(LI^{n}:_R\; x)\cap I^{n-1}}{y((LI^{n-1}:_R\; x)\cap I^{n-2})} @>>> \frac{(0:_{R_L(I)}\; z_1)\cap I^{n-1}/LI^{n-1}}{z_2((0:_{R_L(I)}\; z_1)\cap I^{n-2}/LI^{n-2})} @>>> 0.
\end{CD}
\end{equation*}
The exactness of the bottom row and Lemma~\ref{successio_U-two-generated} lead to the conclusion.
\qed

\medskip

Although the next lemma is a consequence of Lemma~\ref{successio_U} and 
Proposition~\ref{succexcurta-RL(I)}, we will give an explicit proof
by carrying over the ideas of Lemma~\ref{successio_U}. We keep focusing on
$\koszul{1}{x_1t,\ldots,x_st,yt}{\rees(I)}_n$, which by Remark~\ref{explicit-presentation} 
is isomorphic to the module of effective relations of $\rees(I)$, 
$E(I)_n\cong (Q/Q\langle n-1 \rangle)_n$ (see also 
Definition~\ref{definicio-relacions-efectives}).

\begin{Lemma}\label{succexcurta}
Let $(R,\fm)$ be a Noetherian local ring and let $I$ be an ideal of
$R$. Let $x_1,\ldots,x_s,y$ be a minimal generating set of $I$. Set
$J=(x_1,\ldots,x_s)$. Then, for each integer $n\geq 2$, there is a
short exact sequence
\begin{eqnarray*}
0\to \frac{\koszul{1}{x_1t,\ldots,x_st}
    {\rees(I)}_n}{yt\koszul{1}{x_1t,\ldots,x_st}{\rees(I)}_{n-1}}\longrightarrow
  \koszul{1}{x_1t,\ldots,x_st,yt}{\rees(I)}_n\stackrel{\sigma_n}
         {\longrightarrow} \frac{(JI^{n-1}:_R y^n)}{(JI^{n-2}:_R y^{n-1})} \to
         0.
\end{eqnarray*}
Moreover, $\sigma_n$ sends the class of a cycle
$(w_1t^{n-1},\ldots,w_st^{n-1},w_{s+1}t^{n-1})$, $w_{i}\in I^{n-1}$,
to the class of $a\in (JI^{n-1}:_R y^{n})$, where $w_{s+1}=ay^{n-1}+b$
for some $b\in JI^{n-2}$.
\end{Lemma}

\demo Take $\underline{z}=x_1t,\ldots,x_st$ and $\underline{z}'=x_1t,\ldots,x_st,yt$ in
$It\subset \rees(I)=R[It]=\oplus_{n\geq 0}I^n$. Consider the induced
graded long exact sequence of Koszul homology:
\begin{eqnarray*}
\koszul{1}{\underline{z}}{\rees(I)}_{n-1}
  \stackrel{(\rho_1)_n}{\longrightarrow} \koszul{1}{\underline{z}}{\rees(I)}_n
  \longrightarrow \koszul{1}{\underline{z}'}{\rees(I)}_n\to
  \koszul{0}{\underline{z}}{\rees(I)}_{n-1}\stackrel{(\rho_0)_n}
         {\longrightarrow} \koszul{0}{\underline{z}}{\rees(I)}_n,
\end{eqnarray*}
where $(\rho_i)_n$ is just the multiplication by $\pm yt$. We get the
following short exact sequence:
\begin{eqnarray*}
0 \to\textrm{coker}(\rho_1)_n \to
\koszul{1}{\underline{z}'}{\rees(I)}_n
\stackrel{\widetilde{\sigma}_n}{\to} \ker(\rho_0)_n \to 0.
\end{eqnarray*}
Clearly ${\rm coker}(\rho_1)_n=
\koszul{1}{\underline{z}}{\rees(I)}_n/yt\koszul{1}{\underline{z}}{\rees(I)}_{n-1}$. On the
other hand, 
\begin{eqnarray*}
\koszul{0}{\underline{z}}{\rees(I)}_{n-1}=I^{n-1}/JI^{n-2}.
\end{eqnarray*}
 Thus $\ker(\rho_0)_n=I^{n-1}\cap (JI^{n-1}:_R y)/JI^{n-2}$. One can
 check that $\widetilde{\sigma}_n$ maps the homology class of a cycle
 $(w_1t^{n-1},\ldots,w_st^{n-1},w_{s+1}t^{n-1})\in
 \oplus_{i=1}^{s}I^{n-1}t^{n-1}$ to the class of $w_{s+1}\in
 I^{n-1}\cap (JI^{n-1}:_R y)$ modulo $JI^{n-2}$. Finally, consider the
 mapping $\theta_n$,
\begin{eqnarray*}
\ker(\rho_0)_n=\frac{I^{n-1}\cap(JI^{n-1}:_R y)}{JI^{n-2}}
\stackrel{\theta_n}{\longrightarrow}\frac{(JI^{n-1}:_R y^n)}{(JI^{n-2}:_R y^{n-1})},
\end{eqnarray*}
defined as follows: for each $w\in I^{n-1}\cap(JI^{n-1}:_R y)$, since
$I^{n-1}=y^{n-1}R+JI^{n-2}$, take $a\in R$ and $b\in JI^{n-2}$ such
that $w=ay^{n-1}+b$. Clearly $a\in (JI^{n-1}:_R y^{n})$. Let
$\overline{w}$ be the class of $w$ modulo $JI^{n-2}$ and let
$\overline{a}$ the class of $a$ modulo $(JI^{n-2}:_R y^{n-1})$. Set
$\theta_n(\overline{w})=\overline{a}$. An elementary computation shows
that $\theta_n$ is a well-defined isomorphism. Take $\sigma_n=
\theta_n\circ\widetilde{\sigma}_n$.\qed

\medskip

Although the Lemma~\ref{succexcurta} has a version for $n=1$, we are
more interested in an analogous version obtained from considering
$x_{1},\ldots ,x_{s}$ and $x_1,\ldots ,x_s,y$ as sequences of elements
of degree zero in $I\subset R$.

\begin{Remark}\label{syzygies}{\rm
With the assumptions and notations of Lemma~\ref{succexcurta}, there
is an exact sequence
\begin{eqnarray*}
  0\to\frac{\koszul{1}{x_1,\ldots,x_s}{R}}
  {y\koszul{1}{x_1,\ldots,x_s}{R}}\to
  \koszul{1}{x_1,\ldots,x_s,y}{R}\stackrel{\tilde{\sigma}_1}{\to}
  \frac{(J:_R y)}{J}\to 0,
\end{eqnarray*}
where $\tilde{\sigma}_1$ sends the homology class of a cycle
$(w_1,\ldots,w_s,w_{s+1})\in R^{s+1}$ to the class of $w_{s+1}\in
(J:_R y)$. In particular, if $x_1,\ldots,x_s$ is an $R$-sequence,
$\koszul{1}{x_1,\ldots,x_s}{R}=0$ and $\tilde{\sigma}_1$ is an
isomorphism. 
}\end{Remark}

\section{Vanishing of the Koszul homology}

If $x_1,\ldots,x_s$ a sequence of elements belonging to an ideal 
$I$, the next lemma characterises the vanishing of the $n$-th graded 
components $\koszul{1}{x_1t,\ldots,x_it}{\rees(I)}_{n}$ for all
$1\leq i\leq s$, in terms of sequential conditions in the 
$x_1,\ldots,x_s$.

\begin{Lemma}\label{characterisation-of-homology}
Let $(R,\fm)$ be a Noetherian local ring and let $I$ be an ideal of
$R$. Let $x_1,\ldots,x_s$ be a sequence of elements belonging to $I$
and let $n\geq 2$ be an integer. Then, the following conditions are 
equivalent:
\begin{enumerate}
\item[$(i)$] $\koszul{1}{x_1t,\ldots,x_it}{\rees(I)}_{n}=0$, for all $i=1,\ldots,s$;
\item[$(ii)$] $I^{n-1}t^{n-1}\cap ((x_1t,\ldots,x_{i-1}t):_{R(I)} x_it)=I^{n-1}t^{n-1}\cap (x_1t,\ldots,x_{i-1}t)$, 
for all $i=1,\ldots,s$;
\item[$(iii)$] $((x_1,\ldots,x_{i-1})I^{n-1}:_R x_i)\cap I^{n-1}=(x_1,\ldots,x_{i-1})I^{n-2}$, 
for all $i=1,\ldots,s$.
\end{enumerate}
\end{Lemma}
\demo
Let $n$ be fixed and let us prove the equivalence by induction on $s$. If $s=1$ then 
$\koszul{1}{x_1t}{\rees(I)}_n=(0:_R x_1)\cap I^{n-1}$ vanishes if and only if 
$(0:_R x_1)\cap I^{n-1}=0$, hence the claim follows. 

By induction, assume the lemma holds for $s-1>0$. Let us denote 
$J_{0}=0$ and $J_{i}=(x_{1},\ldots ,x_{i})$.
Take $1\leq i\leq s$ and set $\underline{z}=x_1t,\ldots,x_{i-1}t$ and 
$\underline{z}'=x_1t,\ldots,x_{i-1}t,x_it$. By the induction 
hypothesis and using the graded long exact sequence 
of Koszul homology,
\begin{eqnarray*}
0\to\koszul{1}{\underline{z}'}{\rees(I)}_n\to
\koszul{0}{\underline{z}}{\rees(I)}_{n-1}\stackrel{(\rho_0)_n}
       {\to}\koszul{0}{\underline{z}}{\rees(I)}_n
\end{eqnarray*}
is an exact sequence for all $1\leq i\leq s$. 
Therefore,
\begin{align*}
\koszul{1}{x_1t,\ldots,x_it}{\rees(I)}_{n} & \cong \ker(\rho_{0})_{n} = 
\frac{I^{n-1}t^{n-1}\cap ((x_1t,\ldots,x_{i-1}t):_{R(I)} x_{i}t)}{I^{n-1}t^{n-1}\cap (x_1t,\ldots,x_{i-1}t)} \\
& \cong \frac{(J_{i-1}I^{n-1}:_R x_{i})\cap I^{n-1}}{J_{i-1}I^{n-2}},
\end{align*}
whence the claim follows for all $1\leq i\leq s$.
\qed

\medskip

Recall that the sequence $x_{1},\ldots ,x_{s}$ is a $d$-sequence if 
it minimally generates $(x_1,\ldots,x_s)$ and if
$(J_i:_R x_{i+1}x_{j})=(J_i:_R x_{j})$ for all $0\leq i\leq s-1$ and all
$j\geq i+1$ (where $J_0=0$, $J_{i}=(x_1,\ldots ,x_i)$ and $J_s=J$, 
as before). Recall that this second condition is equivalent to 
$(J_i:_R x_{i+1})\cap J= J_i$ for all $0\leq i\leq s-1$ (see 
\ref{def:d-sequence} and subsequent remarks). Clearly, any $R$-sequence
is a $d$-sequence.

\begin{Remark}{\rm
If the equivalent conditions in Lemma \ref{characterisation-of-homology} 
hold for all $n\geq 2$, then $x_1t,\ldots,x_st$ is a $d$-sequence of 
$\rees(I)$. Furthermore, $((x_1t,\ldots,x_{i-1}t):_{R(I)} x_it)\cap\rees(I)_{+}=
(x_1t,\ldots,x_{i-1}t)$ for all $i=1,\ldots,s$.
 
However, being a $d$-sequence in $\rees(I)$ does not imply 
that $\koszul{1}{x_1t,\ldots,x_it}{\rees(I)}_n=0$ for all $i=1,\ldots,s$ 
and $n\geq 2$. Consider the ring $R=k[X,Y]_{(X,Y)}/(XY,Y^2)_{(X,Y)}$ with 
$k$ a field and let $I=(x,y)$ and $J=(x)$. Then $xt$ is a $d$-sequence in 
$\rees(I)$, while $(0:x)\cap (x,y)=(y)\neq 0$.
}\end{Remark}
\demo
Let $x_1,\ldots,x_s$ and $I$ as in Lemma \ref{characterisation-of-homology}. 
Let $\widetilde{J}=(x_1t,\ldots,x_st)\subset \rees(I)_{+}$. The condition
$I^{n-1}t^{n-1}\cap((x_1t,\ldots,x_{i-1}):_{R(I)} x_it)=I^{n-1}t^{n-1}\cap(x_1t,\ldots,x_{i-1})$, 
for all $i=1,\ldots,s$ and all $n\geq 2$, is equivalent 
to $((x_1t,\ldots,x_{i-1}t):_{R(I)} x_it)\cap\rees(I)_{+}=
(x_1t,\ldots,x_{i-1}t)$ for all $i=1,\ldots,s$, which implies 
$((x_1t,\ldots,x_{i-1}t):_{R(I)} x_it)\cap\widetilde{J}=
(x_1t,\ldots,x_{i-1}t)$ for all $i=1,\ldots,s$, i.e., $x_1t,\ldots,x_st$ 
is a $d$-sequence. However, the converse is not true, even if $J$ is a 
reduction of $I$, as it is shown in Example~\ref{not-enough}. Compare this 
with \cite[Lemma 12.7]{hsv}.
\qed

\medskip
In the next lemma we will prove that, for $n\geq 2$, then
$\koszul{1}{x_1t,\ldots,x_st}{\rees(I)}_{n}=0$, provided that
$x_{1},\ldots ,x_{s}$ an $R$-sequence such that $x_{1}^*,\ldots
,x_{s-1}^*$ is a $\graded(I)$-sequence. Notice that we are not saying
that the whole $\koszul{1}{x_1t,\ldots,x_st}{\rees(I)}$ vanishes. In
fact, if $s>1$, $x_{1}(x_{2}t)\in x_{1}t\cdot\rees(I)$ whereas
$x_{1}\not\in x_{1}t\cdot\rees(I)$, hence $x_1t,\ldots,x_st$ is not an
$\rees(I)$-sequence.

Before presenting the lemma, recall the result of Valabrega and Valla in
\cite[Corollary~2.7]{vv} which characterizes being a
$\graded(I)$-sequence.  Let $J\subset I$ be two ideals of $R$ and let
$x_1,\ldots ,x_s$ denote a minimal generating set of $J$. Keeping the
hypotheses and notations in the proof of Lemma \ref{characterisation-of-homology} 
(in particular, $R$ is local so the over-riding condition in \cite[Section 2]{vv} is met), 
write $J_{0}=0$ and $J_{i}=(x_{1},\ldots ,x_{i})$, for $i=1,\ldots ,s$. In particular, 
the initial forms $x_1^*,\ldots ,x_s^*$ in $\graded(I)$ are in $I/I^2$. Then, 
$x_1^*,\ldots ,x_s^*$ is a $\graded(I)$-sequence if and only if 
$x_1,\ldots ,x_s$ is an $R$-sequence and the Valabrega-Valla 
modules $VV_{J_i}(I)_{n}=J_i\cap I^n/J_iI^{n-1}$ vanish
for all $i=1,\ldots ,s$ and all $n\geq 1$.

\begin{Lemma}\label{homologia-zero}
Let $(R,\fm)$ be a Noetherian local ring and let $I$ be an ideal of
$R$. Let $x_1,\ldots,x_s$ be a minimal generating set of $J$, where
$J=(x_1,\ldots,x_s)$ is a reduction of $I$ with reduction number
$r=\rednumber_J(I)$. Assume that $x_1,\ldots,x_s$ is an $R$-sequence
and that $x_1^*,\ldots,x_{s-1}^*$ is a $\graded(I)$-sequence. Then,
for all $n\geq 2$ and all $i=1,\ldots,s$,
\begin{eqnarray*}
\koszul{1}{x_1t,\ldots,x_it}{\rees(I)}_n=0.
\end{eqnarray*}
\end{Lemma}
\demo We will check that condition $(iii)$ of Lemma 
\ref{characterisation-of-homology} is fulfilled for all $n\geq 2$. 
Using the aforementioned result of Valabrega and Valla, 
the hypotheses of the Lemma imply
\begin{eqnarray*}
VV_{J_i}(I)_{n}=J_i\cap I^n/J_iI^{n-1}=0
\end{eqnarray*}
for all $i=1,\ldots ,s-1$ and all $n\geq 2$. Since $x_1,\ldots,x_s$
is an $R$-sequence, then $(J_{i-1}:_R x_{i})=J_{i-1}$ for all
$i=1,\ldots,s$. Therefore 
\begin{eqnarray*}
(J_{i-1}:_R x_{i})\cap I^{n-1}=J_{i-1}\cap I^{n-1}=J_{i-1}I^{n-2},
\end{eqnarray*}
for all $i=1,\ldots, s$ and all $n\geq 2$. 
This identity implies
\begin{eqnarray}\label{cond-ast}
(J_{i-1}I^{n-1}:_R x_{i})\cap I^{n-1}=J_{i-1}I^{n-2},
\end{eqnarray}
for all $i=1,\ldots ,s$ and $n\geq 2$.
But these are precisely the conditions $(iii)$ in 
Lemma \ref{characterisation-of-homology}, 
whence the claim.\qed

Notice that the assumption $J\subseteq I$ being a reduction is not needed 
in the proof of Lemma~\ref{homologia-zero}. However, we will carry this 
condition since it will be relevant in the context of the results sought. 

One can state a different version of Lemma~\ref{homologia-zero} by considering weaker hypotheses 
and hence a weaker thesis, which will in turn lead to a slightly different version 
of the main result of this chapter (see Remark~\ref{different-version-theorem}).

\begin{Remark}\label{different-version-lemma}{\rm
Let $(R,\fm)$ be a Noetherian local ring and let $I$ be an ideal of
$R$. Let $x_1,\ldots,x_s$ be a minimal generating set of $J$, where
$J=(x_1,\ldots,x_s)$ is a reduction of $I$ with reduction number
$r=\rednumber_J(I)$. Assume that
\begin{itemize}
\item[$(a)$] $x_1,\ldots,x_s$ is a $d$-sequence;
\item[$(b)$] $VV_{J_i}(I)_{r+1}=(x_{1},\ldots ,x_{i})\cap
  I^{r+1}/(x_{1},\ldots ,x_{i})I^{r}=0$ for all $i=1,\ldots ,s-1$.
\end{itemize}
Then $\koszul{1}{x_1t,\ldots,x_it}{\rees(I)}_n=0$, for all $n\geq r+2$
and all $i=1,\ldots ,s$. 

Suppose that, in addition,
\begin{itemize}
\item[$(c)$] $x_1,\ldots,x_s$ is an $R$-sequence;
\item[$(d)$] $VV_{J_{i}}(I)_{r}=(x_{1},\ldots ,x_{i})\cap
  I^{r}/(x_{1},\ldots ,x_{i})I^{r-1}=0$ for all $i=1,\ldots ,s-1$.
\end{itemize}
Then $\koszul{1}{x_1t,\ldots,x_it}{\rees(I)}_{r+1}=0$, for all
$i=1,\ldots ,s$.
}\end{Remark}
\demo
Using $(a)$ and $(b)$, since $I^{r+1}=JI^r \subset J$ we get 
$$(J_{i-1}:_R x_{i})\cap I^{r+1}=(J_{i-1}:_R x_{i})\cap J\cap I^{r+1}=
J_{i-1}\cap I^{r+1}=J_{i-1}I^{r},$$ 
for $i=1,\ldots,s$. Therefore, for all $i=1,\ldots ,s$,
\begin{eqnarray*}
(J_{i-1}:_R x_{i})\cap I^{r+1}=J_{i-1}I^{r}.
\end{eqnarray*}
Using the result of Trung in \cite[Proposition~4.7(i)]{trung2},
one has $(J_{i-1}:_R x_{i})\cap I^{n}=J_{i-1}I^{n-1}$, for all $n\geq
r+1$ and all $i=1,\ldots ,s$. These identities imply that
\begin{eqnarray*}
(J_{i-1}I^{n-1}:_R x_{i})\cap I^{n-1}=J_{i-1}I^{n-2},
\end{eqnarray*}
for all $n\geq r+2$ and all $i=1,\ldots ,s$. Therefore, by 
Lemma \ref{characterisation-of-homology}, the claim holds. 
As for the second claim, $(c)$ implies that 
$(J_{i-1}:_R x_{i})\cap I^r=J_{i-1}\cap I^{r}$ for all $i=1,\ldots,s$
(observe that that $(a)$ is not sufficient) and using $(d)$ we get
$(J_{i-1}:_R x_{i})\cap I^r=J_{i-1}I^{r-1}$, hence 
$(J_{i-1}I^r:_R x_{i})\cap I^r=J_{i-1}I^{r-1}$ for all $i=1,\ldots,s$,
thus Lemma \ref{characterisation-of-homology} shows that 
$\koszul{1}{x_1t,\ldots,x_it}{\rees(I)}_n=0$ for $i=1,\ldots,s$ 
and $n\geq r+1$.\qed

\begin{Remark}\label{rednumber=1}{\rm
Observe that if $J\subset I$ is a reduction
with small reduction number $\rednumber_J(I)$, then 
the graded homology vanishes with milder hypotheses 
than those in Lemma~\ref{homologia-zero}.
Let $(x_1,\ldots,x_s)=J\subset I$ be a reduction with 
$\rednumber_J(I)=1$. Suppose that:
\begin{enumerate}
\item[$(a')$] $x_1,\ldots,x_s$ is a $d$-sequence;
\item[$(b')$] $((x_1,\ldots,x_{i-1}):_R x_i)\cap I = (x_1,\ldots,x_{i-1})$, for all $i=1,\dots,s$;
\item[$(c')$] $(x_1,\ldots,x_{i-1})\cap I^2=(x_1,\ldots,x_{i-1})I$, for all $i=1,\ldots,s-1$.
\end{enumerate}
Then, $\koszul{1}{x_1t,\ldots,x_st}{\rees(I)}_n=0$, for $n\geq 2$.
}\end{Remark}
\demo
As in the proof of Lemma \ref{homologia-zero} and Remark \ref{different-version-lemma}, 
the goal is to reach the identities in Lemma \ref{characterisation-of-homology} which are
equivalent to the vanishing of the graded homology. From $(b')$ we are able to prove $\koszul{1}{x_1t,\ldots,x_st}{\rees(I)}_2=0$. On the other hand, $(a')$, $(c')$ and
\cite[Proposition~4.7 (i)]{trung2} lead to $\koszul{1}{x_1t,\ldots,x_st}{\rees(I)}_n=0$, 
for $n\geq 3$.
\qed

\begin{Remark}{\rm
Let $x_1,\ldots, x_s$ be a $d$-sequence contained in $I$ and $J=(x_1,\ldots,x_s)$. 
We know that $J$ is of linear type (see \ref{d-sequences-are-linear-type}), hence 
$\koszul{1}{x_1t,\ldots,x_st}{\rees(J)}_n=0$ for all $n\geq 2$ (see \ref{effective-relations-homological}), 
i.e., the complex 
$$\rees(J)(-2)^{\binom{s}{2}}\stackrel{\partial_2}{\longrightarrow}
\rees(J)(-1)^{s}\stackrel{\partial_1}{\longrightarrow} \rees(J)$$
is exact in degree $\geq 2$. Tensoring by $\rees(I)$ and computing $\ker\partial_1\otimes {\bf 1}_{R(I)}/{\rm Im}\, \partial_2\otimes {\bf 1}_{R(I)}$ yields $\koszul{1}{x_1t,\ldots,x_st}{\rees(I)}$, but notice that we do not necessarily get
$\koszul{1}{x_1t,\ldots,x_st}{\rees(I)}_n=0$ for all $n\geq 2$ 
(see Example \ref{not-enough}). For instance, $\koszul{1}{x_1t,\ldots,x_st}{\rees(I)}_2$ 
is the homology of the complex
\begin{equation*}
R^{\binom{s}{2}}\stackrel{\widetilde{\partial}_{2}}{\longrightarrow} 
J^{\oplus s}\oplus I^{\oplus s} \stackrel{\widetilde{\partial}_{1}}{\longrightarrow}
J^2 \oplus JI
\end{equation*}
where $\widetilde{\partial}_{2}(a_0)=(\partial_{2}(a_0),0)$ and 
$\widetilde{\partial}_{1}(a_1, b_0)=(\partial_{1}(a_1),\partial_{1}(b_0))$. 
It is readily seen that $\ker\widetilde{\partial}_{1}$ may contain 
elements which are not in ${\rm Im}\,\widetilde{\partial}_{2}$.
}\end{Remark}

\section{Theorem A}\label{main-result}

We have now all the ingredients to prove the main result of the
chapter. As in the rest of the chapter, set $V=R[X_1,\ldots,X_s,Y]$ 
and let $Q$ be the kernel of the polynomial presentation 
$\varphi:V\to\rees(I)$ sending $X_i$ to $x_it$ and $Y$ to $yt$.

\begin{Theorem A}
Let $(R,\fm)$ be a Noetherian local ring and let $I$ be an ideal of
$R$. Let $x_1,\ldots,x_s,y$ be a minimal generating set of $I$, where
$J=(x_1,\ldots,x_s)$ is a reduction of $I$ with reduction number
$r=\rednumber_J(I)$. Assume that $x_1,\ldots,x_s$ verify 
the following condition for all $n\geq 2$:
\[ 
((x_1,\ldots,x_{i-1})I^{n-1}: x_i)\cap I^{n-1}=(x_1,\ldots,x_{i-1})I^{n-2}, 
\textrm{ for all } i=1,\ldots,s. \label{eq:trung-condition} \tag{$\mathcal{T}_n$}
\]
Then, for each $n\geq 2$, the map sending $F\in Q_{n}$ to
$F(0,\ldots,0,1)\in (JI^{n-1}:_R y^{n})$ induces an isomorphism of
$R$-modules
\begin{eqnarray*}
\left[ \frac{Q}{Q\langle n-1\rangle}\right]_{n}\cong
\frac{(JI^{n-1}:_R y^{n})}{(JI^{n-2}:_R y^{n-1})}.
\end{eqnarray*}
In particular, $\reltype(I)=\rednumber_J(I)+1$ and there is a form
$Y^{r+1}-\sum X_iF_i\in Q_{r+1}$, with $F_i\in V_r$, such that
$Q=(Y^{r+1}-\sum X_iF_i)+Q\langle r\rangle$. Moreover, if 
$x_1,\ldots,x_s$ is an $R$-sequence and $x_1^*,\ldots,x_{s-1}^*$ is
a $\graded(I)$-sequence, then $x_1,\ldots,x_s$ verify condition $(\mathcal{T}_n)$
for all $n\geq 2$.
\end{Theorem A}

\demo 
By Lemma~\ref{succexcurta}, with $\underline{z}=x_1t,\ldots,x_st$ and
$\underline{z}'=x_1t,\ldots,x_st,yt$,
\begin{eqnarray*}
0\to\frac{\koszul{1}{\underline{z}}{\rees(I)}_{n}}{yt\koszul{1}{\underline{z}}{\rees(I)}_{n-1}}\to
\koszul{1}{\underline{z}'}{\rees(I)}_n\stackrel{\sigma_n}{\to}
\frac{(JI^{n-1}:_R y^{n})}{(JI^{n-2}:_R y^{n-1})} \to 0
\end{eqnarray*}
is an exact sequence for all $n\geq 2$. By Lemma~\ref{characterisation-of-homology},
$\koszul{1}{\underline{z}}{\rees(I)}_{n}=0$ for all $n\geq 2$. Therefore,
$\koszul{1}{\underline{z}'}{\rees(I)}_n\stackrel{\sigma_n}{\cong}
(JI^{n-1}:_R y^{n})/(JI^{n-2}:_R y^{n-1})$, for $n\geq 2$. By
Remark~\ref{explicit-presentation}, we conclude:
\begin{eqnarray*}
\left[ \frac{Q}{Q\langle n-1\rangle}\right]_{n}
\stackrel{\tau_n}{\cong}\koszul{1}{\underline{z}'}{\rees(I)}_n
\stackrel{\sigma_n}{\cong} \frac{(JI^{n-1}:_R y^{n})}{(JI^{n-2}:_R y^{n-1})},
\end{eqnarray*}
for all $n\geq 2$. Given $F\in Q_{n}$, write
$F=\sum_{i=1}^{s}X_{i}F_{i}+YG$, with $F_{i},G\in V_{n-1}$. Then the
morphism $\tau_n$ sends the class of $F$ to the homology class of
$(F_{1}(\underline{z}'),\ldots ,F_{s}(\underline{z'}),
G(\underline{z}'))$. But $G(\underline{z}')=G(x_{1},\ldots
,x_{s},y)t^{n-1}$ and $G(x_{1},\ldots ,x_{s},y)=G(0,\ldots
,0,1)y^{n-1}+b$, for some $b\in JI^{n-2}$. By Lemma~\ref{succexcurta},
$\sigma_n$ sends the homology class of $(F_{1}(\underline{z}'),\ldots
,F_{s}(\underline{z}'), G(\underline{z}'))$ to the class of
$G(0,\ldots,0,1)$ modulo $(JI^{n-2}:_R y^{n-1})$ and notice that
$G(0,\ldots,0,1)=F(0,\ldots,0,1)$.

Since $J$ is a reduction of $I$ with reduction number $\rednumber_J(I)=r$, 
we know that $(JI^{n-1}:_R y^{n})=R$ for all $n\geq r+1$. Therefore $(Q/Q\langle
n-1\rangle)_n=0$ for all $n>r+1$ and $(Q/Q\langle r\rangle)_{r+1}\cong
R/(JI^{r-1}:_R y^r)\neq 0$, since $y^r\not\in JI^{r-1}$. Therefore
$\reltype(I)=r+1$. Notice that the containment $y^{r+1}\in
JI^{r}$ induces an equation of the form $Y^{r+1}-\sum_iX_iF_i$, with
$F_i\in V_r$, which is sent by $\sigma_n\circ\tau_n$ to the class of
$1$ in $R/(JI^{r-1}:_R y^r)$.

Finally, the last claim is nothing but the application of 
Lemma~\ref{homologia-zero}, together with Lemma~\ref{characterisation-of-homology}.
\qed

\medskip

With weaker hypotheses we get the next version of
Theorem~A (see Remark~\ref{different-version-lemma}).

\begin{Remark}\label{different-version-theorem}{\rm
Let $(R,\fm)$ be a Noetherian local ring and let $I$ be an ideal of
$R$. Let $x_1,\ldots,x_s,y$ be a minimal generating set of $I$, where
$J=(x_1,\ldots,x_s)$ is a reduction of $I$ with reduction number
$r=\rednumber_J(I)$. Assume that the following two conditions hold:
\begin{itemize}
\item[$(a)$] $x_1,\ldots,x_s$ is a $d$-sequence;
\item[$(b)$] $VV_{J_i}(I)_{r+1}=(x_{1},\ldots ,x_{i})\cap
  I^{r+1}/(x_{1},\ldots ,x_{i})I^{r}=0$ for all $i=1,\ldots ,s-1$.
\end{itemize}
Then $\reltype(I)=\rednumber_J(I)+1$ and 
$(Q/Q\langle r\rangle)_{r+1}$ contains the non-zero class 
of an equation $Y^{r+1}-\sum X_iF_i$ induced by 
$y^{r+1}\in JI^r$.

\medskip

Suppose, in addition, that the following conditions hold:
\begin{itemize}
\item[$(c)$] $x_1,\ldots,x_s$ is an $R$-sequence;
\item[$(d)$] $VV_{J_{i}}(I)_{r}=(x_{1},\ldots ,x_{i})\cap
  I^{r}/(x_{1},\ldots ,x_{i})I^{r-1}=0$ for all $i=1,\ldots ,s-1$.
\end{itemize}
Then $\reltype(I)=\rednumber_J(I)+1$ and there is a form $Y^{r+1}-\sum
X_iF_i\in Q_{r+1}$ induced by $y^{r+1}\in JI^r$, with $F_i\in V_r$, 
such that $Q=(Y^{r+1}-\sum X_iF_i)+Q\langle r\rangle$.  
}\end{Remark}

\demo 
It follows from the proof of Theorem~A, but using
Remark~\ref{different-version-lemma}.
\qed

\begin{Discussion}\label{general-hypotheses}{\rm 
The hypotheses of Remark~\ref{different-version-theorem} are connected
with the works of Huckaba in \cite[Theorem~1.4]{huckaba2} and Trung in
\cite[Theorem~6.4]{trung2} (see also \cite[Theorem~3.2]{cz},
\cite{huckaba1}, \cite[Theorem~5.3]{gp} and \cite{trung1}). In
\cite[Theorems~1.4, 1.5]{huckaba2}, Huckaba proved that if $I$ is an
ideal with $\spread(I)=\height(I)+1\geq 2$ and such that any minimal
reduction $J$ of $I$ can be generated by a $d$-sequence $x_{1},\ldots
,x_{s}$ with $x_{1}^{*},\ldots ,x_{s-1}^{*}$ being a
$\graded(I)$-sequence ($s=\spread(I)$), then $\reltype(I)\leq
\rednumber_J(I)+1$. If in addition $\mu(I)=\spread(I)+1$, then the
equality $\reltype(I)= \rednumber_J(I)+1$ holds.  In particular,
$r=\rednumber_J(I)$ is independent of $J$. In fact, Trung improved
this last result in \cite[Theorem~6.4]{trung2} by showing that $r$
coincides with the Castelnuovo-Mumford regularity of $\rees(I)$. To
prove $\reltype(I)\geq\rednumber_J(I)+1$, Huckaba showed that the
equality $I^{r+1}=JI^r$ induces an equation of $\rees(I)$ of maximum
degree. In Theorem~A and with a different approach,
we have fulfilled the description of the entire ideal of equations of
$\rees(I)$.}
\end{Discussion}

As a corollary of Theorem~A, one can prove, following
an alternative path, the result of Vasconcelos in 
\cite[Theorem~2.3.3]{vasconcelos1}. The technique used 
here allows us to prove it for any reduction number, 
not necessarily equal to $1$.

\begin{Corollary}\label{crelle-generalisation}
Let $(R,\fm)$ be a Cohen-Macaulay local ring of dimension $d$ and let
$I$ be an $\fm$-primary ideal of $R$. Let $x_{1},\ldots, x_{d},y$ be a
minimal generating set of $I$, with $J=(x_{1},\ldots ,x_{d})\subset I$ 
being a reduction of $I$, with reduction number $\rednumber_J(I)=r$. Suppose
that $x_{1}^{*},\ldots ,x_{d-1}^{*}$ is a $\graded(I)$-sequence. Then
there is a form $Y^{r+1}-\sum X_iF_i\in Q_{r+1}$, induced by $y^{r+1}\in JI^r$, 
with $F_i\in V_r$, such that $Q=(Y^{r+1}-\sum X_iF_i)+Q\langle r\rangle$.  
In particular, $\reltype(I)=r+1$.
\end{Corollary}
\demo
Since $R$ is Cohen-Macaulay and $I$ is $\fm$-primary, 
$x_{1},\ldots ,x_{d}$ is an $R$-sequence and the results 
follows from Theorem~A.
\qed

\medskip

The case $\rednumber_J(I)=1$ has attracted much attention in the literature.
In fact, this was the case covered in \cite[Theorem~2.3.3]{vasconcelos1}.
Using Remark \ref{rednumber=1}, we are able to state an analogue of 
Theorem~A with $\rednumber_J(I)=1$, which permits us to relax other 
hypotheses.

\begin{Proposition}\label{rednumber=1-theorem}
Let $(R,\fm)$ be a Noetherian local ring and let $I$ be an ideal of
$R$. Let $x_1,\ldots,x_s,y$ be a minimal generating set of $I$, where
$J=(x_1,\ldots,x_s)\subset I$ is a reduction of $I$ with reduction number
$\rednumber_J(I)=1$. Suppose that
\begin{enumerate}
\item[$(a')$] $x_1,\ldots,x_s$ is a $d$-sequence;
\item[$(b')$] $((x_1,\ldots,x_{i-1}):_R x_i)\cap I = (x_1,\ldots,x_{i-1})$, for all $i=1,\dots,s$;
\item[$(c')$] $(x_1,\ldots,x_{i-1})\cap I^2=(x_1,\ldots,x_{i-1})I$, for all $i=1,\ldots,s-1$.
\end{enumerate}
Then, $\reltype(I)=2$ and we have the following isomorphisms of $R$-modules: 
\begin{eqnarray*}
\ker \alpha_2 \cong \left[ \frac{Q}{Q\langle 1 \rangle}\right]_{2}\cong
\frac{(JI:_R y^2)}{(J:_R y)}=\frac{R}{(J:_R y)}.
\end{eqnarray*}
In particular, there is a form $Y^2-\sum X_iF_i\in Q_{2}$, 
with $F_i\in V_1$, such that $$Q=(Y^2-\sum X_iF_i)+Q\langle 1\rangle.$$ 
\end{Proposition}
\demo
It follows from the proof of the Theorem~A, using
Remark~\ref{rednumber=1}.
\qed

\medskip

Another corollary of Theorem~A is the following
result due to Heinzer and Kim in \cite[Theorem~5.6]{hk}. Recall that for
an ideal $L$ of $R$ and any standard $R$-algebra $U$, it is verified that 
$\reltype(U\otimes R/L)\leq \reltype(U)$ (see e.g. \cite[Example~3.2]{planas}),
hence $\reltype(\fiber(I))\leq \reltype(\graded(I))\leq \reltype(I)$. 
In fact, for any $n\geq 2$, there is an exact sequence
$E(I)_{n+1}\to E(I)_{n}\to E(\graded(I))_{n}\to 0$. In particular,
$\reltype(I)=\reltype(\graded(I))$ and $E(I)_{N}\cong
E(\graded(I))_{N}$, with $N=\reltype(I)$
(see \cite[Proposition~3.3]{planas}; see also \cite[p.~268]{hku}). Within 
our settings, however, $\reltype(\fiber(I))$ is also equal to $\reltype(I)$,
as proved in the next corollary. Notice that this behaviour is by no means 
a general fact (see the Conjecture of Valla in \cite[\S~2]{hmv}
and a counterexample in \cite[Example~4.4]{suv}; see also \cite[p.~268
and Corollary~2.6]{hku}).

\begin{Corollary}\label{hk-generalisation}
Let $(R,\fm)$ be a Noetherian local ring with infinite residue field
$k=R/\fm$ and let $I$ be an ideal of $R$. Let $x_1,\ldots,x_s,y$ be a
minimal generating set of $I$, where $J=(x_1,\ldots,x_s)\subset I$ is a
reduction of $I$ with reduction number $r=\rednumber_J(I)$. Assume
that $x_1,\ldots,x_s$ is an $R$-sequence and that
$x_1^*,\ldots,x_{s-1}^*$ is a $\graded(I)$-sequence. Then there is a
form $Y^{r+1}-\sum X_iF_i$, with $F_i\in k[X_1,\ldots ,X_s,Y]$ forms
of degree $r$ and $\fiber(I)\cong k[X_1,\ldots ,X_s,Y]/(Y^{r+1}-\sum
X_iF_i)$. In particular, 
$\reltype(\fiber(I))=\rednumber_J(I)+1=\reltype(I)$.
\end{Corollary}
\demo
By Theorem~A, $\reltype(\fiber(I))\leq
\reltype(I)=r+1$. By \cite[Lemma~5.2]{hk}, $E(\fiber(I))_{n}=0$ for
all $2\leq n\leq r$ and $E(\fiber(I))_{r+1}\neq 0$. Thus $\fiber(I)$
has only equations of degree $r+1$ and $\reltype(\fiber(I))=r+1$. By
Theorem~A, $E(I)_{r+1}$ is cyclic and generated by
the equation of $\rees(I)$ induced by the containment $y^{r+1}\in
JI^{r}$. Therefore the same happens with $E(\fiber(I))_{r+1}$ (see
\cite[Proposition~3.2]{gp} or \cite[p.~268]{hku}).
\qed

\medskip
Observe that the constraints on $I$ in Corollary \ref{hk-generalisation} 
automatically hold for $I$ an ideal with $\mu(I)=\spread(I)+1$,
$\grade(\graded(I)_+)\geq \spread(I)-1$ and such that there is a minimal 
reduction $J\subseteq I$ generated by an $R$-sequence (see \cite[Proposition 1.5.12]{bh}). 

\begin{Remark}{\rm
An immediate consequence of Corollary \ref{hk-generalisation} 
is that $\fiber(I)$ is Cohen-Macaulay. In fact, recall that 
\cite[Proposition 5.4]{hk} already prove that, under the 
conditions of Corollary \ref{hk-generalisation}, $\fiber(I)$ 
is a hypersurface if and only if $\fiber(I)$ is Cohen-Macaulay.
}\end{Remark}

\section{Examples and applications}\label{some-examples}

Our purpose in this section is to use Theorem~A in order to
deduce minimal generating sets of the equations of $\rees(I)$. 
If $J$ is a reduction of $I=(J,y)$, the
ascending chain of colon ideals $\{(JI^{n-1}:y^{n})\}_{n\geq 1}$,
which reaches eventually $R$ as $n$ increases, can be calculated 
in any computer algebra system, giving an alternative strategy 
to find out the equations of $\rees(I)$.

The ideal $I$ in the next example, for a specific $p\geq 1$, 
is often used as a paradigm of an ideal of relation type $=p$. 
As said above, Theorem~A will be crucial to our purposes.

\begin{Example}\label{exemple-classic} {\rm
Let $(R,\fm)$ be a Noetherian local ring. Let $a,b$ an
$R$-sequence and $p\geq 2$. Set $x_{1}=a^{p}$, $x_{2}=b^{p}$
and $y=ab^{p-1}$. Let $I$ be the ideal generated by
$x_{1},x_{2},y$. Set $V=R[X_{1},X_{2},Y]$ and let $\varphi:V\to
\rees(I)$ be the presentation of $\rees(I)$ sending $X_{i}$ to $x_it$
and $Y$ to $yt$. Then a minimal generating set of the ideal
$Q=\ker(\varphi)$ is given by:	
\begin{itemize}
\item[$(i)$] two linear forms $F_1(X_1,X_2;Y)=a^{p-1}Y-b^{p-1}X_1$ and
  $G_1(X_1,X_2,Y)=bY-aX_2$;
\item[$(ii)$] for each $2\leq n\leq p$, a unique equation
  $F_n(X_1,X_2,Y)=a^{p-n}Y^n-b^{p-n}X_1X_2^{n-1}$.
\end{itemize}
}\end{Example}
\demo
\noindent We start by proving that $I$ fulfils the hypotheses
of Theorem~A. Clearly $x_{1},x_{2}$ is an $R$-sequence
and $J=(x_1,x_2)$ is a reduction of $I$ since $I^{p}=JI^{p-1}$. By
\cite[Corollary~3]{ks}, a monomial $m$ on $a,b$ belongs to an
ideal generated by monomials $m_{1},\ldots ,m_{r}$ on $a,b$ if
and only if $m$ is a multiple of some $m_{i}$. It follows that
$y^{p-1}\not\in JI^{p-2}$ and $I^{p-1}\not\subset JI^{p-2}$. Thus
$\rednumber_J(I)=p-1$.

\noindent {\bf Claim 1}: $\grade(\graded(I)_+)\geq 1$ and
$x_{1}^{*}$ is $\graded(I)$-regular.  

\noindent {\bf Proof of Claim 1:} By \cite[Corollary~2.7]{vv}, it suffices to
prove 
\begin{eqnarray*}
VV_{(x_1)}(I)_n=x_{1}R\cap I^{n}/x_{1}I^{n-1}=0
\end{eqnarray*} 
for all $n\geq 1$. Fix $n\geq 1$. By \cite[Proposition~1]{ks},
\begin{eqnarray*}
x_{1}R\cap I^{n}=(L_{i,j,k}\; |\; 
i,j,k \text{ positive integers such that } i+j+k=n),
\end{eqnarray*} 
where
$L_{i,j,k}=\lcm(a^{p},a^{ip}b^{jp}(ab^{p-1})^{k})$.
Let us prove that $L_{i,j,k}$ is in $x_1I^{n-1}$.

Indeed, if $i\geq 1$, then
$L_{i,j,k}=a^{ip}b^{jp}(ab^{p-1})^{k}=
a^{p}[a^{(i-1)p}b^{jp}(ab^{p-1})^{k}]\in
x_{1}I^{n-1}$ and we have finished. Hence we can suppose $i=0$ and
$j+k=n$. If $k=0$, then $j=n$ and $L_{0,j,0}=a^{p}b^{jp}\in
x_{1}I^{n-1}$. Suppose $0<k\leq p$. Then
$L_{0,j,k}=a^{p}b^{jp+k(p-1)}=a^{p}(b^{p})^{j+k-1}b^{p-k}\in
x_1I^{n-1}$. Finally, if $k>p$, then
$L_{0,j,k}=b^{jp}(ab^{p-1})^{k}=
a^{p}[a^{k-p}b^{jp+k(p-1)}]=
a^{p}[a^{k-p}b^{(k-p)(p-1)}b^{jp+p(p-1)}]=
a^{p}(b^{p})^{j+p-1}(ab^{p-1})^{k-p}\in x_{1}I^{n-1}$.

\medskip

Note that for $p>2$, then $x_{2}^{*}$ is not a $\graded(I)$-sequence
because
\begin{eqnarray*}
(a^{p+2}b^{p-2})x_{2}=a^{p}(ab^{p-1})^{2}=x_{1}y^{2}\in
I^{3},
\end{eqnarray*}
where $(a^{p+2}b^{p-2})\in I\setminus I^{2}$.

Therefore, we can apply Theorem~A and deduce that,
for all $n\geq 2$,
\begin{eqnarray*}
\left[ \frac{Q}{Q\langle n-1\rangle}\right]_{n}\cong
\frac{JI^{n-1}:y^{n}}{JI^{n-2}:y^{n-1}}.
\end{eqnarray*}
In particular, since $(JI^{p-2}:y^{p-1})\subsetneq
(JI^{p-1}:y^{p})=R$, then $\reltype(I)=\rednumber_J(I)+1$.

\medskip

\noindent {\bf Claim 2}: $(JI^{n-1}:y^{n})=(a^{p-n},b)$
for $2\leq n\leq p-1$. 

\noindent {\bf Proof of Claim 2:}
First note that, for all $2\leq n\leq p-1$,
$b\in (JI^{n-1}:y^{n})$ since
$by^n=a^{n}b^{n(p-1)+1}=ab^p(a^{n-1}b^{(n-1)(p-1)})=
ab^p(ab^{p-1})^{n-1}\in JI^{n-1}$. Since $(JI^{n-1}:y^{n})$ is
generated by monomials in $a,b$ (c.f. Remark~\ref{monomial-ideal-operations}) 
and $b\in (JI^{n-1}:y^{n})$, there is only one possible remaining generator: 
a power of $a$. Since $a^{p-n}y^n=a^{p-n} a^n b^{n(p-1)}=
a^pb^{n(p-1)}=b^{p-n}a^p(b^p)^{n-1}\in JI^{n-1}$, then
$a^{p-n}\in (JI^{n-1}:y^{n})$. However, and using again
\cite[Corollary~3]{ks}, one has $a^{p-n-1}\notin (JI^{n-1}:y^{n})$.

\medskip

Hence $(Q/Q\langle p-1\rangle)_{p}\cong
(JI^{p-1}:y^{p})/(JI^{p-2}:y^{p-1})=R/(a,b)$ and, for $2\leq
n\leq p-1$,
\begin{eqnarray*}
\left[ \frac{Q}{Q\langle n-1\rangle}\right]_{n}\cong
\frac{JI^{n-1}:y^{n}}{JI^{n-2}:y^{n-1}}=
\frac{(a^{p-n},b)}{(a^{p-n+1},b)}\cong
\frac{(a^{p-n})}{(a^{p-n+1},a^{p-n}b)}.
\end{eqnarray*}
In other words, for each $2\leq n\leq p$, $(Q/Q\langle
n-1\rangle)_{n}$ is generated by a single element that corresponds to
the class of $a^{p-n}$ ($1$ if $n=p$). To find this element,
consider the identity $a^{p-n}y^n=b^{p-n}a^p(b^p)^{n-1}$,
which induces the equation
$F_{n}(X_1,X_2,Y)=a^{p-n}Y^n-b^{p-n}X_1X_2^{n-1}\in Q_n$. Since
the isomorphism of Theorem~A sends the class of $F_n$
to the class of $F_n(0,0,1)=a^{p-n}$, we are done. By
Remark~\ref{min-gen-set-of-equations}, $F_n$ is in a minimal
generating set of $Q$, for $2\leq n\leq p$.

To finish, let us find the equations of degree one. Although this is
trivial, we sketch the proof here to show the similarity with the
greater degrees. As before, one shows that $J=(a^{p-1},b)$.
Using Remark~\ref{syzygies},
\begin{eqnarray*}
\koszul{1}{x_1,x_2,y}{R}\stackrel{\tilde{\sigma}_1}{\cong}
\frac{(J:y)}{J}=\frac{(a^{p-1},b)}{(a^p,b^p)}.
\end{eqnarray*}
Identify $Q_1$ with $\cycle{1}{x_1,x_2,y}{R}$ and $B_1=\langle
x_1Y-yX_1, x_2Y-yX_2,x_1X_2-x_2X_1\rangle$ with
$\boundary{1}{x_1,x_2,y}{R}$.  Then $Q_1/B_1$ is minimally generated
by the classes of the two equations corresponding to the classes of
$a^{p-1}$ and $b$. The identities $a^{p-1}y=b^{p-1}a^p$ and
$by=ab^p$ induce the desired equations
$F_1(X_1,X_2,Y)=a^{p-1}Y-b^{p-1}X_1$ and
$G_1(X_1,X_2,Y)=bY-aX_2\in Q_1$, since by $\tilde{\sigma}_1$ their
classes are sent to the classes of $F_1(0,0,1)=a^{p-1}$ and
$G_1(0,0,1)=b$. Clearly, $B_1\subset \langle F_1,G_1\rangle$ and
$F_1,G_1$ are a minimal generating set of $Q_1$.
\qed

\medskip

The following example gives the equations of the Rees algebra of an
$\fm$-primary ideal of almost-linear type having a reduction
generated by a $d$-sequence. We will use Proposition~\ref{rednumber=1-theorem}.

\begin{Example}\label{typical-buchsbaum}{\rm
Let $S=k[[X_1,X_2,U_1,U_2]]$ be the power series ring over $k$ 
with indeterminates $X_1,X_2,U_1,U_2$, and let 
$R=k[[X_1,X_2,U_1,U_2]]/L=k[[x_1,x_2,u_1,u_2]]$, 
where $L=(X_1,X_2)\cap (U_1,U_2)=(X_1U_1,X_1U_2,X_2U_1,X_2U_2)$.
Set $J=(x_1+u_1,x_2+u_2)$ and $I=(J,x_1+x_2)$. Let
$\varphi:R[T_1,T_2,Y]\to\rees(I)$ be the polynomial presentation 
that sends $T_i$ into $x_i+u_i$ and $Y$ into $x_1+x_2$ and 
let $Q=\ker\varphi$. Then $J$ is a reduction of $I$, 
with $\rednumber_J(I)=1$, $\reltype(I)=2$ and 
$\ker\alpha_{I,2}\cong (Q/Q\langle 1\rangle)_2\cong R/(J:(x_1+x_2))=k$
is cyclic, generated by $Y^2-Y(T_1+T_2)$.  
}\end{Example}
\demo
The ring $R$ is Noetherian local, $2$-dimensional and Buchsbaum (in fact, the ideal 
$(X_1,X_2)\cap (U_1,U_2)$ is said to be of type $(1,2)$, 
see \cite[Theorem 3 and Definition 2, p. 742]{sv}), consequently 
every system of parameters is a $d$-sequence. Since $x_1+u_1,x_2+u_2$ 
is a system of parameters of $R$, it is a $d$-sequence. It is also 
readily verified that $(x_1+u_1:x_2+u_2)\cap I=(x_1+u_1)$ and 
$(x_1+x_2)\cap I^2=(x_1+x_2)I$. By Proposition \ref{rednumber=1-theorem}, 
the claim follows.
\qed

\medskip

The next example shows that the hypotheses $(a)$ and $(b)$ in
Remark~\ref{different-version-theorem} alone are not sufficient 
to ensure that there is only one equation of $\rees(I)$ of maximum 
degree.

\begin{Example}\label{not-enough} {\rm 
Let $k$ be a field and $R=k[X,Y]_{(X,Y)}/(XY,Y^{2})_{(X,Y)}$. Set $x$
and $y$ the classes of $X$ and $Y$ in $R$. Let $\fm=(x,y)$ be the
maximal ideal of $R$. Then, $\reltype(\fm)=2$ and there
are two quadratic equations in a minimal generating set 
of the equations of $\rees(\fm)$.}\end{Example}
\demo
Set $J=(x)$. Since $y^2=0\in J\fm$ and $y\not\in J$, then $\fm^2=J\fm$
and $J$ is a reduction of $\fm$ with reduction number $1$. Moreover,
since $(0:x)=(0:x^2)$, $x$ is a $d$-sequence. By
Remark~\ref{different-version-theorem}, $\reltype(\fm)\leq
\rednumber_J(\fm)+1=2$. Set $V=k[S,T]$ and let $\psi:V\to
\graded(\fm)$ be the presentation of $\graded(\fm)$ sending $S$ to
$x+\fm^2$ and $T$ to $y+\fm^2$. For $n\geq 2$, $\fm^n=(x^n)$.  Thus
$\fm^2/\fm^3$ is a $k$-vector space of dimension 1. Therefore
$\ker(\psi_2)\subset V_2$ must have dimension 2. In fact,
$\ker(\psi_2)=\langle ST,T^2\rangle$. Since $\ker(\psi_1)=0$, then
$E(\graded(\fm))_2$ is minimally generated by two elements. 
We finish by using that $E(\fm)_2\cong E(\graded(\fm))_2$
(see \ref{defining-equations-associated-graded} and 
\cite[Proposition~3.3]{planas}). Observe that in this case, 
$J$ and $I$ do not fulfil the condition $(b')$ of 
Proposition \ref{rednumber=1-theorem}:
in fact, $(0:x)\cap \fm = (y) \neq 0$.\qed

\medskip

As it has already been pointed out in Proposition \ref{rednumber=1-theorem}, 
we can relax the hypotheses in Theorem~A if the reduction number is small. 
The next example illustrates this by computing the equations for a 
non-equimultiple ideal $I$ with ${\rm ad}(I)=1$.

\begin{Example}\label{exemple-C4}{\rm
Let $(R,\fm)$ a Noetherian regular local ring with infinite residue field $k$ 
and dimension $4$, with $u,v,z,t$ a regular system of parameters. 
Set $I=(u,z)\cap (v,t)$. Then $I=(J,ut)$ where $J=(uv,vz-ut,zt)$ is a 
minimal reduction of $I$ with reduction number $\rednumber_J(I)=1$ and 
such that $uv,vz-ut,zt$ is a $d$-sequence. 
Let $\varphi:R[X_1,X_2,X_3,Y]\to\rees(I)$ be the polynomial presentation 
relative to the minimal generating set $uv,vz-ut,zt,ut$ of $I$ and let 
$Q=\ker\varphi$ be the ideal of equations of $\rees(I)$. Then 
$\reltype(I)=2$ and $Q$ has only one minimal equation of degree 
$2$: $Y^2+X_2Y-X_1X_3$. Notice that $\grade(I)=2$, $\height(I)=2$, 
$\spread(I)=3$ and $\mu(I)=4$.
}\end{Example}
\demo
Notice that $I$ is generated by $uv,vz,zt,ut$. We know that 
$\fiber(I)\cong k[uv,vz,zt,ut]$ (see \cite[Example~1.90]{vasconcelos}). 
Displaying a Noether normalisation, one can readily check that $uv,vz-ut,zt$ 
generate a minimal reduction $J$ of $I$ (see \cite[Exercise~16, p.69]{am}). 
It is readily seen that $uv,vz-ut,zt$ is a $d$-sequence 
but not an $R$-sequence, that $uv,vz-ut,zt$ is a minimal generating set 
of $J$ and that $JI=I^2$, since $(ut)^2=(uv)(zt)-(vz-ut)(ut)\in JI$, 
hence $\rednumber_J(I)=1$. In order to call on 
Proposition~\ref{rednumber=1-theorem} one may check the 
following identities:
\begin{enumerate}
\item[$(b')$] $(0:uv)\cap I=0$ (for $i=1$), 
\item[] $(uv:vz-ut)\cap I=(uv)$ (for $i=2$), 
\item[] $((uv,vz-ut):zt)\cap I=(uv,vz-ut)$ (for $i=3$);
\item[$(c')$] $(uv)\cap I^2=(uv)I$ (for $i=2$), 
\item[] $(uv, vz-ut)\cap I^2=(uv,vz-ut)I$ (for $i=3$).
\end{enumerate}
Consequently, $Q=Q\langle 2\rangle$, $\reltype(I)=2$ and 
$(Q/Q\langle 1\rangle)_2\cong (JI:(ut)^2)/(J:ut)\cong R/\fm$ 
is cyclic with generator $Y^2+X_2Y-X_1X_3$.
\qed

\begin{Example}\label{exemple-pseudo-classic}{\rm
Let $(R,\fm)$ be a Noetherian local ring. Let $a,b$ an
$R$-sequence, let $p\geq 5$ be an odd integer, and 
set $x_{1}=a^{p}$, $x_{2}=b^{p}$ and $y=a^2 b^{p-2}$. Let 
$J=(x_1,x_2)$ and $I=(x_1,x_2,y)$, $V=R[X_{1},X_{2},Y]$ 
the polynomial ring in three indeterminates over $R$ 
and let $\varphi:V\to\rees(I)$ be the presentation of 
$\rees(I)$ sending $X_{i}$ to $x_it$
and $Y$ to $yt$.
\begin{itemize}
\item[$(i)$] The following equalities hold: 
\begin{eqnarray*}
(JI^n:y^{n+1}) = 
\left\{
\begin{array}{ll}
(a^{p-2n-2},b^2) & 1\leq n\leq \frac{p-3}{2}, \\
(a,b) & \frac{p-1}{2}\leq n\leq p-2, \\
R & n\geq p-1.
\end{array} \right.
\end{eqnarray*}
In particular, $J$ is a reduction of $I$ with $\rednumber_J(I)=p-1$.
\item[$(ii)$] The following equalities hold: $(x_1I^{n-1}:x_2)\cap I^{n-1}/x_1I^{n-2}=0$, for all $n\geq 2$; 
equivalently, $\koszul{1}{x_1t,x_2t}{\rees(I)}_n=0$, for all $n\geq 2$.
\item[$(iii)$] It is verified that $(Q/Q\langle n-1\rangle)_{n}\cong (JI^{n-1}:y^n)/(JI^{n-2}:y^{n-1})$, 
for all $n\geq 2$. In particular $\reltype(I)=p$. Moreover,  $$\mathcal{Q}=\{S_{1,2},S_{1,p-2},F_2,\ldots,F_{(p+1)/2},F_p\}$$ 
is a minimal generating set of $Q$ with
\begin{eqnarray*}
F_n(X_1,X_2,Y)=\left\{
\begin{array}{ll}
a^{p-2n} Y^n - b^{p-2n} X_1 X_2^{n-1} & 2\leq n\leq \frac{p-1}{2}, \\
bY^n - a X_1 X_2^{n-1} & n=\frac{p+1}{2}, \\
Y^p - X_1^2 X_2^{p-2} & n=p,
\end{array} \right.
\end{eqnarray*}
and $S_{1,2}(X_1,X_2,Y)=a^2 X_2 - b^2 Y $ and $S_{1,p-2}(X_1,X_2,Y)=a^{p-2} Y - b^{p-2} X_1$.
\item[$(iv)$] The initial forms $x_1^*,x_2^*$ and $y^*$ are zero-divisors in $\graded(I)$.
\end{itemize}
}\end{Example}
\demo
In order to show $(i)$, notice first that since 
$y^{p}=x_1^2x_2^{p-2}\in J^p\subset JI^{p-1}$, we have
$(JI^{p-1}:y^p)=R$ and $\rednumber_J(I)\leq p-1$. Let us show 
that $y^{n+1}\notin JI^{n}$ for $n\leq p-2$: assuming the opposite, 
then we can write $y^{n+1}=x_1^rx_2^sy^t$ with $r,s,t$ positive 
integer such that $r+s+t=n+1$. If $t\geq 1$ we can simplify the 
expressions, hence we can assume without loss of generality that 
$y^{n+1}=x_1^rx_2^s$ with $r+s=n+1$ and $n\leq p-2$. Then
$a^{2(n+1)}b^{(p-2)(n+1)}=a^{rp}b^{sp}$, thus $2(n+1)=rp$. 
But observe that $n+1\leq p-1<p$, so $r$ has to be $1$, leading to 
$p=2(n+1)$, a contradiction with the assumption that $p$ is odd. 
Consequently, $\rednumber_J(I)=p-1$.

We claim that if $(p-1)/2 \leq n\leq p-2$, then 
$(JI^n:y^{n+1})=(a,b)$. Playing with the exponents we get 
$$ay^{n+1}=a^{p}a^{2n-p+3}b^{(2n-p+3)(p-2)/2}b^{(p-1)(p-2)/2}=bx_1x_2^{(p-3)/2}y^{n-(p-3)/2}\in JI^n,$$ 
with $2n-p+3$ being positive by our assumptions on $n$.
On the other hand,
$$by^{n+1}=aa^pa^{2n-p+1}b^{(2n-p+1)(p-2)/2}b^{p(p-1)/2}=ax_1x_2^{(p-1)/2}y^{n-(p-1)/2}\in JI^n,$$ 
with $2n-p+1$ being positive by our assumptions on $n$. Therefore $(a,b)\subseteq (JI^n:y^{n+1})$ and 
using Remark~\ref{monomial-ideal-operations} the conclusion follows.

Observe that if $n\leq (p-3)/2$ then $b^2\in (JI^n:y^{n+1})$, since 
$b^2y^{n+1}=a^2 b^p a^{2n} b^{n(p-2)}=a^2x_2y^n\in JI^n$. However,
in this case $b\notin (JI^n:y^{n+1})$. If $by^{n+1}\in JI^n$ we can 
assume without loss of generality that $by^{n+1}$ is of the form
$ax_1^rx_2^s$ with $r,s$ positive integers such that $r+s=n+1$, 
but then $2(n+1)=rp+1$ and since $2(n+1)\leq p-2$ we arrive to a 
contradiction.

Likewise, if $n\leq (p-3)/2$ then $a^{p-2n-2}\in (JI^n:y^{n+1})$, 
since $a^{p-2n-2}y^{n+1}=a^pb^{pn}b^{p-2n-2}=b^{p-2n-2}x_1x_2^n\in JI^n$. 
However, if $k$ is a positive integer such that $k\leq p-2n-3$, 
we have $a^k\notin (JI^n:y^{n+1})$. Assuming the opposite, 
i.e., $a^ky^{n+1}\in (JI^n:y^{n+1})$, without loss of generality 
$a^k y^{n+1}=b^k x_1^r x_2^s$ with $r,s$ positive integers 
such that $r+s=n+1$. But then it follows that $2(n+1)+k=rp$,
a contradiction, since $2(n+1)+k\leq p-1$. 

To finish with this description of the generators of $(JI^n:y^{n+1})$, 
let us show that for $n\leq (p-3)/2$ and $k\leq p-2n-3$, 
$a^k b\notin (JI^n:y^{n+1})$. 
If $a^k b y^{n+1} \in JI^n$ then we can assume without loss of generality
that $a^k b a^{2(n+1)}b^{(p-2)(n+1)}=mx_1^rx_2^s$ with $r,s$ positive 
integers such that $r+s=n+1$ and $m$ a monomial in $a,b$ with 
${\rm deg}(m)=k+1$. If $m$ is multiple of $b$ then, simplifying $b$'s, 
we get $a^k\in (JI^n:y^{n+1})$, which is in contradiction with the 
preceding paragraph. Thus $m=a^{k+1}$, but then, simplifying $a$'s, we 
would deduce that $b\in (JI^n:y^{n+1})$, which is a contradiction with 
the claims above. Therefore, when $n\leq (p-3)/2$ the ideal 
$(JI^n:y^{n+1})$ is the monomial ideal generated by 
$a^{p-2n-2}$ and $b^2$.

\medskip

Let us show $(ii)$. We will show the equality $(x_1I^n:x_2)\cap I^n=x_1I^{n-1}$ 
for all $n\geq 1$, and then apply Lemma \ref{characterisation-of-homology}. 
For that purpose we will make use of Remark \ref{monomial-ideal-operations}. Concretely, 
we will show that any monomial generator of $(x_1I^n:x_2)\cap I^n$ belongs to 
$x_1I^{n-1}$. Each monomial generator of $(x_1I^n:x_2)\cap I^n$ is of the form 
\begin{eqnarray*}
G_{\sigma,\tau}={\rm lcm}\left(M_\sigma,x_1^ux_2^vy^w \right),
\end{eqnarray*}
where $M_\sigma={\rm lcm}(x_1^{r+1}x_2^sy^t,x_2)/x_2$, with $r+s+t=n$, 
$u+v+w=n$; $\sigma$ and $\tau$ denote the integer vectors $(r,s,t)\in\mathbb{N}^3$ 
and $(u,v,w)\in\mathbb{N}^3$, respectively. 
Let us show that $G_{\sigma,\tau}\in x_1I^{n-1}$.

\medskip

\noindent {\bf Case:} $s\geq 1$. 

\noindent So $M_\sigma = x_1^{r+1}x_2^{s-1}y^t$, 
and then $G_{\sigma,\tau}$ is a multiple of $M_\sigma$ which can be 
written as $x_1\cdot (x_1^rx_2^{s-1}y^t)\in x_1I^{n-1}$.

\medskip

\newpage

\noindent {\bf Case:} $s=0$. 

\noindent So $M_\sigma={\rm lcm}(x_1^{r+1}y^t,x_2)/x_2={\rm lcm}(a^{p(r+1)+2t}b^{(p-2)t},b^p)/b^p$. At this point
it is convenient to distinguish a couple of subcases: $t=0,1$ and $t\geq 2$.
\newline
{\it Sub-case:} $t=0,1$. 

If $t=0$, $M_\sigma=a^{p(n+1)}$, and then $G_{\sigma,\tau}$ is a multiple of $a^{p(n+1)}=x_1^{n+1}\in x_1 I^{n}\subset x_1I^{n-1}$. 
If $t=1$, $M_\sigma=a^{pn+2}$, and then $G_{\sigma,\tau}={\rm lcm}(a^{pn+2},a^{pu+2w}b^{pv+(p-2)w})$. 
Since $pn+2>pu+2w$ for all $(u,v,w)\in\mathbb{N}^3$ such that $u+v+w=n$, 
$G_{\sigma,\tau}=a^{pn+2}b^{pv+(p-2)w}=a^p a^{p(n-1)+2}b^{pv+(p-2)w}\in x_1 I^{n-1}$.

{\it Sub-case:} $t\geq 2$. 

In this case, since $p\geq 5$, the inequality $(p-2)t>p$ holds, 
and then $$M_\sigma={\rm lcm}(a^{p(r+1)+2t}b^{(p-2)t},b^p)/b^p=a^{p(r+1)+2t}b^{(p-2)t-p}.$$ 
Therefore, we have to check 
$$G_{\sigma,\tau}={\rm lcm}\left( a^{p(r+1)+2t}b^{(p-2)t-p},x_1^ux_2^vy^w \right)\in x_1 I^{n-1}.$$
Since $G_{\sigma,\tau}$ is a multiple of $x_1^ux_2^vy^w$, if $u\geq 1$, then 
$G_{\sigma,\tau}$ will be multiple of $x_1\cdot x_1^{u-1}x_2^vy^w\in x_1I^{n-1}$. So it 
remains to consider the case where $u=0$. In this case,
\begin{eqnarray*}
G_{\sigma,\tau} & = & {\rm lcm}\left( a^{p(r+1)+2t}b^{(p-2)t-p}, a^{2w}b^{pv+(p-2)w}\right)\\
								& = & {\rm lcm}\left( a^{p(n-t+1)+2t}b^{(p-2)t-p}, a^{2w}b^{p(n-w)+(p-2)w}\right).
\end{eqnarray*}
We claim that $2w\leq p(n-t+1)+2t$. The least value that $p(n-t+1)+2t$ can take is $p+2n$, 
whereas the greatest value that $2w$ can take is $2n$, then $2w\leq p(n-t+1)+2t$. So 
$p(n-w)+(p-2)w=pn-2w \geq pn-p(n-t+1)-2t=t(p-2)-p>0$. Hence 
$G_{\sigma,\tau}=a^{p(n-t+1)+2t}b^{p(n-w)+(p-2)w}=a^{p}\cdot a^{p(n-t)+2t}b^{p(n-w)+(p-2)w}$. 
If $w<t$ then $G_{\sigma,\tau}$ is a multiple of $a^{p}\cdot a^{p(n-t)+2t}b^{p(n-t)+(p-2)t}=
x_1x_1^{n-t}x_2^{n-t}y^t\in x_1I^{2n-t}\subset x_1I^{n}\subset x_1I^{n-1}$. If $w\geq t$ then
$G_{\sigma,\tau}$ is a multiple of $a^{p(n-w+1)+2w}b^{p(n-w)+(p-2)w}=
a^{p}\cdot a^{p(n-w)+2w}b^{p(n-w)+(p-2)w}=x_1x_1^{n-w}x_2^{n-w}y^w\in x_1I^{2n-w}\subset x_1I^{n}
\subset x_1I^{n-1}$.

\medskip

To prove $(iii)$ simply use $(ii)$ and Theorem~A.

\medskip
 
In order to prove $(iv)$ it suffices to check that both $x_1^*$ 
and $x_2^*$ are zero-divisors in $\graded(I)$. It will be sufficient 
to find elements $\omega_1\in I^{n_1}\backslash I^{n_1+1}$ and 
$\omega_2\in I^{n_2}\backslash I^{n_2+1}$ such that 
$x_1\omega_1\in I^{n_1+2}$ and $x_2\omega_2\in I^{n_2+2}$. 
Notice that the exponents $n_1,n_2$ may depend on $p$, 
as the ideal $I$ does. Let $n_1(p)=(p-3)/2$, $n_2=0$ and take 
$\omega_1=ab^{(n_1(p)+2)\cdot(p-2)}\in I^{n_1(p)}\backslash I^{n_1(p)+1}$ and 
$\omega_2=a^4b^{p-4}\in R\notin I$. Then $x_1\omega_1=a^{p+1}b^{(n_1(p)+2)(p-2)}$ 
and since $p+1=2(n_1(p)+2)$ we get 
$x_1\omega_1=(a^2b^{p-2})^{n_1(p)+2}=y^{n_1(p)+2}\in I^{n_1(p)+2}$ 
and, on the other hand, 
$x_2\omega_2=a^4b^{2p-4}=(a^2b^{p-2})^2=y^2\in I^2$.
\qed

\medskip

The next example shows the non-surprising fact that, although
the $(Q/Q\langle n-1\rangle)_n$ are isomorphic to quotients of the form
$(JI^{n-1}:y^{n})/(JI^{n-2}:y^{n-1})$, they need not be cyclic in 
degrees $< \rednumber_J(I)+1$. The proof in full detail of Example~\ref{exemple-dim=4}
would be similar to that of Example~\ref{exemple-pseudo-classic}.

\begin{Example}\label{exemple-dim=4}{\rm
Let $(R,\fm)$ be a Noetherian local ring. Let $a,b,c$ an
$R$-sequence, let $p\geq 4$ be an even integer and 
set $x_{1}=a^{p}$, $x_{2}=b^{p}$, $x_{3}=c^p$ 
and $y=abc^{p-2}$. Let $J=(x_1,x_2,x_3)$ and 
$I=(x_1,x_2,x_3,y)$, $V=R[X_{1},X_{2},X_{3},Y]$ the 
polynomial ring in four indeterminates over $R$ and 
let $\varphi:V\to\rees(I)$ be the presentation of 
$\rees(I)$ sending $X_{i}$ to $x_it$
and $Y$ to $yt$. Then,
\begin{itemize}
\item[$(i)$] The ideal $J$ is a reduction of $I$ with $\rednumber_J(I)=p-1$:
concretely, 
\begin{eqnarray*}
(JI^n:y^{n+1}) = 
\left\{
\begin{array}{ll}
(a^{p-n-1},b^{p-n-1},c^2) & 0\leq n < \frac{p}{2}, \\
(a^{p/2},b^{p/2},c^2,a^{p-n-1}b^{p-n-1}) & \frac{p}{2}\leq n < p-1, \\
R & n\geq p-1.
\end{array} \right.
\end{eqnarray*}

\item[$(ii)$] The initial forms $x_1^*,x_2^*$ are a $\graded(I)$-sequence, in particular
$\grade(\graded(I)_+)\geq 2$; since $x_1,x_2,x_3$ are an $R$-sequence, by Lemma~\ref{homologia-zero}, 
we get $\koszul{1}{x_1t,x_2t,x_3t}{\rees(I)}_n=0$, for all $n\geq 2$;

\item[$(iii)$] It is verified that $(Q/Q\langle n-1\rangle)_{n}\cong (JI^{n-1}:y^n)/(JI^{n-2}:y^{n-1})$, 
for all $n\geq 2$. In particular, $\reltype(I)=p$. Moreover, in terms of generators, the graded components 
of the modules $(Q/Q\langle n-1\rangle)_n$ are given by
\begin{eqnarray*}
(Q/Q\langle n-1\rangle)_n=\left\{
\begin{array}{ll}
\langle a^{p-n}Y^{n}-b^{n}c^{p-2n}X_1X_3^{n-1}, & \\
b^{p-n}Y^{n}-a^{n}c^{p-2n}X_2X_3^{n-1} \rangle & 2 \leq n \leq \frac{p}{2}, \\
\langle a^{n-2}b^{n-2}Y^{n}-c^{2n-4}X_1X_2X_3^{n-2} \rangle & n = \frac{p}{2}+1, \\
\langle a^{p-n}b^{p-n}Y^{n}-c^{2(p-n)}X_1X_2X_3^{n-2} \rangle & \frac{p}{2}+1 < n < p, \\
\langle Y^{n}-X_1X_2X_3^{n-2} \rangle & n=p.
\end{array} \right.
\end{eqnarray*}

Recall that $Z_1(x_1,x_2,x_3,y;R)$ and $B_1(x_1,x_2,x_3,y;R)$ stand for the Koszul $1$-cycles and $1$-boundaries, 
respectively, as introduced in Section~\ref{equations-of-u}. Identifying $\cycle{1}{x_1,x_2,x_3,y}{R}$ with $Q_1$ and
$\boundary{1}{x_1,x_2,x_3,y}{R}$ with 
$$P_1=\langle x_1Y-yX_1, x_2Y-yX_2, x_3Y-yX_3, x_2X_3-x_3X_2, x_1X_3-x_3X_1, x_1X_2-x_2X_1\rangle,$$
then we have
$$\koszul{1}{x_1,x_2,x_3,y}{R}\cong\frac{Q_1}{P_1}\cong  \frac{(J:y)}{J}\cong\frac{(a^{p-1},b^{p-1},c^{2})}{(a^{p},b^{p},c^{p})},$$ 
which is minimally generated by the classes of the equations corresponding to the 
classes of $a^{p-1}$, $b^{p-1}$ and $c^2$:
$$\frac{Q_1}{P_1}\cong \langle a^{p-1}Y-bc^{p-2}X_1, b^{p-1}Y-ac^{p-2}X_2,c^2Y-abX_3\rangle.$$
\end{itemize}
}\end{Example}
\demo
Assertion $(i)$ follows by applying the same strategy as in the proof of Example 
\ref{exemple-pseudo-classic}. In order to prove $(ii)$ recall Remark~\ref{valabrega-valla}: 
since $x_1,x_2$ is an $R$-sequence, it suffices to check that $(x_1)\cap I^n=(x_1)I^{n-1}$ 
and $(x_1,x_2)\cap I^n=(x_1,x_2)I^{n-1}$, for all $n\geq 2$. In order to prove $(iii)$ just 
combine $(i)$, $(ii)$, Theorem~A and Remark \ref{syzygies}.
\qed

\medskip

The next example shows that if ${\rm sd}(I)=\mu(I)-\spread(I)\neq 1$, 
there might be several equations of $\rees(I)$ of top degree.

\begin{Example}\label{exemple-m-h=2} {\rm 
Let $(R,\fm)$ be a two dimensional regular local ring and 
let $p\geq 2$ be an integer. Then $\reltype(\fm^p)=2$ and there are 
$\binom{p}{2}$ quadratic equations in a minimal generating set of 
equations of $\rees(\fm^p)$.
}\end{Example}
\demo
Let $x,y$ be a regular system of parameters of $R$, $V=R[X,Y]$ and let
$\varphi:V\to\rees(\fm)$ be the presentation of $\rees(\fm)$ sending
$X$ to $xt$ and $Y$ to $yt$. Since $x,y$ is an $R$-sequence,
$\ker(\varphi)=\langle xY-yX\rangle$. Set $V(p)=R[X^p,X^{p-1}Y,\ldots
  ,XY^{p-1},Y^p]$ and $\varphi(p):V(p)\to \rees(\fm^p)$ the $p$-th
Veronese transform of $\varphi$. Note that
$\ker(\varphi(p)_n)=\ker(\varphi_{pn})=\langle xY-yX\rangle V_{pn-1}$.

Set $W=R[T_0,T_1,\ldots ,T_p]$ and let $\psi:W\to V(p)$ be the
polynomial presentation of $V(p)$ sending $T_i$ to $X^{p-i}Y^{i}$. It
is known (see e.g. \cite[Proposition~2.5]{jk}) that the kernel of
$\psi$ is the determinantal ideal generated by the $2\times 2$ minors
of the matrix $\cm$, where $\cm$ is
\begin{displaymath}
\left( \begin{array}{cccc} T_0 & T_1 & \ldots & T_{p-1} \\ T_1 & T_2 &
  \ldots & T_p \end{array} \right).
\end{displaymath}
In particular, $\ker(\psi)=\ker(\psi)\langle 2\rangle$.

Consider $\Phi=\varphi(p)\circ\psi:R[T_0,T_1,\ldots ,T_p]\to
\rees(\fm^p)$, the polynomial presentation of $\rees(\fm^p)$ sending
$T_i$ to $x^iy^{p-i}$ and let $Q=\ker(\Phi)$ be the ideal of equations
of $\rees(\fm^p)$. Let us see that $Q_n=\ker(\psi_n)+W_{n-1}Q_1$, for
all $n\geq 2$. Indeed, given $F\in Q_n$, since
$\ker(\varphi(p)_n)=\langle xY-yX\rangle V_{pn-1}$, one can find
$F_i\in W_{n-1}$ and $G_i\in Q_1$ such that $\psi _n(F)=\psi _n(\sum
F_iG_i)$. Therefore, $(F-\sum F_iG_i)\in \ker(\psi_n)$ and
$F\in\ker(\psi_n)+W_{n-1}Q_1$.

Since $\ker(\psi)=\ker(\psi)\langle 2\rangle$, then
$Q_n=\ker(\psi_n)+W_{n-1}Q_1=W_{n-2}\ker(\psi_2)+W_{n-1}Q_1\subset
W_{n-2}Q_2\subset Q_n$, for all $n\geq 2$. Therefore $Q=Q\langle
2\rangle$, $\reltype(\fm^n)=2$ and $\mu(Q_2/W_1Q_1)\leq\binom{p}{2}$.
Consider the diagram

\begin{equation*}
\begin{CD}
0 @>>> \ker\psi_2 @>>> W_2 @>\psi_2>> V(p)_2 @>>> 0\\
&& @VVV @VVV @VVV \\
&& \ker\psi_2 \otimes k @>>> W_2\otimes k @>\psi_2\otimes {\bf 1}_k>> V(p)_2\otimes k @>>> 0
\end{CD}
\end{equation*}

Let us prove that $Q_2/W_1Q_1$ is minimally generated by
$\binom{p}{2}$ elements, which are precisely the classes of the
$2\times 2$ minors of $\cm$. Since $Q_1\subset \fm W_1$, then
$W_1Q_1\subset \fm W_2$. Setting $L=\ker(\psi_2)$, $M=W_2$ and
$N=V(p)_2$, we have
\begin{eqnarray*}
\frac{Q_2}{W_1Q_1}\otimes k\cong \frac{Q_2}{\fm
  Q_2+W_1Q_1}=\frac{\ker(\psi_2)+W_1Q_1} {\fm\ker(\psi_2)+W_1Q_1}
\cong \frac{L}{\fm L+(L\cap W_1Q_1)}.
\end{eqnarray*}
On the other hand, there is a natural epimorphism 
\begin{eqnarray*}
\frac{L}{\fm L+(L\cap W_1Q_1)}\to \frac{L}{L\cap \fm M},
\end{eqnarray*}
where $L/(L\cap\fm M)\cong (L+\fm M)/\fm M=\ker(\psi_2\otimes
1_k)$. Hence, 
\begin{eqnarray*}
\mu\left( \frac{Q_2}{W_1Q_1}\right) \geq \dim \frac{L}{(L\cap \fm
  M}=\dim \ker(\psi_2\otimes 1_k),
\end{eqnarray*}
which clearly is $\binom{p}{2}$. Therefore $\mu((Q/Q\langle
1\rangle)_2)=\binom{p}{2}$ and, by
Remark~\ref{min-gen-set-of-equations}, there are $\binom{p}{2}$
equations of degree 2 in a minimal generating set of equations of
$\rees(\fm^p)$.

Note that $J=(x_1^p,x_2^p)$ is a reduction of $\fm^p$ with reduction
number $\rednumber_J(\fm^p)=1$, that $x_1^p,x_2^p$ is an $R$-sequence
and that $(x_1^p)^*$ is a $\graded(\fm^p)$-sequence.
\qed
\chapter{The injectivity of the canonical blowing-up morphisms}\label{article1}

\section{Introduction}\label{introduction}

Let $R$ be a commutative ring and let $I=(x_1,...,x_s,y)$ be an 
ideal of $R$. Let $\alpha_{I}:\symmetric(I)\rightarrow \rees(I)$ 
be the canonical morphism from the symmetric algebra of $I$ to the Rees algebra of
$I$. We will write $\alpha_{I,p}$ to denote the $p$-th graded component of 
$\alpha_I$. The aim of this chapter is to show that if $I$ is an ideal of 
almost-linear type, the injectivity of a single graded component 
$\alpha_{I,p}$ implies the injectivity of the lower graded 
components.

In \cite[Example~1.4]{kuhl} (see Example~\ref{exemple-kuhl}), 
K\"uhl gave an example of a finitely generated ideal $I$ with 
$\alpha_{I,n}$ being an isomorphism for $n$ sufficiently large, 
in particular, 
$\widetilde{\alpha}_{I}:{\rm Proj}(\rees(I))\to {\rm Proj}(\symmetric(I))$ 
is an isomorphism of schemes (see Remark~\ref{p-linear-type-giral}), 
but such that $\alpha_{I}$ is not an isomorphism of $R$-algebras. 
On the other hand, \cite{tchernev} asked whether $\alpha_{I,p}$ being an 
isomorphism implies that $\alpha_{I,n}$ is an isomorphism for each 
$2\leq n\leq p$. 

The main purpose of this chapter is to prove the following result.

\begin{Theorem B}
Let $R$ be a commutative ring, let $I=(x_1,\ldots,x_s,y)$ 
be an ideal of $R$ and let $p\geq 2$ be an integer. Suppose that the ideal 
$J=(x_1,\ldots ,x_s)$ verifies that $\alpha_{J,n} :\symmetric_n(J)\rightarrow J^n$ 
is an isomorphism for all $2\leq n\leq p$. 
Then the following conditions are equivalent:
\begin{enumerate}
\item[$(i)$] $\alpha_{I,p}:\symmetric_p(I)\to I^p$ is an isomorphism;
\item[$(ii)$] $\alpha_{I,n} :\symmetric_n(I)\rightarrow I^n$ is an
isomorphism for each $2\leq n\leq p$.
\end{enumerate}
\end{Theorem B}

Notice that the ideals of almost-linear type fulfil the hypotheses of Theorem~B, 
hence the equivalence of $(i)$ and $(ii)$ holds true for such ideals.

\medskip

For each integer $p\geq 2$, we also display an example of an ideal $I$ such
that $\alpha_{I,n}$ is an isomorphism for all $n\geq p+1$, whereas
$\alpha_{I,p}$ is not.

\begin{summary:exemple-article1}{\rm Let $k$ be a field and let $p\geq
2$. Let $S=k[U_0,\ldots,U_p,X,Y]$ be a polynomial ring and let $Q$ be
the ideal of $S$ defined as $Q=Q_1 + (U_0X^p)$, where
\begin{eqnarray*}
Q_1=(U_0Y,U_0X-U_1Y,U_1X-U_2Y,\ldots, U_{p-1}X-U_pY, U_pX).
\end{eqnarray*}
Let $R$ be the factor ring $S/Q=k[u_0,\ldots,u_p,x,y]$ and consider
the ideal $I=(x,y)\subset R$. Then $\alpha_{I,n}$ is an isomorphism
for all $n\geq p+1$, whereas $\alpha_{I,p}$ is not.}
\end{summary:exemple-article1}

\medskip

As in Section~\ref{sec:blowing-up-algebras}, we denote by $\beta_I$ and $\gamma_{I,\mathfrak{m}}$ 
the blowing-up morphisms $\alpha\otimes \textbf{1}_{R/I}:\symmetric(I/I^2)\to \graded(I)$ and 
$\alpha\otimes \textbf{1}_{R/\mathfrak{m}}:\symmetric(I/\mathfrak{m}I)\to\fiber(I)$, respectively. 
In order to slim the notations, we will write $\alpha$, $\beta$ and $\gamma$ when the ideals 
involved are clear from the context. Let us begin with some basic remarks about the influence of 
the injectivity of a single graded component of the canonical blowing-up morphisms on the rest of 
the graded components.

Notice first that in the case of $\gamma_{I,\mathfrak{m}}:\symmetric(I/\fm I)\to \fiber_\mathfrak{m}(I)$ 
the vanishing of a single component of $\ker\gamma_{I,\mathfrak{m}}$ imposes a tight condition to the lower
graded components, for if $\gamma_{I,\mathfrak{m},p+1}$ is an isomorphism, then
$\gamma_{I,\mathfrak{m},p}$ is also an isomorphism. In fact, if there
is a non-zero element $a\in\ker(\gamma_{I,\mathfrak{m},p})$, since
every non-zero element in $\symmetric_1(I/\fm I)$ is a non-zerodivisor 
of $\symmetric(I/\fm I)$, we can pick a non-zero element $\omega$ in 
$\symmetric_1(I/\fm I)$, such that $\omega\cdot a\in\ker(\gamma_{I,\mathfrak{m},p+1})$ 
is non-zero.

Considering this propagation of the injectivity of
the graded components of $\gamma_{I,\mathfrak{m}}$ 
towards the lower components, it is natural to ask 
for the same property for the morphisms $\alpha_I$ and 
$\beta_I$. In the case of $\alpha_I$ it turns out to not
hold in general. Indeed, in \cite[Example~1.4]{kuhl}, 
K\"uhl constructed an ideal $I$ with $\alpha_{I,p}$ 
being an isomorphism for $p\geq 3$ while $\alpha_{I,2}$ 
was not (see Example \ref{exemple-kuhl}). However, under 
certain settings the injectivity of the graded 
components of $\alpha$ propagates downwards, 
as Tchernev proved in \cite{tchernev}: 
if $\alpha_{I,p}$ is an isomorphism, then $\alpha_{I,n}$ 
is an isomorphism for $2\leq n\leq p$ when either ${\rm pd}_R(I)\leq 1$, 
or $I$ is perfect Gorenstein of grade $3$, 
or $I$ is a maximal ideal (see \cite[Section 5]{tchernev}). It is worth noting that conditions 
of this form ($\alpha_{n}$ being an isomorphism for all $2\leq n\leq p$) 
have already been considered across the commutative algebra landscape 
(see, for instance, \cite[Introduction, page~754]{concavalla}, 
\cite[Theorem~3]{costa1} or \cite[Definition~1.9 and Exercise~5.84]{vasconcelos}).

Theorem~B provides a new class of ideals $I$ for
which this pattern of the injectivity of $\alpha_{I}$ holds. Concretely, 
if $I=(J,y)$ with $\alpha_{J,n}:\symmetric_n(J)\to J^n$ an isomorphism 
for each $2\leq n\leq p$, then $\alpha_{I,p}:\symmetric_p(I)\to I^p$ 
being an isomorphism implies that $\alpha_{I,n}:\symmetric_n(I)\to I^n$ 
is an isomorphism for each $2\leq n\leq p$. In particular, if $J$ is an 
ideal of linear type, i.e., if $I$ is an ideal of almost-linear type, 
the hypothesis is fulfilled. Remark that writing 
an ideal $I$ as $I=(J,y)$ and assuming hypotheses over $J$ and the relation
between $J$ and $y$, is a common approach in this context
(see e.g. \cite{costa2}, \cite[Theorem~4.7]{hmv},
\cite[Proposition~3.9]{hsv}, \cite[Theorem~2.3 and
Proposition~2.5]{valla}). On the other hand, the study of the 
morphisms $\alpha$, $\beta$ and $\gamma$ has attracted a great deal of attention so 
far (see, for instance, \cite{costa1}, \cite{costa2}, \cite{hmv}, \cite{hsv}, 
\cite{huneke2}, \cite{huneke3}, \cite{johnsonmacloud},
\cite{keel}, \cite{kuhl}, \cite{micali}, \cite{tchernev}, 
\cite{valla}, \cite{vasconcelosblue} and \cite{vasconcelos}).

In Example~\ref{exemple} we will give a family of ideals 
$\{I_p\}_{p\geq 2}$ such that $\alpha_{I_p,n}$ is an isomorphism 
for each $n\geq p+1$, whereas $\alpha_{I_p,p}$ is not. 
Furthermore, $\beta_{I_p,n}$ will be an isomorphism if and only if 
$n\neq p$. More explicitely, for each positive integer $p$ we construct 
a $k[X,Y]$-algebra $\mathcal{A}_p$ and a $2$-generated ideal 
$I_p\subset\mathcal{A}_p$ such that $\beta_{I_p,n}$ is an isomorphism if and only if 
$n\neq p$. Furthermore, the example admits a ``glueing'',
in the sense that it is possible to construct, for any finite subset
$S\subset \mathbb{N}$, an ideal $I_S$ of a $k[X,Y]$-algebra $\mathcal{A}_S$ 
such that $\beta_{I_S,n}$ is an isomorphism if and only if $n\notin S$, 
by tensor-multiplying the collection of above defined $k[X,Y]$-algebras 
$\mathcal{A}_i$, for all $i\in S$. Beyond extending the example given by K\"uhl 
(focused on the case $p=2$), our example illustrates the fact that the morphism 
$\beta$ may have only one single non-injective component.

Many of the results presented in this chapter can also be found in the published work 
\cite{mp}: F.~Muiños, F.~Planas-Vilanova. \emph{On the injectivity of blowing-up ring 
morphisms}. J. Algebra {\bf 320} (2008), no. 8, 3365--3380.

\section{Injectivity does not propagate upwards}

In general, one cannot expect the injectivity of a single component 
of $\alpha$, $\beta$ or $\gamma$ to propagate upwards, i.e., under mild 
conditions, having $\ker\alpha_p=0$ (same for $\beta$ and $\gamma$) 
does not imply $\ker\alpha_n=0$ for $n>p$ (and similarly for $\ker\beta_n$ and $\ker\gamma_n$). 

Let us take $R=k[[X,Y,Z]]$ and $\fp\subset R$ be a prime ideal belonging 
to the family defined in \cite{moh} with $\mu(\fp)=r>3$. 
Such an ideal verifies that ${\rm pd}_R(\fp)=1$, $\height(\fp)=2$ and 
$\fp$ is generically a complete intersection ($\fp R_{\mathfrak{p}}$ is 
generated by a regular sequence, since $R_{\mathfrak{p}}$ is regular). 
Using \cite[Proposition~2.7 and Remark thereof]{hsv} one can conclude that $\fp$ 
is syzygetic (i.e., $\ker\alpha_2=0$). The analytic spread $\spread(\fp)$ satisfies 
$\spread(\fp)\leq\dim(R)=3<r$. For $n$ large enough, $\kdim{\fp^n/\fm\fp^n}$ 
is a polynomial in $n$ of degree $\spread(\fp)-1$, the Hilbert polynomial of 
$\fiber_{\mathfrak{m}}(I)$, whilst $\kdim{{\bf S}_n(\fp/\fm\fp)}=\binom{n+r-1}{r-1}$
is a polynomial in $n$ of degree $r-1$. Therefore, for $n$ large
enough, $\gamma_{\mathfrak{p},\mathfrak{m},n}:{\bf S}_n(\fp/\fm\fp)\to
\fp^n/\fm\fp^n$ is not an isomorphism, and consequently, neither
$\alpha_{\mathfrak{p},n}$ nor $\beta_{\mathfrak{p},n}$ are isomorphisms, 
whereas $\alpha_{\mathfrak{p},2}$, $\beta_{\mathfrak{p},2}$ and
$\gamma_{\mathfrak{p},\mathfrak{m},2}$ are isomorphisms.

\begin{Remark}{\rm
Recall that an ideal is projectively of linear type (of $p$-linear type, for short) if and only if
$\ker\alpha_{I,n}=0$ for $n$ sufficiently large (see claim~\ref{p-linear-type-giral}). 
Then if $I$ is of $p$-linear type it is a necessary condition that 
$\ker\gamma_n$ vanish for large $n$, then it follows that 
$\spread(I)=\mu(I)$, i.e., the second analytic deviation of $I$
must be equal to zero. However, the converse is not true unless 
we consider some extra conditions. For instance, consider the 
classical example $I=(X^3,X^2Y,Y^2Z)\subset k[[X,Y,Z]]$ (see \cite[p.16]{hmv}): 
$\mu(I)=\spread(I)=3$, while the class of $ZT_2^2-XT_1T_3\in\ker\alpha_{I,2}$ 
induces a non-zero element of $\ker\alpha_{I,n}$,
for each $n\geq 2$ (e.g., multiply repeatedly by $T_2$). 
}\end{Remark}

The proof of Theorem~B is based on the characterization of the 
vanishing of the effective relations in terms of the vanishing 
of the Koszul homology, which is in turn translated into a set 
of conditions in terms of colon ideals (see previous 
Chapter~\ref{article2} and \cite{planas}).

Sections~\ref{presentacio} and~\ref{effective} are
devoted to recall and refine some basic results concerning $\ker\alpha$
and the module of effective relations. In Sections~\ref{first} and 
\ref{second} we prove preliminary propositions which will be assembled 
in Section~\ref{mainresults}, where the main result is shown.

\section{The presentation of 
\texorpdfstring{$\ker\alpha$}{the kernel of alpha}}\label{presentacio}

The following lemma provides us with a presentation of 
the kernel of $\alpha$.

\begin{Lemma}\label{keralpha}
Let $R$ be a commutative ring and let $I=(x_1,\ldots ,x_s)$ be a finitely
generated ideal. Let $V=R[X_1,\ldots ,X_s]$ be a polynomial ring
and let $\varphi:V\to \rees(I)$ be the induced surjective graded 
morphism sending $X_{i}$ to $x_{i}t$. Let 
$\alpha_{I}:\symmetric(I)\to\rees(I)$ be the canonical morphism. 
Then, for $n\geq 2$,
\begin{eqnarray*}
{\rm ker}(\alpha_{I,n})\cong\frac{{\rm ker}(\varphi_{n})}{V_{n-1}\, {\rm
  ker}(\varphi_{1})}.
\end{eqnarray*}
\end{Lemma}

\demo Let $S(\varphi_1):V\rightarrow \symmetric(I)$ be the induced
surjective graded morphism sending $X_{i}$ to $x_{i}$ in degree one. We have 
the following commutative diagram:
\begin{equation*}
\xymatrix{
V \ar@{->>}[rd]^\varphi \ar@{->>}[d]_{S(\varphi_1)}\\
\symmetric(I) \ar@{->>}[r]_\alpha & \rees(I).
}
\end{equation*}

The resulting isomorphism is a consequence of relating the kernels of the
commutative diagram in the short exact sequence 
$0 \to \ker S(\varphi_1) \to \ker\varphi \to \ker\alpha \to 0$.
\qed

\begin{Example}\label{casprincipal}{\rm
Let $R$ be a commutative ring and let $I=(x)$ be a principal ideal. Let
$\alpha_{I}:\symmetric(I)\to\rees(I)$ be the canonical morphism. 
Then $\ker(\alpha_{I,n})\cong(0:x^n)/(0:x)$, for all $n\geq 2$.
In particular, one has the injective morphisms ${\rm
ker}(\alpha_{I,2})\hookrightarrow {\rm
ker}(\alpha_{I,3})\hookrightarrow \ldots \hookrightarrow {\rm
ker}(\alpha_{I,p})$. Consequently, if $\alpha_{I,p}$ is an isomorphism
for some $p\geq 2$, then $\alpha_{I}$ is an isomorphism.
}\end{Example}

\demo
Let $V=R[T]$ be the polynomial ring in one variable and let
$\varphi:V\to \rees(I)$ be the induced surjective graded morphism
sending $T$ to $xt$. Observe that
$\ker(\varphi_n)=(0:x^n)T^n$ and
$V_{n-1}\ker(\varphi_1)=RT^{n-1}(0:x)T=(0:x)T^n$. Now the assertion
follows by Lemma \ref{keralpha}. The last conclusion is a consequence 
of the rigidity of $\{(0:x^n)\}_{n\geq 1}$.
\qed

\section{Obstructions to injectivity: general settings}\label{effective}

Let $R$ be a commutative ring and let $I=(x_1,\ldots ,x_s)$ be a finitely
generated ideal. Let $\alpha_{I}:\symmetric(I)\to\rees(I)$ be the canonical
morphism. For $n\geq 2$, the {\em module of effective $n$-relations}
of $I$ is defined as $E(I)_n=\ker(\alpha_n)/\symmetric_1(I)\cdot\ker(\alpha_{n-1})$.
It is not difficult to interpret the effective relations as the fresh
relations of $\ker\alpha$, i.e., those relations which are not combination 
of the relations of lower degree. Recall that for the purpose of this work, 
as we did in Chapter~\ref{article2}, it will be useful to describe 
$E(I)_n$ in terms of a graded Koszul homology (see \cite{planas}):
$$E(I)_n \cong \koszul{1}{x_1t,\ldots,x_st}{\rees(I)}_n.$$

\begin{Remark}{\rm 
Let $V=R[X_1,\ldots ,X_s]$ be a polynomial ring over $R$, 
let $\varphi:V\to \rees(I)$ be the polynomial presentation 
sending $X_{i}$ to $x_{i}t$ and let $Q=\ker\varphi$. 
By Lemma~\ref{keralpha} and using $Q\langle 1\rangle_n=V_{n-1}\ker(\varphi_{1})\subset 
V_{1}\ker(\varphi_{n-1})=Q\langle n-1\rangle_n$ for $n\geq 2$, it
follows that $E(I)_n$ can be calculated as
\begin{eqnarray*}
E(I)_{n}\cong \frac{{\rm ker}(\varphi_{n})}{V_{1}\,{\rm
  ker}(\varphi_{n-1})} = \left[\frac{Q}{Q\langle n-1\rangle}\right]_n.
\end{eqnarray*}
}\end{Remark}

The following straightforward remark will be useful throughout:

\begin{Remark}\label{effectivealpha}\rm
Let $R$ be a commutative ring and let $I=(x_1,\ldots ,x_s)$ be a finitely
generated ideal. Let $\alpha_{I}:\symmetric(I)\to\rees(I)$ be the canonical
morphism and let $p\geq 2$ be an integer. Then, the following two conditions 
are equivalent:
\begin{enumerate}
\item[$(i)$] $\alpha_{I,n}:\symmetric_n(I)\to I^n$ is an isomorphism for
each $2\leq n\leq p$;
\item[$(ii)$] $E(I)_{n}=0$ for every $2\leq n\leq p$.
\end{enumerate}
\end{Remark}

In the particular cases of principal ideals and two-generated ideals,
the next two lemmas refine \cite[Example~3.1]{planas} and
\cite[Proposition~4.5]{planas}, respectively. These two lemmas clearly 
follow from Proposition \ref{succexcurta-RL(I)}, although here we 
provide alternative proofs.

\begin{Lemma}\label{lemma:effective-principal}
Let $R$ be a commutative ring and let $I=(x)$ be a principal ideal. For $n\geq
2$, there is an isomorphism $E(I)_{n}\cong(0:_R x^{n})/(0:_R x^{n-1})$.
\end{Lemma}

\demo
Let $V=R[T]$ and $\varphi:V\to\rees(I)$ be the surjective graded
morphism sending $T$ to $xt$. Note that
$\ker(\varphi_n)=(0:_R x^n)T^n$ and $V_1\ker(\varphi_{n-1})=(0:_R x^{n-1})T^n$. 
Since $E(I)_n\cong\ker(\varphi_n)/V_1\ker(\varphi_{n-1})$, 
the assertion follows.
\qed

\begin{Lemma}\label{exactsequence}
Let $R$ be a commutative ring and let $I=(x,y)$ be a two-generated ideal. For
every $n\geq 2$, there exists a short exact sequence of $R$-modules
\begin{eqnarray*}
0\to \frac{(0:_R x)\cap I^{n-1}}{y((0:_R x)\cap I^{n-2})}\to E(I)_{n}\to
\frac{(xI^{n-1}:_R y^n)}{(xI^{n-2}:_R y^{n-1})} \to 0.
\end{eqnarray*}
\end{Lemma}
\demo
Let $V=R[X,Y]$ be the polynomial ring in two variables and let
$\varphi:V\to\rees(I)$ be the surjective graded morphism mapping $X$
to $xt$ and $Y$ to $yt$, so that $E(I)_{n}\cong{\rm
ker}(\varphi_{n})/V_{1}{\rm ker}(\varphi_{n-1})$. Let
$$\frac{(0:_R x)\cap I^{n-1}}{y((0:_R x)\cap I^{n-2})}\buildrel
f\over\longrightarrow E(I)_{n}$$ be defined in the following way: for
$a\in (0:_R x)\cap I^{n-1}$, writing $a=F(x,y)$ with $F(X,Y)\in V_{n-1}$,
$f$ maps the class of $a$ to the class of $XF(X,Y)$ in
$\ker(\varphi_n)/V_1\ker(\varphi_{n-1})$. Let
$$E(I)_{n}\stackrel{g}{\longrightarrow}
\frac{(xI^{n-1}:_R y^n)}{(xI^{n-2}:_R y^{n-1})}$$ be defined as follows:
for $G(X,Y)\in \ker(\varphi_{n})$, $g$ maps the class of $G$ to
the class of $G(0,1)$. It is easily checked that $f$ and $g$ are
well-defined morphisms and determine a short exact sequence of
$R$-modules.
\qed

\begin{Discussion}\label{recall-exact-sequence}{\rm
Notice that the two previous lemmas could have been proved using
Proposition~\ref{succexcurta-RL(I)}. Provided $\mu(I)=1$ or $2$, 
the homological conditions can be effortlessly expressed in terms 
of colon ideals. However, in general, we are not able to produce a 
strightforward conversion into explicit ideal theoretic expressions. 
Recall that in general we have a short exact sequence as in 
Lemma~\ref{succexcurta}:
\begin{equation*}\label{succexcurta-article1}
0\to \frac{\koszul{1}{x_1t,\ldots,x_st}{\rees(I)}_n}{yt\koszul{1}{x_1t,\ldots,x_st}{\rees(I)}_{n-1}}
  \longrightarrow  E(I)_n\stackrel{\sigma_n}
         {\longrightarrow} \frac{(JI^{n-1}:_R y^n)}{(JI^{n-2}:_R y^{n-1})} \to
         0,
\end{equation*}
where $J=(x_1,\ldots,x_s)$ and $I=(J,y)$.
}\end{Discussion}

The following result (see \cite[Corollary~4.8]{planas}), 
will be fundamental as regards our approach to the proof 
of Theorem~B. In order to state it properly, let us 
introduce the following notation: we will denote by 
$I_{i_1,...,i_r}$ the ideal generated by the $x_j$ with 
$j\notin\{i_1,...,i_r\}$.

\begin{Proposition}\label{cambridge}
Let $I=(x_1,...,x_s)$ be an ideal generated by $s\geq
3$ elements and let $p\geq 2$. Suppose that $E(I_s)_p=0$. Then
$E(I)_p=0$ if and only if the following two conditions hold:
\begin{align*}
(a) & \hskip 10pt O_1(I,\underline{x},p,i):= \frac{(I_iI^{p-1}:_R x_i^p)}{(I_iI^{p-2}:_R x_i^{p-1})}=0,\; \textrm{for each } 1\leq i\leq s;\\
(b) & \hskip 10pt O_2(I,\underline{x},p):= \frac{((\sum_{1\leq i<j\leq
s-1}{x_ix_jI_s^{p-2}}) :_R x_s)\cap I^{p-1}}{
\sum_{i=1}^{s-1}{x_i((I_{i,s}I_s^{p-2}:_R x_s)\cap I^{p-2})}}=0.
\end{align*}
\end{Proposition}

\medskip

When $s\geq 3$, the vanishing of $E(I)_{p}$ can be characterized 
by means of three ideal theoretic explicit conditions. However, under the 
assumption that $E(I_s)_p=0$, one of the conditions is 
satisfied trivially and only $(a)$ and $(b)$ in 
Proposition \ref{cambridge} survive (see \cite[Theorem~4.7 and
Corollary~4.8]{planas}).

Now consider $I=(x_1,...,x_s)$ and suppose that
$J=I_s=(x_1,\ldots,x_{s-1})$ is such that
$\alpha_{J,n}:\symmetric_n(J)\to J^n$ is an isomorphism for $2\leq n\leq
p$, as in the hypothesis of Theorem~B. By
Remark~\ref{effectivealpha}, $E(J)_{n}=0$ for every $2\leq n\leq
p$. Therefore, and again by Remark~\ref{effectivealpha}, to prove the
implication $(i)\Rightarrow (ii)$ in Theorem~B, it
will be enough to fulfil obstructions $(a)$ and $(b)$ in
Proposition~\ref{cambridge} for every $2\leq n\leq p$. Each of
these conditions will be studied separately in the next sections.

\section{Vanishing of the first obstruction}\label{first}

The purpose of this section is to study the vanishing of the condition 
$(a)$ in Proposition~\ref{cambridge}. The following proposition
makes use of an explicit construction.

\begin{Proposition}\label{surjective}
Let $R$ be a commutative ring, let $I=(x_1,...,x_s)$ be a finitely generated
ideal and let $p\geq 2$ be an integer. Let
$\alpha_{I}:\symmetric(I)\to\rees(I)$ be the canonical morphism. Then, for
each $1\leq i\leq s$, there exists a surjective morphism
\begin{eqnarray*}
g_{p,i}:{\rm ker}(\alpha_{I,p})\twoheadrightarrow
\frac{(I_{i}I^{p-1}:_R x_{i}^{p})}{(I_{i}:_R x_{i})}.
\end{eqnarray*}
In particular, if $\alpha_{I,p}$ is an isomorphism, then
\begin{eqnarray*}
\frac{(I_iI^{n-1}:_R x_i^n)}{(I_iI^{n-2}:_R x_i^{n-1})}=0
\end{eqnarray*}
for each $1\leq i\leq s$ and for each $2\leq n\leq p$.
\end{Proposition}

At this point let us introduce some useful notations. For
$\sigma=(\sigma_1,...,\sigma_s)\in\mathbb{N}^s$ set
$\xsigma=x_1^{\sigma_1}...x_s^{\sigma_s}$ and
$\Xsigma=X_1^{\sigma_1}...X_s^{\sigma_s}$. Define the {\it support of}
$\sigma\in\mathbb{N}^s$ as the subset ${\rm supp}(\sigma)\subset
\{1,\ldots,s\}$ such that $i\in{\rm supp}(\sigma)$ if and only if
$\sigma_i\neq 0$. Let $|\sigma|$ denote $\sum_{i=1}^s{\sigma_i}$. We
prove now Propostion~\ref{surjective}.

\demo Let $V=R[X_1,\ldots ,X_s]=R[{\bf X}]$ be the polynomial ring
in $s$ indeterminates over $R$ and let $\varphi:V\to \rees(I)$ be 
the induced surjective graded morphism sending $X_{i}$ to $x_{i}t$. 
Using Lemma~\ref{keralpha}, present ${\rm ker}(\alpha_{I,p})$ as
the quotient ${\rm ker}(\varphi_{p})/V_{p-1}{\rm ker}(\varphi_{1})$.
Let $g_{p,i}:{\rm ker}(\alpha_{I,p})\rightarrow
(I_{i}I^{p-1}:_R x_{i}^{p})/(I_{i}:_R x_{i})$ be defined as follows: for
$H({\bf X})\in{\rm ker}(\varphi_{p})$, $g_{p,i}$ maps the class of
$H({\bf X})$ to the class of $H(e_{i})$ in
$(I_{i}I^{p-1}:_R x_{i}^{p})/(I_{i}:_R x_{i})$, where $e_{i}$ stands for the
$i$-th vector of the natural basis of $R^{s}$. To see that $g_{p,i}$
is well-defined, take $H({\bf X})\in V_{p-1}\ker(\varphi_1)$. Since
any element of $V_{p-1}$ can be written as a linear combination of the
natural basis given by the monomials of degree $p-1$, one can write
$H({\bf X})=\sum_{|\sigma|=p-1}(a_{\sigma,1}X_1+\ldots
+a_{\sigma,s}X_s)\Xsigma$, where
$a_{\sigma,1}x_1+\ldots+a_{\sigma,s}x_s=0$, $a_{\sigma,j}\in R$. Then
$H(e_i)=a_{(p-1)e_i\, ,i}$, which belongs to $(I_i:_R x_i)$. It is
readily seen that $g_{p,i}$ is surjective. Finally, since
$(I_i:_R x_i)\subset (I_iI:_R x_i^2)\subset\ldots\subset
(I_iI^{p-1}:_R x_i^p)$, if $\alpha_{I,p}$ is an isomorphism, then
$(I_iI^{n-1}:_R x_i^n)/(I_{i}I^{n-2}:_R x_{i}^{n-1})=0$ for each $2\leq
n\leq p$.  \qed

\section{Theorem B with stronger hypotheses}

In this section we state and prove a version of Theorem~B that works 
under more restrictive assumptions than the version that we will prove 
later. The special version uses key results from Chapter~\ref{article2} 
and can be proved with little effort, whereas the more general version 
will be more demanding in terms of preparatory results.

\begin{Theorem}\label{maintheorem-weak}
Let $(R,\fm)$ be a Noetherian local ring, let $I=(x_1,\ldots,x_s,y)$ be an ideal of $R$. 
Suppose that $x_1,\ldots,x_s$ verify the following condition, for all $n\geq 2$:
\[
((x_1,\ldots,x_{i-1})I^{n-1}:_R x_i)\cap I^{n-1}=(x_1,\ldots,x_{i-1})I^{n-2},\;\; \textrm{for all } i=1,\ldots,s.
\label{trung-module} \tag{$\mathcal{T}_n$}
\]
Then the following conditions are equivalent:
\begin{enumerate}
\item[$(i)$] $\alpha_{I,p}:\symmetric_p(I)\to I^p$ is an isomorphism;
\item[$(ii)$] $\alpha_{I,n}:\symmetric_n(I)\to I^n$ is an isomorphism, 
for each $2\leq n\leq p$;
\end{enumerate}
\end{Theorem}
\demo
Let $J=(x_1,\ldots,x_s)$. Recall from Lemma~\ref{characterisation-of-homology} that the 
conditions in (\ref{trung-module}) hold if and 
only if $\koszul{1}{x_1t,\ldots,x_it}{\rees(I)}_{n}=0$, 
for all $i=1,\ldots,s$. Then using the exact sequence 
in Discussion~\ref{recall-exact-sequence} we get the 
isomorphisms: $$E(I)_n\cong \frac{(JI^{n-1}:_R I^n)}{(JI^{n-2}:_R I^{n-1})},
\;\textrm{for all } 2\leq n\leq p.$$ 
Since $\ker\alpha_{I,p}=0$, from Proposition~\ref{surjective} 
we know that $E(I)_n=0$, for all $2\leq n\leq p$. 
The claim follows by Remark~\ref{effectivealpha} .\qed

\begin{Corollary}
Let $(R,\fm)$ be a Noetherian local ring, let $I=(x_1,\ldots,x_s,y)$ be an ideal of $R$ 
with $J=(x_1,\ldots,x_s)$ a reduction of $I$, $x_1,\ldots,x_s$ 
an $R$-sequence and $x_1^*,\ldots,x_{s-1}^*$ a $\graded(I)$-sequence. 
Then the following conditions are equivalent:
\begin{enumerate}
\item[$(i)$] $\alpha_{I,p}:\symmetric_p(I)\to I^p$ is an isomorphism;
\item[$(ii)$] $\alpha_{I,n}:\symmetric_n(I)\to I^n$ is an isomorphism, 
for each $2\leq n\leq p$;
\end{enumerate}
\end{Corollary}
\demo
The result follows directly from Lemma \ref{homologia-zero} and 
Theorem \ref{maintheorem-weak}.
\qed

\section{Vanishing of the second obstruction}\label{second}

Let $I=(x_1,\ldots ,x_s)$ be a finitely generated ideal of $R$ and
let $\alpha_{I}:\symmetric(I)\to\rees(I)$ be the canonical morphism. Using
Lemma~\ref{keralpha} and its terminology, for each $p\geq 2$, present
${\rm ker}(\alpha_{I,p})$ as ${\rm ker}(\varphi_{p})/V_{p-1}{\rm
ker}(\varphi_{1})$. If $1\leq n\leq p-1$, since $V_{p-1}{\rm
ker}(\varphi_{1})\subset V_{p-n}{\rm ker}(\varphi_{n})$, one has a
natural surjective morphims:
\[ 
\ker(\alpha_{I,p})\cong\frac{{\rm ker}(\varphi_{p})}{V_{p-1}{\rm
ker}(\varphi_{1})}\stackrel{\pi_{p,n}}{\twoheadrightarrow} \frac{{\rm
ker}(\varphi_{p})}{V_{p-n}{\rm ker}(\varphi_{n})}.
\label{equation-ker-alpha-epi} \tag{$\Pi_{p,n}$}
\]
Set $J=(x_{1},\ldots ,x_{s-1})$ and suppose that
$\alpha_{J,n}:\symmetric_{n}(J)\to J^n$ is an isomorphism for each $2\leq
n\leq p$. In this section we will construct morphisms $f_{p,n}$, with
\[
O_2(I,\underline{x},n+1)=\frac{((\sum_{1\leq i<j\leq s-1}{x_ix_jI_s^{n-1}}) :_R x_s)\cap I^{n}}{
\sum_{i=1}^{s-1}{x_i((I_{i,s}I_s^{n-1}:_R x_s)\cap I^{n-1})}} \stackrel{f_{p,n}}{\longrightarrow}
          \frac{\ker \varphi_{p}}{V_{p-n} \ker\varphi_n},
\]
from the modules which appear in condition $(b)$ of
Proposition~\ref{cambridge} to ${\rm ker}(\varphi_{p})/V_{p-n}{\rm
ker}(\varphi_{n})$ (see Propositions~\ref{2injective} and
\ref{fnps}). Thus we will have the diagram:
\begin{equation*}
\small
\begin{CD}
& & \ker\alpha_{I,p} \\
& & @VV{\pi_{p,n}} V \\
O_2(I,\underline{x},n+1) @>f_{p,n}>> 
          \ker \varphi_{p}/V_{p-n} \ker\varphi_n.
\end{CD}
\end{equation*}
\normalsize

\medskip

We will show that, for $p\geq 2$ and $1\leq n\leq p-1$, $f_{p,n}$ is
injective provided that $I$ is two-generated and $(0:x)=(0:x^2)$, where
$I=(x,y)$ (see Proposition~\ref{2injective}) or, in the general case, if
$\alpha_{I,p}$ is an isomorphism (see Proposition~\ref{fnps} and the
proof of Proposition~\ref{c2zero}). In particular, if $\alpha_{I,p}$
is an isomorphism, the modules of the bottom row in the diagram above
are zero and the equations $(b)$ in Proposition~\ref{cambridge} hold.

\subsection*{\bf Two-generated case}

For the sake of a better comprehension, first we treat the particular
case of $I$ being a two-generated ideal. Remark that here, unlike in
the general case, we do not need the hypothesis ``$\alpha_{J,n}$
isomorphic'' to ensure that the morphisms $f_{p,n}$ are well-defined.

\begin{Proposition}\label{2injective}
Let $R$ be a commutative ring, let $I=(x,y)$ be a two-generated ideal and let
$p\geq 2$ be an integer. Let $V=R[X,Y]$ be the polynomial ring in
two variables and let $\varphi:V\to \rees(I)$ be the induced surjective
graded morphism sending $X$ to $xt$ and $Y$ to $yt$.
Then, for each $1\leq n\leq p-1$, there exists a
well-defined morphism
\begin{eqnarray*}
f_{p,n}:\frac{(0:y)\cap I^n}{x((0:y)\cap I^{n-1})}\longrightarrow
\frac{\ker(\varphi_p)}{V_{p-n}\ker(\varphi_n)}.
\end{eqnarray*}
Moreover, if $(0:x)=(0:x^{2})$, then $f_{p,n}$ is injective.
\end{Proposition}

\demo Given $a\in (0:y)\cap I^n$, $a=A(x,y)$ with $A(X,Y)\in V_n$, let
$f_{p,n}$ be the map which sends the class of $a$ to the class of
$Y^{p-n}A(X,Y)$ in $\ker(\varphi_p)/V_{p-n}\ker(\varphi_{n})$.

\medskip

\noindent {\bf Well-defined}. The definition of $f_{p,n}$ is
independent of the choice of $A(X,Y)$, since if $a=B(x,y)$ for another
$B(X,Y)\in V_n$, then $Y^{p-n}(A(X,Y)-B(X,Y))$ is in
$V_{p-n}\ker(\varphi_n)$. In addition, if $a\in x((0:y)\cap I^{n-1})$,
then $a=xK(x,y)$ where $K(X,Y)\in V_{n-1}$ and $yK(x,y)=0$. Thus the
image of the class of $a$ is the class of
$Y^{p-n}XK(X,Y)=Y^{p-n-1}XYK(X,Y)\in V_{p-n}\ker(\varphi_n)$.

\medskip

\noindent {\bf Injectivity of $f_{p,n}$}. Let
$a=\sum_{j=0}^n{\lambda_jx^{j}y^{n-j}}$ with $ay=0$. Suppose that
$f_{p,n}$ maps the class of $a$ to the zero class. Then we can
write
\begin{eqnarray*}
Y^{p-n}\Bigl(\sum_{j=0}^{n}{\lambda_j
X^{j}Y^{n-j}}\Bigr)=\sum_{r=0}^{p-n}G_r(X,Y)X^{r}Y^{p-n-r}\in V_{p-n}\ker\varphi_n,
\end{eqnarray*}
where $G_r(X,Y)\in\ker(\varphi_n)$. Set $G_r(X,Y)=\sum_{t=0}^n
u_{r,t}X^{t}Y^{n-t}$. Note that $G_r(X,Y)$ and the $u_{r,t}$ depend on $n$. 
Equating the coefficients
\begin{eqnarray*}
\lambda_j & = & \sum_{r+t=j}u_{r,t} \; \mbox{ for }\; 1\leq j\leq n \, ,
\end{eqnarray*}
where $0\leq r\leq p-n$ and $0\leq t\leq n$. Therefore,
\begin{eqnarray*}
a & = & \sum_{j=0}^n{\Bigl(\sum_{r+t=j} u_{r,t}\Bigr) x^{j}y^{n-j}} =
  u_{0,0}y^n + \sum_{j=1}^n{\Bigl(\sum_{r+t=j} u_{r,t}\Bigl)
  x^{j}y^{n-j}}.
\end{eqnarray*}
Since $G_0(x,y)=0$, $u_{0,0}y^n=-\sum_{j=1}^n{u_{0,j}x^jy^{n-j}}$. Then
\begin{eqnarray*}
a & = & \sum_{j=1}^n{\biggl(\Bigl(\sum_{r+t=j} u_{r,t}\Bigr) -
  u_{0,j}\biggr) x^{j}y^{n-j}} =\\ & = &
  x\sum_{j=1}^n{\biggl(\Bigl(\sum_{r+t=j} u_{r,t}\Bigr)
  -u_{0,j}\biggr) x^{j-1}y^{n-j}}=xz,
\end{eqnarray*}
with $z=\sum_{j=1}^n{\Bigl(\bigl(\sum_{r+t=j} u_{r,t}\bigr) -u_{0,j}\Bigr)
x^{j-1}y^{n-j}}\in I^{n-1}$. It remains to show that $z\in (0:y)$.
First, consider the case $n=1$. Note that in this case $z=u_{1,0}$. Then, since $n=1$, 
$u_{1,0}y=-u_{1,1}x$ and $0=ay=xzy=u_{1,0}xy=-u_{1,1}x^2$. Thus
$u_{1,1}\in (0:x^2)$, which by hypothesis is equal to
$(0:x)$. Therefore $u_{1,0}y=-u_{1,1}x=0$ and we conclude that
$a=xz=xu_{1,0}\in x(0:y)$, proving the injectivity of $f_{p,1}$. 

From now on, let $n\geq 2$. Then,
\begin{eqnarray*}
zy & = & \sum_{j=1}^n{\biggl(\Bigl(\sum_{r+t=j} u_{r,t}\Bigr)-u_{0,j}\biggr)
   x^{j-1}y^{n-j+1}}\\ & = & u_{1,0}y^n+
   \sum_{j=2}^n{\biggl(\Bigl(\sum_{r+t=j} u_{r,t}\Bigr)-u_{0,j}\biggr)
   x^{j-1}y^{n-j+1}}.
\end{eqnarray*}
Since $G_1(x,y)=0$, $u_{1,0}y^n=-\sum_{j=1}^n u_{1,j}x^jy^{n-j}$. Then,
\begin{eqnarray*}
zy & = & \sum_{j=2}^n{\biggl(\Bigl(\sum_{r+t=j}
   u_{r,t}\Bigr)-u_{0,j}\biggr) x^{j-1}y^{n-j+1}}-\sum
   _{j=1}^{n}u_{1,j}x^{j}y^{n-j}\\ & = & x \biggl[
   \sum_{j=2}^n{\biggl(\Bigl(\sum_{r+t=j} u_{r,t}\Bigr)-u_{0,j}\biggr)
   x^{j-2}y^{n-j+1}} -\sum _{j=1}^{n}u_{1,j}x^{j-1}y^{n-j}\biggr].
\end{eqnarray*}
Setting $w=\sum_{j=2}^n{\Bigl(\bigl(\sum_{r+t=j} u_{r,t}\bigr)-u_{0,j}\Bigr)
x^{j-2}y^{n-j+1}}-\sum _{j=1}^{n}u_{1,j}x^{j-1}y^{n-j}$, we have
$zy=xw$. Since $0=ay=xzy=x^2w$, $w\in(0:x^2)$ which, by hypothesis, is
equal to $(0:x)$. Thus $zy=xw=0$ and $a=xz\in x((0:y)\cap I^{n-1})$,
hence $f_{p,n}$ is injective.
\qed

\subsection*{\bf General case}

Now let $I=(x_{1},\ldots ,x_{s})$, $s\geq 3$, and set $J=(x_{1},\ldots
,x_{s-1})$. For $p\geq 2$ and $1\leq n\leq p-1$, we are going to
construct the morphisms $f_{p,n}$ under the assumption that
$\alpha_{J,n}:\symmetric_{n}(J)\to J^n$ is an isomorphism for $2\leq n\leq
p$. For convenience, we say that a form
$A(X_1,\ldots,X_s)\in V_n$ {\it represents} an
element $a\in I^n$ if $A(x_1,\ldots,x_s)=a$.

\begin{Proposition}\label{fnps}
Let $R$ be a commutative ring, let $I=(x_1,\ldots,x_s)$ with $s\geq 3$, 
and let $p\geq 2$ be an integer.  Let
$V=R[X_1,\ldots,X_s]$ be a polynomial ring and let
$\varphi:V\to \rees(I)$ be the induced surjective graded morphism
sending $X_i$ to $x_it$. Suppose that
$J=(x_{1},\ldots,x_{s-1})$ verifies that $\alpha_{J,n}:\symmetric_{n}(J)\to
J^n$ is an isomorphism for each $2\leq n\leq p$. Then, for each $1\leq
n\leq p-1$, there exists a well-defined morphism
\begin{eqnarray*}
f_{p,n}:O_2(I,\underline{x},n+1)=\frac{\Bigl(\bigl(\sum_{1\leq i<j\leq
s-1}{x_ix_jI_s^{n-1}}\bigr) :x_s\Bigr)\cap I^{n}}{
\sum_{i=1}^{s-1}{x_i\bigl(\bigl(I_{i,s}I_s^{n-1}:x_s\bigr)\cap
I^{n-1}\bigr)}} \longrightarrow\frac{\ker\varphi_p}{V_{p-n}\ker{\varphi_n}}.
\end{eqnarray*}
Moreover $f_{p,1}$ is injective.
\end{Proposition}

\demo Pick an element $$a\in\Bigl( \bigl(\sum_{1\leq i<j\leq
s-1}{x_ix_jI_s^{n-1}}\bigr) :x_s\Bigr)\cap I^{n}.$$ 
Thus, $x_sa=\sum_{1\leq i<j\leq s-1}{x_ix_ja_{ij}}$ with $a_{ij}\in
I_s^{n-1}$. Let $A({\bf X})\in V_n$ be a form representing $a$ and let
$A_{ij}({\bf X}')\in V_{n-1}$ be a form representing $a_{ij}$, where
${\bf X}'=\{X_1,\ldots,X_{s-1}\}$. Let $f_{p,n}$ be the map sending
the class of $a$ to the class of
\begin{eqnarray*}
X_s^{p-n-1}\Bigl(X_sA({\bf X})-\sum_{1\leq i<j\leq
s-1}{X_iX_jA_{ij}({\bf X}')}\Bigr)
\end{eqnarray*}
in $\ker(\varphi_p)/V_{p-n}\ker{\varphi_n}$.

\medskip

\noindent {\bf Well-defined:} Set $W=R[X_1,\ldots,X_{s-1}]$ and define
$\psi:W\to\rees(J)$ sending $X_{i}$ to $x_{i}t$. In
particular, the following diagram is commutative:

\begin{equation*}
\begin{CD}
W @>{\psi}>> \rees(J) \\ 
@VVV @VVV \\
V @>{\varphi}>> \rees(I).
\end{CD}
\end{equation*}

In order to check that $f_{p,n}$ is well-defined, we first prove that
the image of the class of $a$ does not depend on the choices of
$A({\bf X})$ and $A_{ij}({\bf X}')$. Let $A({\bf X}),B({\bf X})\in
V_n$ be two forms representing $a$. Suppose that for $a_{ij},
b_{ij}\in I_s^{n-1}$
\begin{eqnarray*}
x_sa=\sum_{1\leq i<j\leq s-1}{x_ix_ja_{ij}}=\sum_{1\leq i<j\leq
s-1}{x_ix_jb_{ij}}.
\end{eqnarray*}
Let $A_{ij}({\bf X}')$, $B_{ij}({\bf X}')\in W_{n-1}$ be two forms
representing $a_{ij}$ and $b_{ij}$, respectively. Then the class of
\begin{eqnarray*}
X_s^{p-n-1}\biggl(X_s\Bigl(A({\bf X})-B({\bf X})\Bigr)-\sum_{1\leq i<j\leq
s-1}{X_iX_j\Bigl(A_{ij}({\bf X}')-B_{ij}({\bf X}')\Bigr)}\biggr)
\end{eqnarray*}
in $\ker{\varphi_p}/V_{p-n}\ker{\varphi_n}$ is the zero class, since
$A({\bf X})-B({\bf X})\in\ker(\varphi_n)$ and
\begin{equation*}
\sum_{1\leq i<j\leq
s-1}{X_iX_j\bigl(A_{ij}({\bf X}')-B_{ij}({\bf X}')\bigr)}\in\ker(\psi_{n+1}),
\end{equation*}
which by hypothesis is equal to $W_1\ker(\psi_n)\subset
V_1\ker(\varphi_n)$.

Assume now that $a\in
\sum_{i=1}^{s-1}{x_i\bigl((I_{i,s}I_s^{n-1}:x_s)\cap I^{n-1}\bigr)}$.
Write $a=\sum_{i=1}^{s-1}{x_ia_i}$, where $a_i\in I^{n-1}$,
$x_sa_i=\sum_{ {j=1},{j\neq i}}^{s-1}{x_jc_{ij}}$ and $c_{ij}\in
I_s^{n-1}$. Hence if $A_i({\bf X})\in V_{n-1}$ is a form representing
$a_i$, then $\sum_{i=1}^{s-1} X_iA_i({\bf X})$ is a form representing
$a$. Note that $x_sa=\sum_{i=1}^{s-1}{x_i\sum_{{j=1},{j\neq
i}}^{s-1}{x_jc_{ij}}}$ which can be written as $\sum_{1\leq i<j\leq
s-1}x_ix_jd_{ij}$, where $d_{ij}=c_{ij}+c_{ji}$ for $1\leq i<j\leq
s-1$. Thus $f_{p,n}$ maps the class of $a$ to the class of
\begin{eqnarray*}
X_s^{p-n-1}\Bigl( X_s\sum_{i=1}^{s-1}{X_iA_i({\bf
X})}-\sum_{1\leq i<j\leq r-1}X_iX_jD_{ij}({\bf X}')\Bigr)=\\
= X_s^{p-n-1}\Bigl( X_s\sum_{i=1}^{s-1}{X_iA_i({\bf
X})}-\sum_{i=1}^{s-1}{X_i\sum_{{j=1}\atop{j\neq
i}}^{s-1}{X_jC_{ij}({\bf X}')}}\Bigr)=\\ =X_s^{p-n-1}
\sum_{i=1}^{s-1}{ X_i \Bigl( X_sA_i({\bf X}) - \sum_{{j=1} \atop
{j\neq i}}^{s-1}{X_jC_{ij}({\bf X}')} \Bigr) },
\end{eqnarray*}
where $C_{ij}({\bf X}')\in W_{n-1}$ is a form representing $c_{ij}$
and $D_{ij}({\bf X}')=C_{ij}({\bf X}')+C_{ji}({\bf X}')$. Since
$X_sA_i({\bf X}) - \sum_{{j=1},{j\neq i}}^{s-1}{X_jC_{ij}({\bf
X}')}\in\ker(\varphi_n)$, then $f_{p,n}$ maps the class of $a$ to
the zero class.

\medskip

\noindent {\bf Injectivity of $f_{p,1}$}. Let
$a\in\Bigl(\bigl(\sum_{1\leq i<j\leq s-1}{x_ix_jR}\bigr):x_s\Bigr)\cap
I$. Then $a=\sum_{i=1}^s{\lambda_ix_i}$, with $\lambda_i\in R$, and
$x_sa=\sum_{1\leq i<j\leq s-1}{x_ix_ja_{ij}}$, $a_{ij}\in R$. Suppose
that $f_{p,1}$ maps the class of $a$ to the zero class. Then we can
write
\begin{equation}\label{linear}
X_s^{p-2}\Bigl(X_s\sum_{i=1}^s{\lambda_iX_i}-\sum_{1\leq i<j\leq
s-1}{X_iX_ja_{ij}}\Bigr)=\\
\sum_{|\sigma|=p-1}{(u_{\sigma,1}X_1+\ldots+u_{\sigma,s}X_s)\Xsigma},
\end{equation}
with $u_{\sigma,1}x_1+\ldots+u_{\sigma,s}x_s=0$, $u_{\sigma,j}\in R$,
where recall that if
$\sigma=(\sigma_1,\ldots,\sigma_s)\in\mathbb{N}^s$,
$x^\sigma=x_1^{\sigma_1}\ldots x_s^{\sigma_s}$ and
$X^\sigma=X_1^{\sigma_1}\ldots X_s^{\sigma_s}$. For convenience, denote
\begin{displaymath}
u_{k,t} = \begin{cases}
u_{e_k+(p-2)e_s,t}& \mbox{ for }1\leq k\leq s-1, \\
u_{(p-1)e_s,t}& \mbox{ for }k=s.
\end{cases}
\end{displaymath}
where $\{e_1,\ldots,e_s\}$ stands for the canonical basis of
$\mathbb{N}^s$. With these notations and from equation (\ref{linear}),
we get the following identities:
\begin{displaymath}
\begin{cases}
\lambda_i=u_{s,i}+u_{i,s}& \mbox{ for }1\leq i\leq s-1, \\
\lambda_s=u_{s,s}. & 
\end{cases}
\end{displaymath}
Using the fact that $u_{s,1}x_1+...+u_{s,s}x_s=0$, we get
\begin{eqnarray*}
a & = & \sum_{i=1}^s{\lambda_ix_i}=
\sum_{i=1}^{s-1}{(u_{s,i}+u_{i,s})x_i}+u_{s,s}x_s=\\
& = & \sum_{i=1}^{s-1}{(u_{s,i}+u_{i,s})x_i}-\sum_{i=1}^{s-1}u_{s,i}x_i=
\sum_{i=1}^{s-1}{u_{i,s}x_i}.
\end{eqnarray*}
Since $u_{i,1}x_1+\ldots+u_{i,s}x_s=0$ for $1\leq i\leq s-1$, then
$u_{i,s}x_s=-\sum_{j=1}^{s-1}u_{i,j}x_j$ and
\begin{eqnarray*}
ax_s=\sum_{i=1}^{s-1}{u_{i,s}x_sx_i}=-\sum_{i,j=1}^{s-1}{u_{i,j}x_ix_j}.
\end{eqnarray*}
But $ax_s=\sum_{1\leq i<j\leq s-1}{x_ix_ja_{ij}}$. Hence
\begin{eqnarray}\label{R1}
R({\bf X})=\sum_{i,j=1}^{s-1}{u_{i,j}X_iX_j}+\sum_{1\leq i<j\leq
s-1}{X_iX_ja_{ij}}
\end{eqnarray}
belongs to $\ker(\psi_2)$. Since $\alpha_{J,2}$ is an isomorphism,
$R({\bf X})$ can be written as
\begin{eqnarray}\label{R2}
R({\bf X})= \sum_{i=1}^{s-1}{(w_{i,1}X_1+\ldots+w_{i,s-1}X_{s-1})X_i}
\end{eqnarray}
with $w_{i,1}x_1+\ldots+w_{i,s-1}x_{s-1}=0$, $w_{i,j}\in R$. In
particular $w_{i,i}\in (I_{i,s}:x_i)$. Equating the two expressions of
$R({\bf X})$ in (\ref{R1}) and (\ref{R2}) we get the identities
$u_{i,i}=w_{i,i}\in (I_{i,s}:x_i)$ for each $1\leq i\leq s-1$. Then,
for $1\leq i\leq s-1$, $u_{i,s}x_{s}=-\sum_{j=1}^{s-1}u_{i,j}x_{j}=
-\sum_{j=1,j\neq i}^{s-1}u_{i,j}x_{j}-u_{i,i}x_{i}\in I_{i,s}$ and
$u_{i,s}\in (I_{i,s}:x_s)$. Therefore
$a=\sum_{i=1}^{s-1}u_{i,s}x_i\in\sum_{i=1}^{s-1}{x_i(I_{i,s}:x_s)}$,
proving the injectivity of $f_{p,1}$.  \qed

The next proposition proves the vanishing of equations $(b)$ in
Proposition~\ref{cambridge} under our hypotheses.

\begin{Proposition}\label{c2zero}
Let $R$ be a commutative ring, let $I=(x_1,\ldots,x_s)$ with $s\geq 3$ 
and let $p\geq 2$ be an integer. Suppose that
$J=(x_1,...,x_{s-1})$ verifies that $\alpha_{J,n}:\symmetric_{n}(J)\to J^n$
is an isomorphism for each $2\leq n\leq p$. If $\alpha_{I,p}$ is an
isomorphism, then
\begin{eqnarray*}
O_2(I,\underline{x},n+1):=\frac{((\sum_{1\leq i<j\leq s-1}{x_ix_jI_s^{n-1}}) :x_s)\cap I^{n}}{
\sum_{i=1}^{s-1}{x_i((I_{i,s}I_s^{n-1}:x_s)\cap I^{n-1})}}=0 \;\;\;
{\rm for}\; 1\leq n\leq p-1.
\end{eqnarray*}
\end{Proposition}

\demo 
Let $V=R[X_1,...,X_s]$ be the polynomial ring in $s$ variables
and let $\varphi:V\to \rees(I)$ be the induced surjective graded
morphism sending $X_i$ to $x_it$. If $\alpha_{I,p}$ is an
isomorphism, using the natural surjective morphism 
(\ref{equation-ker-alpha-epi}) in the beginning of Section~\ref{second}, 
one deduces that
\begin{eqnarray*}
\frac{\ker(\varphi_p)}{V_{p-n}\ker(\varphi_n)}=0\; \mbox{ for each }1\leq
n\leq p-1.
\end{eqnarray*}
For $1\leq n\leq p-1$, consider the morphisms $f_{p,n}$ defined in
Proposition~\ref{fnps}:
\begin{eqnarray*}
f_{p,n}:O_2(I,\underline{x},n+1)\to\frac{\ker(\varphi_p)}{V_{p-n}\ker{\varphi_n}}=0.
\end{eqnarray*}
It suffices to show that $f_{p,n}$ is injective for $1\leq n\leq
p-1$. This will be done by induction on $n$, $1\leq n\leq p-1$. By
Proposition~\ref{fnps}, $f_{p,1}$ is injective. Suppose now that
$p\geq 3$, fix $n$ with $2\leq n\leq p-1$ and assume that $f_{p,q}$
are injective for $1\leq q\leq n-1$. In particular,
\begin{eqnarray*}
O_2(I,\underline{x},q+1)=0\, \mbox{ for } 1\leq q\leq n-1.
\end{eqnarray*}
Since $\alpha_{J,n}:\symmetric_{n}(J)\to J^n$ is an isomorphism
for each $2\leq n\leq p$, by Remark~\ref{effectivealpha}, $E(J)_t=0$
for $2\leq t\leq p$. Since by hypothesis ${\rm ker}(\alpha_{I,p})=0$,
Proposition~\ref{surjective} implies that
\begin{eqnarray*}
\frac{(I_{i}I^{t-1}:x^{t}_{i})}{(I_{i}I^{t-2}:x_{i}^{t-1})}=0\,
\mbox{ for } 1\leq i\leq s\, \mbox{ and for }2\leq t\leq p.
\end{eqnarray*}
Using Proposition~\ref{cambridge}, we get $E(I)_2=0,\ldots ,E(I)_n=0$. 
In particular, by Remark~\ref{effectivealpha}, $\alpha_{I,n}$ is an
isomorphism and $\ker(\varphi_n)=V_{n-1}\ker(\varphi_1)$.

Let us prove that $f_{p,n}$ is injective. Take
\begin{eqnarray*}
a=\sum_{|\sigma|=n}{\lsigma\xsigma}\in \biggl(\Bigl(\sum_{1\leq
i<j\leq s-1}{x_ix_jI_s^{n-1}}\Bigr) :x_s\biggr)\cap I^{n}.
\end{eqnarray*}
Suppose that $f_{p,n}$ maps the class of $a$ to the zero class. We
will show that $a$ can be written as $a=\sum_{i=1}^{s-1}{x_ia_i}$
where $a_i\in (I_{i,s}I_s^{n-1}:x_s)\cap I^{n-1}$ for $1\leq i\leq
s-1$.

By hypothesis, $x_{s}a$ can be written as
\begin{eqnarray}\label{xra}
x_s a=\sum_{1\leq i<j\leq s-1}{x_ix_ja_{ij}}\, ,\, \mbox{ with
}a_{ij}\in I_s^{n-1}.
\end{eqnarray}
Let $A_{ij}({\bf X}')$ be a form of degree $n-1$ in the variables
${\bf X}'=\{X_1,\ldots,X_{s-1}\}$ representing $a_{ij}$. As in the
proof of Proposition~\ref{fnps}, set $W=R[X_1,\ldots,X_{s-1}]$ and
$\psi:W\to\rees(J)$ sending $X_i$ to $x_it$. Since
$f_{p,s}$ maps the class of $a$ to the zero class in
$\ker(\varphi_p)/V_{p-n}\ker(\varphi_n)$ and
$\ker(\varphi_n)=V_{n-1}\ker(\varphi_1)$, we can write
\begin{multline}\label{lambdasus}
X_s^{p-n-1}\Bigl(X_s\sum_{|\sigma|=n}{\lsigma\Xsigma}-\sum_{1\leq
i<j\leq s-1}{X_iX_jA_{ij}({\bf X}')}\Bigr)= \\
\sum_{|\omega|=p-1}{(u_{\omega,1}X_1+\ldots+u_{\omega,s}X_s)X^{\omega}}
\end{multline}
where $u_{\omega,1}x_1+\ldots+u_{\omega,s}x_s=0$, $u_{\omega,i}\in R$ and
$\omega=(\omega_{1},\ldots,\omega_{s})\in\mathbb{N}^{s}$. The
following convention will be useful: if $\omega\in\mathbb{N}^{s}$ with
$|\omega|=p-1$ is of the form $\omega=\sigma+(p-n-1)e_{s}$, for a
certain $\sigma\in\mathbb{N}^{s}$ with $|\sigma|=n$, we will denote
$u_{\sigma+(p-n-1)e_s,i}$ by $v_{\sigma,i}$, for $1\leq i\leq s$. In
particular,
\begin{eqnarray}\label{vlin}
v_{\sigma,1}x_1+\ldots+v_{\sigma,s}x_s=0\; \mbox{ for }\;
\sigma\in\mathbb{N}^{s}, \;  |\sigma|=n.
\end{eqnarray}
Comparing the two polynomial expressions in equation
(\ref{lambdasus}), we observe that the monomial term $\lsigma
X^{\sigma+(p-n)e_{s}}$ of the left hand part is equal to the sum of
the following monomial terms of the right hand part:
$u_{\omega,i}X^{\omega+e_{i}}$, for $1\leq i\leq s$, where
$\sigma+(p-n)e_{s}=\omega+e_{i}$. If $1\leq i\leq s-1$ and
$i\notin{\rm supp}(\sigma)$, there is no $\omega$ such that
$\sigma+(p-n)e_{s}=\omega+e_{i}$. If $1\leq i\leq s-1$ but $i\in{\rm
supp}(\sigma)$ then $\omega=\sigma -e_{i}+e_{s}+(p-n-1)e_{s}$. If
$i=s$, $\omega=\sigma+(p-n-1)e_{s}$. Therefore,
\begin{eqnarray*}
\lsigma = \sum_{1\leq i\leq s-1\atop i\in {\rm
supp}(\sigma)}v_{\sigma-e_i+e_s,\, i}+ v_{\sigma ,\, s}\,.
\end{eqnarray*}
Then,
\begin{align*}
a&=\sum_{|\sigma|=n}\lsigma
x^{\sigma}=\sum_{|\sigma|=n}\biggl(\sum_{1\leq i\leq s-1\atop i\in
{\rm supp}(\sigma)}{v_{\sigma-e_i+e_s,\, i}}\biggr)\xsigma +
\sum_{|\sigma|=n}v_{\sigma ,\, s}\xsigma =\\
&=\sum_{i=1}^{s-1}\sum_{|\sigma|=n\atop \sigma_{i}\neq
0}v_{\sigma-e_i+e_s,\, i}\, \xsigma + \sum_{|\sigma|=n}v_{\sigma ,\,
s}\, \xsigma .
\end{align*}
For a fixed $1\leq i\leq s-1$, if $\sigma\in\mathbb{N}^{s}$ is such
that $|\sigma|=n$ and $\sigma_{i}\neq 0$, then $\theta$ defined as
$\theta=\sigma-e_{i}+e_{s}$ has $|\theta|=n$, $\theta_{s}\neq 0$ and
$\sigma=\theta+e_{i}-e_{s}$. Conversely, if $\theta\in\mathbb{N}^{s}$
is such that $|\theta|=n$ and $\theta_{s}\neq 0$, then $\sigma$
defined as $\sigma=\theta+e_{i}-e_{s}$ has $|\sigma|=n$,
$\sigma_{i}\neq 0$ and $\theta=\sigma-e_{i}+e_{s}$. Therefore,
\begin{align*}
a&=\sum_{i=1}^{s-1}\sum_{|\theta|=n\atop \theta_{s}\neq 0}v_{\theta,\,
i} x^{\theta+e_i-e_s} + \sum _{|\theta|=n}v_{\theta ,\,
s}x^{\theta}=\\ &= \sum_{|\theta|=n \atop \theta_s\neq
0}\sum_{i=1}^{s-1}v_{\theta,\, i} x^{\theta+e_i-e_s} +
\sum_{|\theta|=n \atop \theta_s\neq 0}{v_{\theta,s}x^{\theta}}+
\sum_{|\theta|=n \atop \theta_s=0}v_{\theta,s}x^{\theta}= \\
&=\sum_{|\theta|=n \atop \theta_s\neq 0}\Bigl(
\sum_{i=1}^{s-1}v_{\theta,\, i} x^{\theta+e_i-e_s} + v_{\theta, \,
s}x^{\theta}\Bigr) + \sum_{|\theta|=n \atop
\theta_s=0}v_{\theta,s}x^{\theta}.
\end{align*}
Using equation (\ref{vlin}), for a $\theta\in\mathbb{N}^{s}$,
$|\theta|=n$, then $v_{\theta,1}x_1+\ldots+v_{\theta,s}x_s=0$. If
$\theta_{s}\neq 0$, then $v_{\theta,s}x^{\theta}=
v_{\theta,s}x_{s}x^{\theta-e_{s}}=-\sum_{i=1}^{s-1} v_{\theta
,i}x_{i}x^{\theta -e_{s}}=-\sum_{i=1}^{s-1} v_{\theta ,i}x^{\theta
+e_{i}-e_{s}}$. Therefore,
\begin{eqnarray}\label{asenzill}
a & = & \sum_{|\theta|=n \atop \theta_s=0}{v_{\theta,s}x^{\theta}}.
\end{eqnarray}
From now on let $\mathcal{U}=\{\theta\in\mathbb{N}^s \mid |\theta|=n\;
\mbox{ and } \theta_s=0\}$. Then, by equation (\ref{vlin}) again,
\begin{eqnarray*}
x_s a & = & \sum_{\theta\in\mathcal{U}}{v_{\theta,s}x_sx^{\theta}}
=\sum_{\theta\in\mathcal{U}}
{(-v_{\theta,1}x_1-\ldots-v_{\theta,s-1}x_{s-1})x^{\theta}}.
\end{eqnarray*}
Using equation (\ref{xra}) and the subsequent definitions, it is clear
that the form
\begin{equation}\label{gdex1}
G({\bf X}') = \sum_{\theta\in\mathcal{U}}
{(v_{\theta,1}X_1+\ldots+v_{\theta,s-1}X_{s-1})X^{\theta}}
+\sum_{1\leq i<j\leq s-1}{X_iX_jA_{ij}({\bf X}')},
\end{equation}
belongs to $\ker(\psi_{n+1})$. Since
$\alpha_{J,2},\ldots,\alpha_{J,p}$ are isomorphisms, then
$\ker(\psi_{n+1})=V_n\ker(\psi_1)$ and $G({\bf X}')$ can be written as
\begin{eqnarray}\label{gdex2}
G({\bf X}')=\sum_{\theta\in\mathcal{U}}{(w_{\theta,1}X_1+\ldots
+w_{\theta,s-1}X_{s-1})X^{\theta}}
\end{eqnarray}
with $w_{\theta,1}x_1+\ldots+w_{\theta,s-1}x_{s-1}=0$, $w_{\theta,j}\in
R$. Comparing the two polynomial expressions in (\ref{gdex1}) and
(\ref{gdex2}), we observe that the term $v_{ne_{i},\,i}X^{(n+1)e_{i}}$
is equal to $w_{ne_i,\,i}X^{(n+1)e_{i}}$, for $1\leq i\leq s-1$. In
particular, $v_{ne_i,\,i}=w_{ne_i,\,i}\in(I_{i,s}:x_i)$ for each
$1\leq i\leq s-1$.

Let $A({\bf X}')=\sum_{\theta\in\mathcal{U}}
{v_{\theta,s}X^{\theta}}\in W_{n}$, which by equation (\ref{asenzill})
represents $a$.

Set $\mathcal{U}_{1}=\{ \theta\in\mathcal{U}\mid \theta_{1}\neq 0\}$
and $\mathcal{U}_i=\{\theta\in\mathcal{U} \mid \theta_j=0, \; 1\leq
j\leq i-1 \; \mbox{and} \; \theta_i\neq 0\}$ for each $2\leq i\leq
s-1$. Clearly, the $\mathcal{U}_{i}$ define a partition of
$\mathcal{U}$.

Let $A_i({\bf X}')=\sum_{\theta\in\mathcal{U}_i}
{v_{\theta,s}X^{\theta-e_i}}\in W_{n-1}$. Then we have $A({\bf X}')=
X_1A_1({\bf X}')+\ldots+X_{s-1}A_{s-1}({\bf
X}')$. Evaluating at $x_{1},\ldots, x_{s-1}$, one obtains the equality
$a = A(x_{1},\ldots ,x_{s-1})=\sum_{i=1}^{s-1}{x_ia_i}$, where
$a_i=A_i(x_i,\ldots,x_{s-1})=
\sum_{\theta\in\mathcal{U}_i}v_{\theta,s}x^{\theta-e_i}\in
I^{n-1}$. We claim that $a_i\in (I_{i,s}I_s^{n-1}: x_s)\cap I^{n-1}$
for each $1\leq i\leq s-1$. By equation (\ref{vlin}), 
\begin{eqnarray*}
x_sa_i=\sum_{\theta\in \mathcal{U}_i}{
v_{\theta,s}x_sx^{\theta-e_i}}=\sum_{\theta\in
\mathcal{U}_i}{(-v_{\theta,1}x_1-\ldots-v_{\theta,s-1}x_{s-1})x^{\theta-e_i}}.
\end{eqnarray*}
Fix $\theta\in\mathcal{U}_i$. If ${\rm supp}(\theta)$ contains $j\neq
i$, then
$v_{\theta,i}x_ix^{\theta-e_i}=v_{\theta,i}x_ix_jx^{\theta-e_j-e_i}=
v_{\theta,i}x_jx^{\theta-e_j}\in I_{i,s}I_s^{n-1}$; if ${\rm
supp}(\theta)=\{i\}$, then $\theta=ne_i$ and
$v_{\theta,i}x_ix^{\theta-e_i}=v_{ne_i,i}x_i^{n}\in
(I_{i,s}:x_i)x_i^n\subset I_{i,s}I_s^{n-1}$. So
$a=\sum_{i=1}^{s-1}{x_ia_i}$, where $a_i\in (I_{i,s}I_s^{n-1}:x_s)\cap
I^{n-1}$, proving the injectivity of $f_{p,n}$ as desired.\qed

\section{Theorem B}\label{mainresults}

We have now all the ingredients to prove the main result of the chapter.

\begin{Theorem B}\label{maintheorem}
Let $R$ be a commutative ring, let $I=(x_1,\ldots ,x_s)$ be an ideal
of $R$ and let $p\geq 2$ be an integer. Suppose that $J=(x_1,\ldots ,x_{s-1})$ 
verifies that $\alpha_{J,n}:\symmetric_{n}(J)\to J^n$ is an isomorphism for 
each $2\leq n\leq p$.  Then the following conditions are equivalent:
\begin{enumerate}
\item[$(i)$] $\alpha_{I,p}:\symmetric_{p}(I)\to I^p$ is an isomorphism;
\item[$(ii)$] $\alpha_{I,n}:\symmetric_{n}(I)\to I^n$ is an isomorphism for
each $2\leq n\leq p$.
\end{enumerate}
\end{Theorem B}

\demo Suppose $(i)$ and assume $s\geq 3$ (add repeated generators if
necessary). By Proposition~\ref{surjective},
$(I_iI^{n-1}:x_i^n)/(I_iI^{n-2}:x_i^{n-1})=0$ for each $1\leq
i\leq s$ and $2\leq n\leq p$. By Proposition~\ref{c2zero},
\begin{equation*}
O_2(I,\underline{x},n+1)=0 \textrm{ for each } 1\leq n\leq p-1.
\end{equation*}
By Remark \ref{effectivealpha}, $E(J)_{n}=0$ for every $2\leq n\leq
p$.  Using Proposition~\ref{cambridge}, we have
$E(I)_2=0,\ldots,E(I)_{p}=0$ and $(ii)$ follows by Remark
\ref{effectivealpha}.  \qed

\begin{Corollary}\label{Jx}
Let $R$ be a commutative ring, let $J\subset R$ be a finitely
generated ideal of linear type and let $I=(J,y)$, i.e., $I$ is 
an ideal of almost-linear type. Then the following conditions 
are equivalent:
\begin{enumerate}
\item[$(i)$] $\alpha_{I,p}:\symmetric_{p}(I)\to I^p$ is an isomorphism;
\item[$(ii)$] $\alpha_{I,n}:\symmetric_{n}(I)\to I^n$ is an isomorphism for
each $2\leq n\leq p$.
\end{enumerate}
In particular, this equivalence holds if either
\begin{enumerate}
\item[1)] $R$ is reduced and $I$ is two-generated, or
\item[2)] $R$ is an integrally closed domain and $I$ is three-generated.
\end{enumerate}
\end{Corollary}
\demo If $I=(x,y)$ is a two generated ideal in a reduced ring, $J=(x)$ 
is of linear type. If $I=(x,y,z)$ is a three-generated ideal in an 
integrally closed domain, the ideal $J=(x,y)$ is of linear type by 
\cite[Theorem~3]{costa1}.  \qed

\begin{Corollary}
Let $\{x_1,\ldots,x_s\}$ be a sequence of elements of a commutative
ring $R$ and let $J_i=(x_1,\ldots,x_i)$. Let $p\geq 2$ be an integer. 
Suppose that $\alpha_p:\symmetric_{p}(J_i)\to J_i^p$ is an isomorphism for each
$1\leq i\leq s$.  Then $\alpha_n:\symmetric_{n}(J_i)\to J_i^n$ is an
isomorphism for each $1\leq i\leq s$ and for each $2\leq n\leq p$.
\end{Corollary}
\demo Proceed by induction on $s\geq 1$. Since $\alpha_{J_1,p}$ is
an isomorphism, Example~\ref{casprincipal} implies that
$\alpha_{J_1,n}$ is an isomorphism for each $n\geq 2$. Let $s\geq 2$
and suppose that $\alpha_{J_t,n}$ is an isomorphism for
each $2\leq n\leq p$ and each $2\leq t\leq s-1$. Since by hypothesis
$\alpha_{J_s,n}$ is an isomorphism,
Theorem~B implies that $\alpha_{J_s,n}$
is an isomorphism for each $2\leq n\leq p$.  \qed

\begin{Remark}\label{summary-known-cases}{\rm
To summarize, let us gather the results of \cite{tchernev} with those 
presented above. These are the situations in which the 
injectivity of the graded components of the canonical morphism $\alpha$ 
propagates downwards: $(i)$ $I$ such that $\symmetric(I)$ contains a 
regular element in degree one; $(ii)$ $I=\fm$ maximal 
ideal (see \cite[Theorem~5.5]{tchernev}); $(iii)$ $I$ such that ${\rm pd}(I)\leq 1$ 
(see \cite[Theorem~5.1]{tchernev}); 
$(iv)$ $I$ perfect Gorenstein of grade $3$ (see \cite[Theorem~5.3]{tchernev}); 
$(v)$ $I=(J,y)$ with $J$ is of linear type, i.e., $I$ almost-linear type 
(see Corollary~\ref{Jx}), in particular, either when $I=(x,y)$ and $R$ 
reduced or when $I=(x,y,z)$ and $R$ integrally closed domain.
}\end{Remark}

\begin{Remark}
\rm K\"uhl proved that in a Noetherian ring, linear type and 
projectively of linear type are equivalent conditions for an ideal 
$I$ provided that $I$ is generated by a proper sequence (see \cite[Corollary~2.4]{kuhl}). 
We have achieved new cases in which this equivalence holds: in fact, by Corollary~\ref{Jx}, 
if $I$ is an ideal of the form $I=(J,y)$ with $J$ of linear type, then $I$ 
is of linear type if and only if $I$ is projectively of linear type. We 
can state this claim together with other cases considered so far in the 
following corollary.
\end{Remark}

\begin{Corollary}\label{p-linear-type=linear-type}
Let $R$ be a Noetherian ring and $I$ an ideal of $R$. Consider any of the 
situations considered in Remark~\ref{summary-known-cases}:
\begin{enumerate}
\item[$(a)$] $\symmetric_1(I)$ contains a non-zero-divisor of $\symmetric(I)$;
\item[$(b)$] $I=(J,y)$ with $J$ an ideal of linear type;
\item[$(c)$] ${\rm pd}(I)\leq 1$;
\item[$(d)$] $I$ is perfect Gorenstein of grade $3$;
\item[$(e)$] $I=\fm$ is a maximal ideal.
\end{enumerate}
In any of the situations $(a)-(e)$, $I$ is of linear type if and only if $I$ is of 
$p$-linear type, i.e., the canonical morphism $\alpha_I:\symmetric(I)\to\rees(I)$ 
is an isomorphism if and only if the induced morphism
$\widetilde{\alpha}_I:{\rm Proj}(\rees(I))\to {\rm Proj}(\symmetric(I))$
is an isomorphism of schemes.	
\end{Corollary}
\demo
The claim is a straightforward consequence of Remark~\ref{p-linear-type-giral} and 
Remark~\ref{summary-known-cases} above.
\qed

\section{Examples and applications}\label{sec:exemple-article1}

The following example shows that if the hypothesis of $\alpha_{J,n}$
being an isomorphism for each $2\leq n\leq p$ is not fulfilled, then
the conclusion in Theorem~B may fail. Concretely, for
each $p\geq 2$, we construct a two-generated ideal $I=(x,y)$, with
neither $(x)$ nor $(y)$ of linear type, such that $\alpha_{I,n}$ is an
isomorphism for each $n\geq p+1$, whereas $\alpha_{I,p}$ is not (see
also \cite[Example~1.4]{kuhl}). Experimentation with
Singular~\cite{singular} was useful in the process of
generating this as well as other examples.

\begin{Example}\label{exemple}
Let $k$ be a field and let $p\geq 2$. Let $S=k[U_0,\ldots,U_p,X,Y]$ be
a polynomial ring and let $Q$ be the ideal of $S$ defined as $Q=Q_1 +
(U_0X^p)$, where
\begin{eqnarray*}
Q_1=(U_0Y,U_0X-U_1Y,U_1X-U_2Y,\ldots, U_{p-1}X-U_pY, U_pX).
\end{eqnarray*}
Let $R$ be the factor ring $S/Q=k[u_0,\ldots,u_p,x,y]$ and consider
the ideal $I=(x,y)\subset R$. Then $\alpha_{I,n}$ is an isomorphism
for all $n\geq p+1$, whereas $\alpha_{I,p}$ is not. Moreover,
$\beta_{I,n}$ is an isomorphism for all $n\neq p$, whereas
$\beta_{I,p}$ is not.
\end{Example}

\demo Let $A=k[U_0,\ldots,U_p]$ and $S=A[X,Y]$. We can endow $S$ with
a grading by letting ${\rm deg}(X)=1$, ${\rm deg}(Y)=1$ and ${\rm
deg}(w)=0$ for each $w\in A$. Then $R=S/Q=A[x,y]$ has a standard
graded $A$-algebra structure, with irrelevant ideal given by
$I=R_{+}=(x,y)$.  Note that $\graded(I)=\oplus_{p\geq 0}I^p/I^{p+1}\cong
R$ and $\symmetric(R/I){I/I^2}\cong\symmetric_A{R_1}$. Consider the free
presentation of $R_{1}$
\begin{eqnarray*}
0\to\widetilde{Q}_1\to Ae_x\oplus Ae_y\stackrel{f}{\to} R_1\to 0,
\end{eqnarray*}
where $f$ sends $e_x$ to $x$, $e_y$ to $y$ and $\widetilde{Q}_1$ is
the submodule of $Ae_x\oplus Ae_y$ generated by
$U_0e_y,U_0e_x-U_1e_y,...,U_{p-1}e_x-U_pe_y,U_pe_x$. Applying the
symmetric functor, it gives rise to the top row in the following
commutative diagram of exact rows:
\begin{eqnarray*}
\xymatrix{0 \ar[r] & Q_1 \ar[d] \ar[r] & A[X,Y] \ar@{=}[d] \ar[r] &
          \symmetric_{A}(R_1) \ar@{->>}[d]^{\;\; \beta_I} \ar[r] & 0 \; \\
          0 \ar[r] & Q \ar[r] & A[X,Y] \ar[r] & R \ar[r] & 0. }
\end{eqnarray*}
Therefore, $\ker\beta_I\cong Q/Q_{1}$, which is generated by
the class of $U_0X^p$ in $A[X,Y]/Q_1$. 

Let us see that $U_0X^p\notin Q_1$. Indeed, if $U_0X^p\in Q_1$ then
$U_0X^p=F_0U_0Y+F_1(U_0X-U_1Y)+\ldots+F_p(U_{p-1}X-U_pY)+F_{p+1}U_pX$,
for some $F_{i}\in S=k[U_0,\ldots,U_p,X,Y]$. We can assume that the
$F_i$ are homogeneous polynomials in $X,Y$ of degree $p-1$. Since
$(F_0Y-X^p+F_1X)U_0+(F_2X-F_1Y)U_1+\ldots+(F_{p+1}X-F_pY)U_p=0$, we
get $F_{i+1}X=F_{i}Y$ for each $1\leq i\leq p$. In particular
$F_{1}Y^{p}=F_{2}XY^{p-1}=\ldots =F_{p+1}X^{p}$ and hence $F_1\in
X^pS$. But since ${\rm deg}(F_1)=p-1$, $F_{1}$ must be zero, and hence
$F_{2}=\ldots=F_{p+1}=0$. Thus $F_{0}Y=X^{p}$, a contradiction.

Then $\ker(\beta_{I,p})\neq 0$, hence $\alpha_{I,p}$ is not an
isomorphism. Clearly $\ker(\beta_{I,n})=0$ for each $2\leq n\leq
p-1$. On the other hand, from the polynomial equalities given by
\begin{eqnarray*}
U_0X^{p+1} &=& X^p(U_0X-U_1Y)+ X^{p-1}Y(U_1X-U_2Y)+\ldots+\\ &&
XY^{p-1}(U_{p-1}X-U_pY)+Y^p(U_pX)\in Q_1,\\ U_0X^pY & = & X^p(U_0Y)\in
Q_1,
\end{eqnarray*}
we can conclude that $\ker(\beta_{I,n})=0$ also for $n\geq p+1$. Hence
$\ker(\beta_{I,n})=0$ for each $n\neq p$.

Consider the downgrading homomorphism
$\lambda:\symmetric_{n+1}(I)\to\symmetric_{n}(I)$, which maps the homogeneous
element $a_1\cdot a_2\cdot\ldots\cdot a_{n+1}$ to
$a_1(a_2\cdot\ldots\cdot a_{n+1})$. Note that in $a_1(a_2\cdot\ldots\cdot a_{n+1})$, 
$a_1$ is seen as an element in $R$ and the product on the left 
is the one given by the $R$-module structure of $\symmetric(I)$. 
It is well-known that $\ker(\beta_{I,n})\cong\ker(\alpha_{I,n})/\lambda(\ker(\alpha_{I,n+1}))$
(see \cite[Section~3]{hsv}). In our case, it implies that
$\ker(\alpha_{I,n})=\lambda^t(\ker(\alpha_{I,n+t}))$ for all $n\geq
p+1$ and $t\geq 1$.

Let us see now that $\ker(\alpha_{I,n})=0$ for all $n\geq p+1$. 
Since $R$ is noetherian, ${\rm ker}(\alpha_{I})$ is a finitely
generated ideal of $\symmetric(I)$. Therefore, there exists an integer
$r\geq 1$ such that ${\rm ker}(\alpha_{I,r+t})=\symmetric_{1}(I)^{t}\cdot
{\rm ker}(\alpha_{I,r})$ for all $t\geq 1$ (recall that the smallest of such
integers $r\geq 1$ is $\reltype(I)$, the relation type of $I$; see
e.g. \cite{planas} or \cite{vasconcelos}). If $n\geq {\max}\{p+1,r\}$,
then
\begin{equation*}
{\rm ker}(\alpha_{I,n})=\lambda ^{n}({\rm
ker}(\alpha_{I,2n}))=\lambda^{n}(\symmetric_{1}(I)^{n}\cdot {\rm
ker}(\alpha_{I,n}))=I^{n}{\rm ker}(\alpha_{I,n}).
\end{equation*}
But $I^{n}{\rm ker}(\alpha_{I,n})=0$ since if $w\in {\rm
ker}(\alpha_{I,n})$, then $I^nw\subset \alpha(w)\symmetric_{n}(I)=0$.  On
the other hand, if $p+1\leq n<r$, observe that
$\ker(\alpha_{I,n})=\lambda^r(\ker(\alpha_{I,n+r}))$ and $n+r\geq {\rm
max}\{p+1,r\}$.  Hence $\ker(\alpha_{I,n})=0$ by the previous case.
\qed

\begin{Remark}{\rm
Alternatively, one can use the Artin-Rees Lemma to finish the proof of
the former example. Indeed, since $\ker(\beta_{I,n})=0$ for each
$n\geq p+1$, we have
$\ker(\alpha_{I,n})=\lambda^s(\ker(\alpha_{I,n+s}))\subset
I^s\symmetric_{n}(I)\cap\ker(\alpha_{I,n})$ for each $s\geq 1$ and each
$n\geq p+1$. Let $s_n$ be the Artin-Rees number relative to the pair
$\ker(\alpha_n)\subset\symmetric_{n}(I)$ and the ideal $I$. Then the
Artin-Rees Lemma asserts that $I^s\symmetric_{n}(I)\cap\ker(\alpha_{I,n})=
I^{s-s_n}(I^{s_n}\symmetric_{n}(I)\cap\ker(\alpha_{I,n}))\subset
I^{s-s_n}\ker(\alpha_{I,n})$ for each $s\geq s_n$. Taking $s$ large
enough, we see that $\ker(\alpha_{I,n})\subset
I^n\ker{\alpha_{I,n}}$, which is zero by the argument above (see also
\cite[Corollary~1.3]{kuhl}). Therefore $\ker(\alpha_{I,n})=0$ for each
$n\geq p+1$.  }\end{Remark}

\begin{Remark}\label{many-non-zeros}{\rm
Let $n_{1},\ldots ,n_{r}\geq 2$ be integers. Then one can construct a
Noetherian ring $A$ and an ideal $I$ of $A$ such that $\beta_{I,n}$ is
an isomorphism if and only if $n\neq n_{1},\ldots ,n_{r}$. This
follows using a construction similar to Example~\ref{exemple},  
taking as many sets of variables, with the corresponding relations, 
as there are integers $n_{1},\ldots,n_{r}$. See Example \ref{tensored-example} 
for the details.
}\end{Remark}

\begin{Remark}\rm
Remark that there are ideals $I$ verifying the hypotheses of Theorem~B 
and such that $\symmetric(I)$ has no regular elements 
in degree one. Take, e.g., $R=k[X,Y,Z]/(XZ,YZ)=k[x,y,z]$ and $I=(x,y)\subset R$.
\end{Remark}

\begin{Remark}
\rm Graph ideals provide us with valuable and manageable examples. 
Let $R=k[X_{1},\ldots,X_{2p}]$ be the polynomial ring with 
variables $X_{1},\ldots,X_{2p}$ over a field $k$. Let $J$ be the 
graph ideal $(X_{1}X_{2},X_{2}X_{3},\ldots ,X_{2p-1}X_{2p})$
corresponding to a path of length $2p-1$. It is well-known that
$J$ is of linear type, i.e., $\alpha_{J,n}:\symmetric_{n}(J)\rightarrow J^{n}$ 
is an isomorphism for all $n\geq 2$ (see \cite[Proposition~3.1]{villarreal}). 
Set $I=(J,X_{2p}X_{1})$, the graph ideal corresponding to the closed walk of 
length $2p$. It can be shown that $\alpha_{I,p}:\symmetric_{p}(I)\rightarrow I^{p}$ 
is not an isomorphism (see \cite[Proposition~3.1]{villarreal}). Then
$\alpha_{I,n}$ is not an isomorphism for any $n\geq p$, since otherwise 
$\alpha_{I,p}$ would be an isomorphism by Theorem~B, 
leading to a contradiction.
\end{Remark}
\chapter{Further examples and future work}\label{oddsandends}

We present a set of tools behind many constructions and key ideas used 
in the previous chapters. Many results, examples and remarks are 
outcomes of these techniques. 

The equations of the Rees algebra of an ideal modulo a regular element have been
an object of study since the influential works of \cite[Proposition 3.3]{hsv} 
and \cite{valla} on ideals of linear type. Here we tackle the 
question whether the equations of the Rees algebra of a given ideal 
remain invariant when taking the ideal modulo a regular element. 
We prove Proposition~\ref{bounds-injectivity}, a result that generalises
\cite[Proposition~3.3]{hsv} and \cite[Theorem~2.1]{valla}, relying on the module 
of effective relations, and we are able to deduce at once yet another proof 
of a classical result recorded in the claim~\ref{d-sequences-are-linear-type}, 
namely, that an ideal generated by a $d$-sequence is of linear type.

When studying the equations of the Rees algebras of ideals
in terms of its generating sets and reductions, it is very 
useful to rely on examples where the generators, basic ideal 
operations, computation of reduction and equations of 
blowing-up algebras are intuitive enough to allow simple and
flexible examples and to avoid burdensome experimentation. 

Efforts have been devoted to construct families of ideals 
with prescribed equations. We have explored the obstructions 
arising in the construction of examples of ideals with a 
predefined pattern of equations. We have proved a method for
constructing ideals with predefined top degree equations; 
as well, we provide a general version of Example~\ref{exemple}. 

Finally, a tentative procedure for expanding sets of equations 
of Rees algebras is presented. It is a variant of a well-known determinantal 
procedure leading to the notion of \emph{expected equations}. While the procedure 
presented has been used by the author as a way to expand sets of non-trivial 
equations of Rees algebras, it is by no means well understood why it works and 
why it should work under more general settings: it remains an appealing object of 
future study. Surprisingly enough, in the examples examined, if we try to expand 
a generating set of the first syzygies, the procedure leads, with few adjustements, 
to full generating sets of equations.

\section{The equations of Rees algebras of ideals modulo a regular element}

The interplay between the equations of the Rees algebra of an
ideal $I$ of $R$, the Koszul homology of generating sets of $I$ and 
conditions on these generating sets, yield results that allow for 
setting bounds on the relation type of $I$.

First of all, let us recall a result from \cite{planas-crelle}:

\begin{Lemma}[see \cite{planas-crelle}]\label{bounds-exact-sequence}
Let $R$ be a Noetherian ring, $I$ an ideal of $R$ and $y\in I$. Then, 
for each integer $p\geq 2$, there is an exact sequence of $R$-modules
$$E(I)_p\stackrel{\sigma_p}{\to} E(I/(y))_p\to \frac{(y)\cap I^p}{I((y)\cap I^{p-1})}\to 0.$$
\end{Lemma}
\demo
Use \cite{planas-crelle} with $N=(y)\subseteq M=R$.
\qed

\medskip

The next result generalises \cite[Proposition 3.3]{hsv} and \cite[Theorem 2.1]{valla}. 
Note that this sort of statement is already tacked in \cite[Theorem 2.3]{huckaba2}.

\begin{Proposition}\label{bounds-injectivity}
Let $R$ be a Noetherian ring and let $I$ be an ideal of $R$. 
Let $y\in I\backslash I^2$ such that $(0:y)\cap I=0$ and such that the 
Valabrega-Valla module $$VV_{(y)}(I)_{p-1}=(y)\cap I^{p-1}/yI^{p-2}$$ 
is equal to zero, for some $p\geq 2$. Then, there is an exact sequence 
of $R$-modules $$0 \to E(I)_p \stackrel{\sigma_p}{\to} E(I/(y))_p \to VV_{(y)}(I)_p \to 0.$$
In particular, if $y^*$ is $\graded(I)$-regular, then $E(I)_p\cong E(I/(y))_p$ 
for all $p\geq 2$ and $\reltype(I)=\reltype(I/(y))$.
\end{Proposition}
\demo
Let $I=(x_1,\ldots,x_s)$ and let $R'=R/(y)$, $W=R[T_1,\ldots,T_s]$ and $W'=R'[T_1,\ldots,T_s]$. 
Let $\varphi$ be the polynomial presentation from $W$ to $\rees(I)$ sending $T_i$ to $x_it$ and let 
$\psi$ be the induced polynomial presentation from $W'$ to $\rees(I/(y))$, respectively, and let 
$\pi_S$ and $\pi_R$ the natural projections, as in the following diagram,
\begin{equation*}
\begin{CD}
0 @>>> Q @>>> \;\; W @>\varphi>> \rees(I) @>>> 0\\
&& @VVV  {\pi_S}@VVV {\pi_R}@VVV \\
0 @>>> Q' @>>> \;\; W' @>\psi>> \rees(I/(y)) @>>> 0.
\end{CD}
\end{equation*}
It is readily seen that, for each $p\geq 2$, $\pi_S$ induces an 
homomorphism $$\widetilde{\pi}_p:(Q/Q\langle p-1\rangle)_p \to (Q'/Q'\langle p-1\rangle)_p$$ 
such that the following is a commutative diagram:
\begin{equation*}
\begin{CD}
(Q/Q\langle p-1\rangle)_p @>{\widetilde{\pi}_p}>> \;\; (Q'/Q'\langle p-1\rangle)_p \\
@V{\cong}VV @V{\cong}VV \\
E(I)_p @>{\sigma_p}>> \;\; E(I/(y))_p 
\end{CD}
\end{equation*}
In fact, $\widetilde{\pi}_p$ sends the class of a given form $F\in Q_p$ 
to the class $\overline{F}$ of $F$ modulo $(y)$ in $(Q'/Q'\langle p-1\rangle)_p$. 
If $F\in Q\langle p-1\rangle \cap W_p$, then $F=\sum_{i=1}^s T_i F_i$, 
with $F_i\in Q\langle p-1\rangle \cap W_{p-1}$ and it follows $F'\in Q'\langle p-1\rangle\cap W'_p$, 
hence $\widetilde{\pi}_p$ is well-defined.

\medskip

We will show that $\widetilde{\pi}_p$ is injective, hence we will conclude that 
$\sigma_p$ is injective.

\medskip

Take $F\in Q_p$ a representative of an element in the kernel 
of $\sigma_p$. Then $\overline{F}\in W'$ can be written as 
$\overline{F}=\sum_{i=1}^s T_i \overline{G}_i$, where
$\overline{G}_i\in Q_{p-1}'$, thus $F=\sum_{i=1}^s T_iG_i + yK$ with $K\in W_p$, 
$G_i\in W_{p-1}$, verifying 
\begin{eqnarray}\label{condition-gi}
G_i(x_1,\ldots,x_s)\in (y)\cap (x_1,\ldots,x_s)^{p-1}=(y)\cap I^{p-1}=yI^{p-2}.
\end{eqnarray}
Observe that we can assume without loss of generality that $F=\sum_{i=1}^s T_iG_i$ 
with the $G_i$ satisfying the last condition \ref{condition-gi}, since $K$ can be 
written as $\sum_{i=1}^s T_i K_i$, thus we get $F=\sum_{i=1}^s T_i(G_i+yK_i)$ and 
$\overline{G_i+yK_i}=\overline{G}_i$.

So let us assume that $F=\sum_{i=1}^s T_iG_i$ with the $G_i$ fulfilling condition 
(\ref{condition-gi}). Let $y=\sum_{j=1}^s \lambda_jx_j$. For each $1\leq i\leq s$, there exist 
$H_i\in W$, ${\rm deg}(H_i)=p-2$ such that $G_i-(\sum_{j=1}^s \lambda_jT_j)H_i\in Q_{p-1}$. 
Therefore we have:
\begin{align*}
F & = \sum_{i=1}^s T_i(G_i-(\sum_{j=1}^s \lambda_jT_j)H_i) + \sum_{i=1}^s (\sum_{j=1}^s \lambda_jT_j)T_iH_i \\
  & = \sum_{i=1}^s T_i(G_i-(\sum_{j=1}^s \lambda_jT_j)H_i) + (\sum_{j=1}^s \lambda_jT_j)\sum_{i=1}^s T_iH_i.
\end{align*}

Consequently, $y(\sum_{i=1}^s x_iH_i(x_1,\ldots,x_s))=0$. Since $I\cap (0:y)=0$, then
$$\sum_{i=1}^s x_iH_i(x_1,\ldots,x_s)=0,$$ thus $F\in Q\langle p-1\rangle$, showing that 
the homomorphism $\sigma_p$ is injective. The remaining claims follow from a straightforward 
application of Lemma \ref{bounds-exact-sequence}.
\qed

\medskip

Contrary to what one may guess in the light of some well-known results, 
such as \cite[Proposition 3.3]{hsv} and \cite[Theorem 2.1]{valla}, the 
inequality $\reltype(I)\leq \reltype (I/(y))$ does not hold in general, 
even if $y$ is a regular element. In the following example, we observe 
that if $y\in I$ is an $R$-regular element, but $y^*$ is not 
$\graded(I)$-regular, the preceding inequality may not hold. 
\begin{Example}{\rm
Let $k$ be a field, let $R=k[t^3,t^4,t^5]$ and let $I=(t^3,t^4)$. 
Then $\reltype(I)=3 > \reltype(I/(t^3))=2$.
}\end{Example}
\demo
For convenience we will let $x=t^3$ and $y=t^4$. The effective 
relations of the ideals $I\subset R$ and $I/(x)\subset R/(x)$
can be computed using either \ref{succexcurta-RL(I)} or Lemmas 
\ref{lemma:effective-principal} and \ref{exactsequence}. In fact, 
since $R$ is a domain we have $E(I)_n\cong (xI^{n-1}:_R y^n)/(xI^{n-2}:_R y^{n-1})$ 
and $E(I/(x))\cong (x:_ R y^n)/(x:_R y^{n-1})$. Straightforward computations lead to 
the result. Observe that $x^*$ is not a regular element in $\graded(I)$, since
$(I^2:_R x)=(t^3,t^4,t^5)$ but $t^5\notin I$.
\qed

\medskip

As a consequence of the previous results and considering the elementary 
properties of $d$-sequences (see \cite{huneke4} and \cite{hsv}), we are able 
to recover the well-known fact that ideals generated by $d$-sequences are of linear 
type (see \cite{huneke2}, \cite{valla} and Section \ref{d-sequences-are-linear-type}). 
Let us split the statement and proof of this result, stating first a key lemma. 

\begin{Lemma}[see \cite{hsv}]
Let $J\subset I$ be ideals of a Noetherian ring $R$. If $I/J$ is an ideal of linear
type in $R/J$, then $J\cap I^n=JI^{n-1}$, for all $n\geq 0$.
\end{Lemma}

\begin{Proposition}\label{d-sequences-linear-type}
Let $R$ be a Noetherian ring and $I$ an ideal of $R$. If $I$ is generated
by a $d$-sequence, then $I$ is of linear type, i.e., $\reltype(I)=1$.
\end{Proposition}
\demo
Let $x_1,\ldots,x_s$ be a $d$-sequence and $I$ the ideal it generates. 
Let us proceed by induction on $s$. If $s=1$ then being a $d$-sequence 
is equivalent to $(0:x_1)\cap (x_1)=0$. Then $E(I)_p\cong (0:x_1^p)/(0:x_1^{p-1})=0$,
for all $p\geq 2$ and the result follows.

Suppose that the result holds for $d$-sequences of length $s-1$.
If $x_1,\ldots,x_s$ is a $d$-sequence in $R$, then $x'_2,\ldots,x'_s$ is
also a $d$-sequence in $R'=R/(x_1)$ (see \cite[Definition 1.1 and Remarks]{huneke4}). 
By the induction hypothesis, $I/(x_1)$ is of linear type. 
From \cite[Proposition 2.1]{huneke4} we know that $(0:x_1)\cap I = 0$. Using the
previous lemma it follows that $(x_1)\cap I^n=(x_1)I^{n-1}$, for all $n\geq 1$ 
(see also \cite[Theorem 2.1]{huneke4}). Using Proposition~\ref{bounds-injectivity} 
we get $E(I)_p\cong E(I/(x_1))_p$, for all $p\geq 2$. Since $I/(x_1)$ is of linear type 
for all $p\geq 2$ we have $E(I/(x_1))=0$ thus $E(I)_p=0$ for all $p\geq 2$ and the claim follows.
\qed

\section{Rees algebras with prescribed equations}

Two methods are presented to obtain ideals with a predefined set of 
equations. The first method yields two-generated irrelevant ideals of standard 
algebras with prescribed relation type and top degree equations. The 
second method stretches the arguments presented in Example~\ref{exemple} to show 
ideals with prescribed non-vanishing patterns of the graded components of 
$\ker\beta_I$.

\subsection{Prescribed relation type and top degree equations}\label{sec:prescribed-relation-type}

Let $A=k[U_1,\ldots, U_t]$ where $k$ is a field and 
$U_1,\ldots,U_s$ are indeterminates over $k$. Let 
$X_1,\ldots,X_s$ be indeterminates over $A$ and endow the 
polynomial ring $A[X_1,\ldots,X_s]$ with a grading by setting 
${\rm deg}(X_i)=1$ for all $i=1,\ldots,s$ and ${\rm deg}(w)=0$ 
for each $w\in A$. Let $L$ be an homogeneous ideal of 
$A[X_1,\ldots,X_s]$, generated by elements of positive degree. 
Set $R=A[X_1,\ldots,X_s]/L=\oplus_{n\geq 0} R_{n}=A[x_1,\ldots,x_s]$, 
where $x_i$ stands for the class of $X_i$. The ring $R$ has a 
standard graded $A$-algebra structure, with irrelevant ideal 
given by $I=R_{+}=(x_1,\ldots,x_s)$. 

Notice that $\graded(I)=\oplus_{n\geq 0}I^n/I^{n+1}$ is isomorphic
to $R$ as $A$-algebras and $\symmetric_{R/I}(I/I^2)\cong\symmetric_A(R_1)$. 
Let $\varphi:A^s\to R_1$ be a free presentation of $R_1$, sending each basis 
element $e_i$ to $X_i$. Then there is a short exact sequence
\begin{eqnarray*}
0\to\widetilde{L}_1\to Ae_1\oplus \ldots \oplus Ae_s\stackrel{f}{\to} R_1\to 0.
\end{eqnarray*}
Notice that there is an isomorphism of $A$-modules $\widetilde{L}_1\cong L_1$ 
induced by the map $A^s\to A[X_1,\ldots,X_s]_1$ sending $e_i$ to $X_i$.  
Applying the symmetric functor to the short exact sequence 
(see \cite[Algebra, Ch. 3, Section 6.2]{bourbaki}), one gets the 
following commutative diagram with exact rows:
\begin{eqnarray*}
\xymatrix{0 \ar[r] & L\langle 1\rangle \ar[d] \ar[r] & A[X_1,\ldots,X_s] \ar@{=}[d] \ar[r] &
          \symmetric_{A}(R_1) \ar@{->>}[d]^{\;\; \alpha} \ar[r] & 0 \; \\
          0 \ar[r] & L \ar[r] & A[X_1,\ldots,X_s] \ar[r]^{\;\;\; \varphi} & R \ar[r] & 0.}
\end{eqnarray*}
By the claims of the preceding paragraph and according to the notations adopted throughout, 
the canonical morphism $\alpha$ in the diagram above is no other than $\beta_I$. An immediate 
application of the snake lemma yields ${\rm ker}(\beta_I)\cong L/L\langle 1\rangle$. 

According to the claim \ref{defining-equations-associated-graded}, 
$\reltype(I)=\reltype(\graded(I))$ (set $r$ equal to this value) and $E(I)_r\cong E(\graded(I))_r$.
On the other hand, using the isomorphism $\ker\beta_I\cong L/L\langle 1\rangle$, we get 
$$E(\graded(I))_r = \ker\beta_{I,r}/\symmetric_1(I/I^2)\cdot \ker\beta_{I,r-1}\cong
(L/L\langle r-1\rangle)_r.$$

\subsection{Non-zero \texorpdfstring{$\ker\beta$}{kernel of beta} at prescribed degrees}\label{tensored-example}

In Section \ref{sec:exemple-article1} we constructed a ring $R_p$ and an ideal 
$I_p\subset R_p$ such that $\ker\beta$ is non-zero only 
at a given degree $p$. Moreover, $\ker\beta_{I,p}$ is cyclic.

As it is asserted in Remark \ref{many-non-zeros}, we will prove that 
given a finite set $S=\{p_{1},\ldots ,p_{r}\}$ of integers $\geq 2$, 
there exists a Noetherian ring $\mathcal{A}_S$ and an ideal 
$I_S\subset \mathcal{A}_S$ such that $\beta_{I_S,p}$ is an isomorphism 
if and only if $p\notin S$.

The idea consists in loading a new set of indeterminates and equations 
in these indeterminates for each new fresh equation we want to add to 
the kernel of $\beta_{I_S}$. If the set $S$ is finite, this will be nothing 
but a tensor product of algebras. On the other hand, considering $S$ 
infinite, would lead us to a non-Noetherian scenario that will not be 
tackled here. 

\begin{Example}\label{many-non-zeros-developed}{\rm
Let $U_{p,i}$ be indeterminates for integers $p\geq 2$ and $0\leq i\leq p$.
As in Example~\ref{exemple}, define $Q_p$ as the ideal of 
$A_p=k[X,Y][U_{p,0},\ldots,U_{p,p}]$ given by  
$$Q_p=(U_{p,0}Y,U_{p,0}X-U_{p,1}Y,U_{p,1}X-U_{p,2}Y,\ldots, U_{p,p-1}X-U_{p,p}Y, U_{p,p}X, U_{p,0}X^p).$$
Notice that $A_p/Q_p$ are $k[X,Y]$-modules. Given $S=\{p_{1},\ldots ,p_{r}\}$ a finite set of integers 
satisfying $2\leq p_1<p_2<\ldots<p_r$, define the ring $\mathcal{A}_S$ as
$$\mathcal{A}_S=(A_{p_1}/Q_{p_1})\otimes_{k[X,Y]} (A_{p_2}/Q_{p_2})\otimes_{k[X,Y]} \cdots \otimes_{k[X,Y]} (A_{p_r}/Q_{p_r}).$$
Then the ideal $I=(x,y)\subset \mathcal{A}_S$ verifies that 
\begin{equation*}
\ker\beta_{I,n}= 
\left\{
\begin{array}{ll}
\langle U_{n,0}T_1^n \rangle & n\in S, \\
0 & n\notin S.
\end{array} \right.
\end{equation*}
In particular, $I$ is of $p$-linear type and $\reltype(I)=p_r$, thus not of linear type.
}\end{Example}
\demo
In order to not too much overload the notation, we will prove the claim for a set 
$S=\{p,q\}$ with $2\leq p<q$; then the argument can easily be generalised 
to finite sets of integers. Let 
\begin{align*}
Q_{p} & = (U_{0}Y,U_{0}X-U_{1}Y,U_{1}X-U_{2}Y,\ldots, U_{p-1}X-U_{p}Y, U_{p}X, U_{0}X^{p})\\
Q_{q} & = (V_{0}Y,V_{0}X-V_{1}Y,V_{1}X-V_{2}Y,\ldots, V_{q-1}X-V_{q}Y, V_{q}X, V_{0}X^{q}),
\end{align*}
and let us set 
\begin{align*}
\mathcal{A}_S & =(A_{p}/Q_{p})\otimes_{k[X,Y]} (A_{q}/Q_{q})
= \frac{k[X,Y,U_0,\ldots,U_{p}]}{Q_{p}}\otimes_{k[X,Y]}\frac{k[X,Y,V_0,\ldots,V_{q}]}{Q_{q}} \\
& \cong \frac{k[X,Y,U_0,\ldots,U_{p},V_0,\ldots,V_{q}]}{(Q_{p},Q_{q})}= k[x,y,u_0,\ldots,u_{p},v_0,\ldots,v_{q}],
\end{align*}
$I=(x,y)\subset \mathcal{A}_S$. Let $\beta_I:\symmetric(I/I^2)\to \graded(I)$ be the canonical 
epimorphism onto the associated graded ring of $I$ and consider 
$\psi:W=\mathcal{A}_S[T_1,T_2]\to \graded(I)$, the polynomial presentation of $\graded(I)$ 
sending $T_1$ to $x^*$ and $T_2$ to $y^*$. We know that 
$\ker\beta_I\cong \ker\psi/W_+\ker\psi_1$. Using this correspondence, by the previous 
discussion in Section~\ref{sec:prescribed-relation-type}, we known that $\ker\beta_I$ is generated 
by the classes of $u_{0}T_1^{p}$ and $v_{0}T_1^{q}$ modulo the ideal generated by linear forms
\begin{eqnarray*}
\ker\psi\langle 1\rangle=W_+\ker\psi_1=(u_{0}T_2,u_{0}T_1-u_{1}T_2,u_{1}T_1-u_{2}T_2,\ldots, u_{p-1}T_1-u_{p}T_2, u_{p}T_1,\\
v_{0}T_2,v_{0}T_1-v_{1}T_2,v_{1}T_1-v_{2}T_2,\ldots, v_{q-1}T_1-v_{q}T_2, v_{q}T_1).
\end{eqnarray*}
Repeating the same arguments as in the proof of Example~\ref{exemple}, it is clear that
any combination $\mu u_{0}T_1^{p}+\nu v_{0}T_1^{q}$ with 
${\rm deg}_{T_1,T_2}(\mu)={\rm deg}_{T_1,T_2}(\nu)\geq 1$ belongs to $\ker\psi\langle 1\rangle$,
consequently vanishes in $\ker\beta_I$. Moreover, $\reltype(I)=q$ and from $\ker\beta_{I,n}=0$ 
for all $n\geq q+1$ one deduces $\ker\alpha_{I,n}=0$ for all $n\geq q+1$, 
thus $I$ is of $p$-linear type.
\qed

\section{$g$-linear type but not $p$-linear type}

Recall from Section \ref{geometric-conditions-linear-type} that we have the following 
chain of implications:
$$ \textrm{linear type} \Rightarrow \textrm{$p$-linear type} \Rightarrow \textrm{$g$-linear type}.$$

While the converse implications do not hold in general, we can see that the first one is valid
for principal ideals, that is, the notions of linear type and $p$-linear type coincide for
principal ideals.

\begin{Lemma}
Let $R$ be a Noetherian ring. Let $I$ be a principal ideal. If $I$ is of $p$-linear type, 
then $I$ is of linear type. Furthermore, any principal ideal $I$ with $\ker\alpha_{I,n}=0$ 
for some $n\geq 2$ is of linear type.
\end{Lemma}
\demo
Let $I=(x)$ be a principal ideal of $R$. By Example~\ref{casprincipal},
$\ker\alpha_{I,n}\cong (0:x^n)/(0:x)$. By Remark \ref{p-linear-type-giral}, 
$\ker\alpha_{I,n}=0$ for $n$ large enough. It follows that 
$(0:x^2)=(0:x)$.
\qed

\begin{Remark}{\rm
Recall that we have already seen examples of ideals of $p$-linear type with arbitrarily large 
relation type: recall the previous examples \ref{exemple} and \ref{many-non-zeros-developed}.
}\end{Remark}

In the next example it is proved that $g$-linear type does not imply $p$-linear type,
even for principal ideals.

\begin{Example}\label{principal-g-linear-but-not-p-linear}{\rm
Let $k$ be a field and $R=k[X,Y]/(X^2Y,XY^2)=k[x,y]$. Then $I=(y)$ is 
an ideal of $g$-linear type but not $p$-linear type. 
}\end{Example}
\demo We will show that $\ker\alpha$ is nilpotent whereas $\ker\alpha_n$ 
does not vanish for any $n\geq 2$. Let $\varphi:R[T]\to\rees(I)$ sending 
$T$ to $yt$, $Q=\ker\varphi$. Then $Q/Q\langle 1\rangle\cong\ker\alpha$ is generated 
by the class of $xT^2$. In fact, $\ker\alpha_n\neq 0$ for any $n\geq 2$, 
since the class of $xT^n$ is not zero in $\ker\alpha_n$. On the other hand, 
the class of $(xT^n)^2=x^2T^{2n}=(x^2T)\cdot T^{2n-1}$ in $\ker\alpha_n$ is 
the zero class, since $x^2T\in Q_1$.

\section{Computing the equations of Rees algebras by means of a determinantal closure}

In this section we develop some examples that illustrate 
a tentative determinantal procedure to expand sets of 
equations of the Rees algebra of an ideal $I$. This procedure 
arises as a natural extension of a determinantal method that 
computes the maximal minors of the Jacobian dual of a free 
presentation of $I$. A quick overview of the basic settings 
and description of the procedure is 
presented (see also e.g. \cite{suv}, \cite{morey} and 
\cite{vasconcelos1}). Then a brief explanation of the 
generalisation of the procedure is given.

The aim is to understand whether the equations of 
$\rees(I)$ may arise as an ideal of minors of a matrix 
$B(\psi)$ which is iteratively obtained, starting from the 
first syzygies of $I$ and leading to a sort of determinantal closure. 
Such a procedure has its roots in the computational insights 
already stressed in the literature relating the 
equations of Rees algebras and the first syzygies
(see for instance \cite{chw}, \cite{cox} and \cite{dac}; see also 
\ref{rees-elimination}). 

Let us begin with an introductory example:

\begin{Example}{\rm
The following is a particular case of the family considered in 
Example~\ref{exemple-classic}. Let $I=(s^3,t^3,st^2)\subset R=k[s,t]$. 
Let $\varphi:R[T_1,T_2,T_3]\to\rees(I)$ be the polynomial presentation sending
$T_1$ to $s^3$, $T_2$ to $t^3$ and $T_3$ to $st^2$ (viewing $s^3$, $t^3$ and $st^2$ 
in $\rees(I)_1$), and let $Q=\ker\varphi$. 
We get the following Hilbert-Burch presentation
$$R(-4) \oplus R(-5) \stackrel{\scriptsize\left[\begin{array}{cc}   0 & t^2  \\
																										    s & 0    \\
																											 -t & -s^2 \end{array}\right]}{\longrightarrow} R(-3)^3 
										 \stackrel{\scriptsize\left[\begin{array}{ccc}   s^3 & t^3 & st^2 
										 								 \end{array}\right]}{\longrightarrow} I\longrightarrow 0.$$ 
Let the linear forms $S_1,S_2\in R[T_1,T_2,T_3]$ be given by the following identity:
$$[S_1,S_2]=[T_1,T_2,T_3]\cdot\left[\begin{array}{cc}   0 & t^2  \\
																										    s & 0    \\
																											 -t & -s^2 \end{array}\right].$$
In fact, $S_1,S_2$ is a minimal generating set of the module of first syzygies of $I$ using the
obvious correspondence, i.e., $Q\langle 1\rangle =(S_1,S_2)$. 
Observe that we have another matrix identity to describe $S_1,S_2$:
$$[S_1,S_2]=[s,t]\cdot\left[\begin{array}{cc}  T_2 & -sT_3 \\
														-T_3 &  tT_1 \end{array}\right].$$
The determinant $tT_1T_2-sT_3^2$ of the matrix of the right hand part is 
an element of $Q$ which is not in $Q\langle 1\rangle$. Now we iterate the 
process by adding a column corresponding to the new element:
\begin{equation*}
[S_1, S_2, tT_1T_2-sT_3^2]=[s,t]\cdot\left[\begin{array}{ccc} T_2 & -sT_3 & -T_3^2 \\
																							 -T_3 &  tT_1 & T_1T_2 \end{array}\right].
\end{equation*} 
Among the ($2\times 2$)-minors of the matrix of the right hand part, all of them belonging to $Q$, 
$T_1T_2^2-T_3^3$ is the only one that is not contained in the ideal $(S_1, S_2, tT_1T_2-sT_3^2)$. 
Since $T_1T_2^2-T_3^3\notin (s,t)R[T_1,T_2,T_3]$, we can not play the same game again 
and stop. 

So far we have obtained four elements of $Q$: $S_1$, $S_2$ of degree $1$ and two more 
homogeneous elements of $Q$ of degrees $2$ and $3$: $tT_1T_2-sT_3^2, T_1T_2^2-T_3^3$. 
By Example~\ref{exemple-classic}, we already know that they are enough to generate $Q$.\qed
}\end{Example}

\medskip

For the rest of the section, consider the following framework: let $R$ be the 
polynomial ring $k[Z_1,\ldots,Z_r]$ over a field $k$; let $I=(x_1,\ldots,x_s)$ 
be an ideal of $R$ with $R$-free presentation:
$$R^m\stackrel{\psi}{\longrightarrow} R^s\stackrel{\widetilde{\varphi}_1}{\longrightarrow} I\to 0;$$
let $\varphi:R[T_1,\ldots,T_s]\to \rees(I)$ be the polynomial presentation of $\rees(I)$ sending $T_i$
to $x_it$, and let $Q=\ker\varphi$.

\medskip

If $M$ stands for the matrix of $\psi$, $Q\langle 1\rangle$ is generated by the linear forms 
$S_1,\ldots,S_m$ defined as $[S_1,\ldots,S_m]=[T_1,\ldots,T_s]\cdot M$. Suppose that the ideal 
$L$ generated by the entries of $M$ is included in $(Z_1,\ldots,Z_r)$. Suppose also that 
$m=\mu(Q\langle 1\rangle)\geq r$. Let $B(\psi)$ be an $(r\times m)$-matrix 
with entries in $R[T_1,\ldots,T_s]$ satisfying the equation
\begin{equation*}\label{jacobian-dual-definition}
[S_1,\ldots,S_m]=[Z_1,\ldots,Z_r]\cdot B(\psi).
\end{equation*}

\begin{Remark}\label{minors-in-Q}{\rm
If $\Delta$ is an $(r\times r)$-minor of $B(\psi)$, then $\Delta\in Q$. Consequently, the  
ideal $I_r(B(\psi))$ generated by the $(r\times r)$-minors of $M$ is contained in $Q$.  
}\end{Remark}
\demo
Using Cramer's rule we get $\Delta\cdot Z_i\in (S_1,\ldots,S_m)=Q\langle 1\rangle\subset Q$, for
$i=1,\ldots,r$. Within our assumptions $\rees(I)$ is a domain, thus $Q$ is a prime ideal and 
we have $\Delta\in Q$ since $Z_i\notin Q$.
\qed

\medskip

Unless the entries of $M$ are linear in $Z_1,\ldots,Z_r$, 
there will be several choices for $B(\psi)$. By Remark \ref{minors-in-Q}, 
it is clear that $(S_1,\ldots,S_m)+I_r(B(\psi))\subseteq Q$. 
When the equality holds, the ideal $I$ is said to have the \emph{expected equations}. 
For a taste on such classes of ideals see for instance \cite{vasconcelos1}.

\begin{Example}{\rm
The following is a particular case of the family considered in 
Example~\ref{exemple-pseudo-classic}. Let $I=(s^5,t^5,s^2t^3)\subset R=k[s,t]$. 
We get the following Hilbert-Burch presentation: 
$$R(-7)\oplus R(-8)\stackrel{\scriptsize\left[\begin{array}{cc}  0  & t^3   \\
												 	  													 s^2 &  0    \\
												 															-t^2 & -s^3  \end{array}\right]}{\longrightarrow} R(-5)^3 
										\stackrel{\scriptsize \left[\begin{array}{ccc}  s^5  & t^5 & s^2t^3  \end{array}\right]}{\longrightarrow} I\to 0.$$ 
Let $S_1,S_2\in R[T_1,T_2,T_3]$ be the linear forms given by the following identity:
$$[S_1,S_2]=[T_1,T_2,T_3]\cdot \left[\begin{array}{cc}  0  & t^3   \\
												 	  													 s^2 &  0    \\
												 															-t^2 & -s^3  \end{array}\right].$$
In fact, $S_1,S_2$ is a minimal generating set of $Q\langle 1\rangle$.
Observe that we have another matrix identity to describe $S_1,S_2$ in terms of a row vector $[s,t]$ 
and a matrix of the form $B(\psi)$:
$$[S_1,S_2]=[s,t]\cdot\left[\begin{array}{cc}  sT_2 & -s^2T_3 \\
														  -tT_3 & t^2T_1  \end{array}\right].$$
Consider the determinant $st^2T_1T_2-s^2tT_3^2$ of the matrix of the right hand part. Decompose it into its prime
components and clear the factors not belonging to $Q$. We get a new element $tT_1T_2-sT_3^2$ of $Q$ 
not belonging to $Q\langle 1\rangle$. Since $tT_1T_2-sT_3^2\in (s,t)R[T_1,T_2,T_3]$, we add a new 
column:

\begin{equation*}
[S_1,S_2,tT_1T_2-sT_3^2]=[s,t]\cdot\left[\begin{array}{ccc} sT_2 & -s^2T_3 & -T_3^2 \\
																							 -tT_3 &  t^2T_1 & T_1T_2 \end{array}\right],
\end{equation*} 
thus, after clearing the $(2\times 2)$-minors of $B_2$ that already belong to $P_1=(S_1,S_2,tT_1T_2-sT_3^2)$ 
and the irreducible factors not belonging to $Q$, we get $sT_1T_2^2-tT_3^3$. 
Since $sT_1T_2^2-tT_3^3\in (s,t)R[T_1,T_2,T_3]$, we add a new column:
\begin{equation*}
[S_1,S_2,tT_1T_2-sT_3^2,sT_1T_2^2-tT_3^3]=[s,t]\cdot\left[\begin{array}{cccc} sT_2 & -s^2T_3 & -T_3^2 & T_1T_2^2 \\
																							  -tT_3 &  t^2T_1 & T_1T_2 & -T_3^3   \end{array}\right],
\end{equation*}
and again after clearing the $(2\times 2)$-minors of $B_3$ that already belong to 
$P_2=(S_1,S_2,tT_1T_2-sT_3^2,sT_1T_2^2-tT_3^3)$, we get the irreducible form $T_3^5-T_1^2T_2^3$.
Since $T_3^5-T_1^2T_2^3\notin (s,t)R[T_1,T_2,T_3]$ the procedure halts. Let 
$P_3=(s^2T_2-t^2T_3,tT_1T_2-sT_3^2,sT_1T_2^2-tT_3^3,T_3^5-T_1^2T_2^3)\subset Q$: 
by Example~\ref{exemple-pseudo-classic} we know that $P_3=Q$.\qed
}\end{Example}

\medskip

The next algorithm encodes a sequence of instructions corresponsing to the steps
completed in the previous discussion and examples. It is a tentative procedure 
to expand a generating set of $Q\langle 1\rangle$ into a generating set of $Q$. 
The reader will observe that the same process could be used to expand any finite 
set of elements in $Q$ and could lead to the definition of a sort of determinantal 
closure for subideals of $Q$.

\begin{Algorithm}\label{algorithm}{\rm
In what follows, $L$ and $K$ will be lists of forms in $R[T_1,\ldots,T_s]$, with all 
the forms in $L$ belonging to $(Z_1,\ldots,Z_r)R[T_1,\ldots,T_s]$, $P$ will stand for 
an ideal in $R[T_1,\ldots,T_s]$ and $B$ will denote a matrix of forms in 
$R[T_1,\ldots,T_s]$.

\begin{algorithmic}
\STATE Compute $S=\{S_{1},\ldots,S_{k}\}$ a minimal generating set of $Q_1$;
\STATE Set $L\leftarrow S$;
\STATE Let $P=(L)$, the ideal generated by the elements of the list $L$;
\STATE Since all the forms of the list $L$ belong to $(Z_1,\ldots,Z_r)R[T_1,\ldots,T_s]$, 
there is an $(r\times k)$-matrix $B_1$ such that $[L]=[Z_1,\ldots,Z_r]\cdot B_1$;
\STATE Set $B\leftarrow B_1$;
\WHILE{$K\neq\emptyset$}
\STATE Compute the $(r\times r)$-minors of $B$ and put them into a list $K$; 
\STATE Set $K\leftarrow K \backslash \{ m\in K\;|\; m\in P\}$, i.e., discard all the 
elements of that belong to $P$;
\FORALL{$f\in K$}
\STATE Compute all the irreducible components $\{f_{i}\}$ of $f$ in $k[Z_1,\ldots,Z_r][T_1,\ldots,T_s]$; 
\STATE Set $K\leftarrow (K\backslash f)\cup \{f_i\;|\; f_i\in Q\}$;
\ENDFOR
\STATE Set $L\leftarrow L\cup K$;
\STATE Set $P=(L)$ the ideal generated by the elements of $P$;
\STATE Set $K\leftarrow K\backslash \{f\in K\;|\; f\notin (Z_1,\ldots,Z_r)R[T_1,\ldots,T_s]\}$;
\STATE Since all the forms of the list $K$ belong to $(Z_1,\ldots,Z_r)R[T_1,\ldots,T_s]$, 
there is a matrix $B_K$ with $r$ rows such that $[K]=[Z_1,\ldots,Z_r]\cdot B_K$;
\STATE Set $B\leftarrow [B\;|\; B_K]$;
\ENDWHILE
\RETURN $P$.\qed

\end{algorithmic}
}\end{Algorithm}

\medskip

We proceed to run Algorithm \ref{algorithm} in one more example, respecting the notations used in the
pseudo-code.

\begin{Example}{\rm
Let $I=(u,z)\cap (v,t)=(uv,vz,zt,ut)\subset k[u,v,z,t]$ be the ideal
of the graph cycle $C_4$ (see also Example~\ref{exemple-C4}).
In this case we have the following a free presentation:
$$R(-3)^4\stackrel{\scriptsize \left[\begin{array}{cccc} z & t & 0 & 0 \\
																												-u & 0 & t & 0 \\
																												 0 & 0 & -v & u \\
																												 0 & -v & 0 & -z \end{array}\right]}{\longrightarrow} R(-2)^4
					\stackrel{\scriptsize\left[\begin{array}{cccc} uv & vz & zt & ut \end{array}\right]}{\longrightarrow} I\to 0.$$
Let $P_0=Q\langle 1\rangle=(S_1,S_2,S_3,S_4)$, where
$$[S_1,S_2,S_3,S_4]=[T_1,T_2,T_3,T_4]\left[\begin{array}{cccc} z & t & 0 & 0 \\
																															-u & 0 & t & 0 \\
																															 0 & 0 & -v & u \\
																															 0 & -v & 0 & -z \end{array}\right].$$
Writing $[S_1,S_2,S_3,S_4]$ as $[u,v,z,t]\cdot B_1$, we have
$$B_1=\left[\begin{array}{cccc} -T_2 &  0   &   0  &  T_3  \\
														   	 0 & -T_4 & -T_3 &  0 \\
															 T_1 &  0   &   0  & -T_4  \\
																 0 & T_1  &  T_2 &  0 \end{array}\right].$$														
The only $(4\times 4)$-determinant of $B_1$ is $(T_1T_4-T_2T_3)^2 \notin P_0$ 
and its only irreducible factor $T_1T_4-T_2T_3$ belongs to $Q$, then $K_1=\{T_1T_4-T_2T_3\}$
and we set $P_1=P_0+(T_1T_4-T_2T_3)$. Since $T_1T_4-T_2T_3 \notin (u,v,z,t)R[T_1,T_2,T_3,T_4]$ 
the procedure halts and returns $$P_1=(zT_1-uT_2,tT_1-vT_4,tT_2-vT_3,uT_3-zT_4,T_1T_4-T_2T_3).$$ 
It is readily seen from the results in \cite{villarreal} that $P_1$ is the ideal of 
equations of $\rees(I)$, i.e., $P_1=Q$.\qed
}\end{Example}







\end{document}